\documentclass[12pt]{article}
\usepackage{latexsym, amssymb}
\usepackage{dsfont}

\DeclareMathAlphabet\gothic{U}{euf}{m}{n}

\makeatletter
\def\eqnarray{\stepcounter{equation}\let\@currentlabel=\theequation
\global\@eqnswtrue
\tabskip\@centering\let\\=\@eqncr
$$\halign to \displaywidth\bgroup\hfil\global\@eqcnt\z@
  $\displaystyle\tabskip\z@{##}$&\global\@eqcnt\@ne
  \hfil$\displaystyle{{}##{}}$\hfil
  &\global\@eqcnt\tw@ $\displaystyle{##}$\hfil
  \tabskip\@centering&\llap{##}\tabskip\z@\cr}

\def\endeqnarray{\@@eqncr\egroup
      \global\advance\c@equation\m@ne$$\global\@ignoretrue}

\def\@yeqncr{\@ifnextchar [{\@xeqncr}{\@xeqncr[5pt]}}
\makeatother

\begin{document}
\bibliographystyle{tom}

\newtheorem{lemma}{Lemma}[section]
\newtheorem{thm}[lemma]{Theorem}
\newtheorem{cor}[lemma]{Corollary}
\newtheorem{voorb}[lemma]{Example}
\newtheorem{rem}[lemma]{Remark}
\newtheorem{prop}[lemma]{Proposition}
\newtheorem{ddefinition}[lemma]{Definition}
\newtheorem{stat}[lemma]{{\hspace{-5pt}}}

\newenvironment{remarkn}{\begin{rem} \rm}{\end{rem}}
\newenvironment{exam}{\begin{voorb} \rm}{\end{voorb}}
\newenvironment{definition}{\begin{ddefinition} \rm}{\end{ddefinition}}

\newcommand{\gota}{\gothic{a}}
\newcommand{\gotb}{\gothic{b}}
\newcommand{\gotc}{\gothic{c}}
\newcommand{\gote}{\gothic{e}}
\newcommand{\gotf}{\gothic{f}}
\newcommand{\gotg}{\gothic{g}}
\newcommand{\gothh}{\gothic{h}}
\newcommand{\gotk}{\gothic{k}}
\newcommand{\gotm}{\gothic{m}}
\newcommand{\gotn}{\gothic{n}}
\newcommand{\gotp}{\gothic{p}}
\newcommand{\gotq}{\gothic{q}}
\newcommand{\gotr}{\gothic{r}}
\newcommand{\gots}{\gothic{s}}
\newcommand{\gotu}{\gothic{u}}
\newcommand{\gotv}{\gothic{v}}
\newcommand{\gotw}{\gothic{w}}
\newcommand{\gotz}{\gothic{z}}
\newcommand{\gotA}{\gothic{A}}
\newcommand{\gotB}{\gothic{B}}
\newcommand{\gotG}{\gothic{G}}
\newcommand{\gotL}{\gothic{L}}
\newcommand{\gotS}{\gothic{S}}
\newcommand{\gotT}{\gothic{T}}

\newcounter{teller}
\renewcommand{\theteller}{(\alph{teller})}
\newenvironment{tabel}{\begin{list}%
{\rm  (\alph{teller})\hfill}{\usecounter{teller} \leftmargin=1.1cm
\labelwidth=1.1cm \labelsep=0cm \parsep=0cm}
                      }{\end{list}}

\newcounter{tellerr}
\renewcommand{\thetellerr}{(\roman{tellerr})}
\newenvironment{tabeleq}{\begin{list}%
{\rm  (\roman{tellerr})\hfill}{\usecounter{tellerr} \leftmargin=1.1cm
\labelwidth=1.1cm \labelsep=0cm \parsep=0cm}
                         }{\end{list}}

\newcounter{tellerrr}
\renewcommand{\thetellerrr}{(\Roman{tellerrr})}
\newenvironment{tabelR}{\begin{list}%
{\rm  (\Roman{tellerrr})\hfill}{\usecounter{tellerrr} \leftmargin=1.1cm
\labelwidth=1.1cm \labelsep=0cm \parsep=0cm}
                         }{\end{list}}

\newcounter{proofstep}
\newcommand{\nextstep}{\refstepcounter{proofstep}\vertspace \par 
          \noindent{\bf Step \theproofstep} \hspace{5pt}}
\newcommand{\firststep}{\setcounter{proofstep}{0}\nextstep}

\newcommand{\Ni}{\mathds{N}}
\newcommand{\Qi}{\mathds{Q}}
\newcommand{\Ri}{\mathds{R}}
\newcommand{\Ci}{\mathds{C}}
\newcommand{\Ti}{\mathds{T}}
\newcommand{\Zi}{\mathds{Z}}
\newcommand{\Fi}{\mathds{F}}

\newcommand{\proof}{\mbox{\bf Proof} \hspace{5pt}} 
\newcommand{\remark}{\mbox{\bf Remark} \hspace{5pt}}
\newcommand{\vertspace}{\vskip10.0pt plus 4.0pt minus 6.0pt}

\newcommand{\simh}{{\stackrel{{\rm cap}}{\sim}}}
\newcommand{\ad}{{\mathop{\rm ad}}}
\newcommand{\Ad}{{\mathop{\rm Ad}}}
\newcommand{\Aut}{\mathop{\rm Aut}}
\newcommand{\arccot}{\mathop{\rm arccot}}
\newcommand{\capp}{{\mathop{\rm cap}}}
\newcommand{\rcapp}{{\mathop{\rm rcap}}}
\newcommand{\diam}{\mathop{\rm diam}}
\newcommand{\divv}{\mathop{\rm div}}
\newcommand{\codim}{\mathop{\rm codim}}
\newcommand{\RRe}{\mathop{\rm Re}}
\newcommand{\IIm}{\mathop{\rm Im}}
\newcommand{\Tr}{{\mathop{\rm Tr \,}}}
\newcommand{\Vol}{{\mathop{\rm Vol}}}
\newcommand{\card}{{\mathop{\rm card}}}
\newcommand{\supp}{\mathop{\rm supp}}
\newcommand{\sgn}{\mathop{\rm sgn}}
\newcommand{\essinf}{\mathop{\rm ess\,inf}}
\newcommand{\esssup}{\mathop{\rm ess\,sup}}
\newcommand{\Int}{\mathop{\rm Int}}
\newcommand{\lcm}{\mathop{\rm lcm}}
\newcommand{\loc}{{\rm loc}}

\newcommand{\mod}{\mathop{\rm mod}}
\newcommand{\spann}{\mathop{\rm span}}
\newcommand{\one}{\mathds{1}}

\hyphenation{groups}
\hyphenation{unitary}

\newcommand{\tfrac}[2]{{\textstyle \frac{#1}{#2}}}

\newcommand{\ca}{{\cal A}}
\newcommand{\cb}{{\cal B}}
\newcommand{\cc}{{\cal C}}
\newcommand{\cd}{{\cal D}}
\newcommand{\ce}{{\cal E}}
\newcommand{\cf}{{\cal F}}
\newcommand{\ch}{{\cal H}}
\newcommand{\ci}{{\cal I}}
\newcommand{\ck}{{\cal K}}
\newcommand{\cl}{{\cal L}}
\newcommand{\cm}{{\cal M}}
\newcommand{\co}{{\cal O}}
\newcommand{\cs}{{\cal S}}
\newcommand{\ct}{{\cal T}}
\newcommand{\cx}{{\cal X}}
\newcommand{\cy}{{\cal Y}}
\newcommand{\cz}{{\cal Z}}

\newlength{\hightcharacter}
\newlength{\widthcharacter}
\newcommand{\covsup}[1]{\settowidth{\widthcharacter}{$#1$}\addtolength{\widthcharacter}{-0.15em}\settoheight{\hightcharacter}{$#1$}\addtolength{\hightcharacter}{0.1ex}#1\raisebox{\hightcharacter}[0pt][0pt]{\makebox[0pt]{\hspace{-\widthcharacter}$\scriptstyle\circ$}}}
\newcommand{\cov}[1]{\settowidth{\widthcharacter}{$#1$}\addtolength{\widthcharacter}{-0.15em}\settoheight{\hightcharacter}{$#1$}\addtolength{\hightcharacter}{0.1ex}#1\raisebox{\hightcharacter}{\makebox[0pt]{\hspace{-\widthcharacter}$\scriptstyle\circ$}}}
\newcommand{\scov}[1]{\settowidth{\widthcharacter}{$#1$}\addtolength{\widthcharacter}{-0.15em}\settoheight{\hightcharacter}{$#1$}\addtolength{\hightcharacter}{0.1ex}#1\raisebox{0.7\hightcharacter}{\makebox[0pt]{\hspace{-\widthcharacter}$\scriptstyle\circ$}}}

\thispagestyle{empty}

\vspace*{1cm}
\begin{center}
{\Large\bf Partial Gaussian bounds \\[2mm]
for  
degenerate differential operators II} \\[5mm]
\large A.F.M. ter Elst$^1$ and E.M.  Ouhabaz$^2$

\end{center}

\vspace{5mm}

\begin{center}
{\bf Abstract}
\end{center}

\begin{list}{}{\leftmargin=1.8cm \rightmargin=1.8cm \listparindent=10mm 
   \parsep=0pt}
\item
Let $A = - \sum \partial_k \, c_{kl} \, \partial_l$ be a degenerate 
sectorial differential operator with complex bounded mesaurable coefficients.
Let $\Omega \subset \Ri^d$ be open and suppose that $A$ is strongly elliptic 
on $\Omega$.
Further, let $\chi \in C_{\rm b}^\infty(\Ri^d)$ be such that an 
$\varepsilon$-neighbourhood of $\supp \chi$ is contained in $\Omega$.
Let $\nu \in (0,1]$ and suppose that the ${c_{kl}}_{|\Omega} \in C^{0,\nu}(\Omega)$.
Then we prove (H\"older) Gaussian kernel bounds for the kernel of 
the operator $u \mapsto \chi \, S_t (\chi \, u)$, where
$S$ is the semigroup generated by $-A$.
Moreover, if $\nu = 1$ and the coefficients are real, then we prove 
Gaussian bounds for the kernel of the operator $u \mapsto \chi \, S_t u$
and for the derivatives in the first variable.
Finally we show boundedness on $L_p(\Ri^d)$ of various Riesz transforms.

\end{list}

\vspace{4cm}
\noindent
January 2012

\vspace{5mm}

\vspace{5mm}

\vspace{15mm}

\noindent
{\bf Home institutions:}    \\[3mm]
\begin{tabular}{@{}cl@{\hspace{10mm}}cl}
1. & Department of Mathematics  & 
  2. & Institut de Math\'ematiques de Bordeaux \\
& University of Auckland   & 
  & Universit\'e Bordeaux 1, UMR 5251,  \\
& Private bag 92019 & 
  &351, Cours de la Lib\'eration  \\
& Auckland 1142 & 
  &  33405 Talence \\
& New Zealand  & 
  & France  \\[8mm]
\end{tabular}

\newpage
\setcounter{page}{1}

\section{Introduction} \label{Slocseh1}

If $A$ is a strongly elliptic second-order operator on $\Ri^d$ in 
divergence form with complex bounded H\"older continuous coefficients,
then it is well known that it generates a holomorphic semigroup $S$ which 
satisfies Gaussian kernel bounds and Gaussian bounds for first 
order derivatives in each of the variables.
If $A$ is merely partially strongly elliptic on an open set $\Omega \subset \Ri^d$
then in general Gaussian bounds on $\Ri^d$ fail, but in a previous
paper \cite{EO1} we showed Gaussian kernel bounds on good parts of $\Omega$
if the coefficients of $A$ are real and measurable.
Precisely, if $\chi \in C_{\rm b}^\infty(\Omega,\Ri)$ and if $A$ is strongly
elliptic on $\supp \chi$, then for all $t > 0$ the operator 
$M_\chi \, S_t \, M_\chi$ has a H\"older continuous kernel
satisfying (H\"older) Gaussian bounds, where $M_\chi$ is the multiplication 
operator with the function~$\chi$.
In this paper we extend this to (H\"older) derivatives of the kernel 
if the coefficients of the operator $A$ are complex H\"older continuous
on $\Omega$ and the distance $d(\supp \chi, \Omega^{\rm c}) > 0$,
that is an $\varepsilon$-neighbourhood of $\supp \chi$ is still in~$\Omega$.
If in addition the coefficients are in $W^{1,\infty}(\Omega)$
and real on $\Ri^d$, then we also show that for all $t > 0$ the operator
$M_\chi \, S_t$ has a kernel $K_t$ satisfying Gaussian bounds.
This is remarkable, since there is no cut-off for the operator 
$M_\chi \, S_t$ on the right.
Moreover, we show that there exists a representative of the kernel $K_t$
such that $(t,x,y) \mapsto K_t(x,y)$ is measurable on $(0,\infty) \times \Ri^d \times \Ri^d$
and $x \mapsto K_t(x,y)$ is once differentiable for all $y \in \Ri^d$ and $t > 0$,
and the derivatives satisfy (H\"older) Gaussian bounds.
This allows to prove boundedness of the Riesz transforms
$\nabla \, M_\chi \, (I + A)^{-1/2}$ on $L_p(\Ri^d)$
for all $p \in (1,\infty)$.

Throughout this paper the field is $\Ci$.
Fix $d \in \Ni$ and for all $k,l \in \{ 1,\ldots,d \} $ let $c_{kl} \colon \Ri^d \to \Ci$
be a measurable bounded function.
Suppose that the matrix $C(x) := (c_{kl}(x))$ is uniformly sectorial for all $x \in \Ri^d$,
i.e., there exists a $\theta \in [0,\frac{\pi}{2})$ such that 
\[
\sum_{k,l = 1}^d c_{kl}(x) \, \xi_k \, \overline{\xi_l} \in \Sigma_\theta
\]
for all $(\xi_1,\dots,  \xi_d) \in \Ci^d$ and  $x \in \Ri^d$, where
\[
\Sigma_\theta
= \{ r \, e^{i \alpha} : r \geq 0 \mbox{ and } \alpha \in [-\theta,\theta] \} 
.  \]
Define the form 
$\gota \colon W^{1,2}(\Ri^d) \times W^{1,2}(\Ri^d) \to \Ci$ by
\[
\gota(u,v)
= \sum_{k,l=1}^d \int_{\Ri^d} c_{kl} \, (\partial_k u) \, \overline{\partial_l v}
.  \]
Then $\gota$ is a densely defined sectorial form.
In general $\gota$ is not closable, but nevertheless one can 
assign a semigroup generator $A$ with~$\gota$ as follows.
If $u,f \in L_2(\Ri^d)$ then $u \in D(A)$ and $A u = f$
if and only if there exist $u_1,u_2,\ldots \in W^{1,2}(\Ri^d)$ 
such that $\lim u_n = u$ in $L_2(\Ri^d)$, $\sup \RRe \gota(u_n) < \infty$ and 
$\lim \gota(u_n,v) = (f,v)$ for all $v \in W^{1,2}(\Ri^d)$.
The operator~$A$ is well defined and is $m$-sectorial by Theorem~1.1 in 
\cite{AE2}.
If $\gota$ is closable then $A$ is the operator associated with the 
closure $\overline \gota$ of the form $\gota$ in the sense of Kato \cite{Kat1}.
We call $A$ the sectorial degenerate differential operator with coefficients $(c_{kl})$.
Formally,  $A = - \sum_{k,l} \partial_l \, c_{kl} \, \partial_k$.
We denote  by $S = (S_t)_{t > 0}$ the contraction semigroup generated by $-A$ on $L_2(\Ri^d)$.
Then $S$ is holomorphic on the sector $\Sigma_{\theta_\gota}^\circ$,
where throughout this paper we define $\theta_\gota = \frac{\pi}{2} - \theta$.
Let $\Omega \subset \Ri^d$ open.
We suppose that the coefficients of $A$ are strongly elliptic on $\Omega$,
that is, there exists a $\mu > 0$ such that 
\[
\RRe \sum_{k,l=1}^d c_{kl}(x) \, \xi_k \, \overline{\xi_l}
\geq \mu \, |\xi|^2
\]
for all $\xi \in \Ci^d$ and a.e.\ $x \in \Omega$.

The main results of this paper are the following.
The first theorem is for complex H\"older continuous coefficients
on $\Omega$, but with a multiplication operator on both sides 
of the semigroup.

\begin{thm} \label{tlocseh101}
Let $A$ be a sectorial degenerate differential operator with coefficients $(c_{kl})$,
where $c_{kl} \colon \Ri^d \to \Ci$ 
is a bounded measurable function for all $k,l \in \{ 1,\ldots,d \} $.
Let $\Omega \subsetneqq \Ri^d$ be open and suppose that $(c_{kl})$ is strongly 
elliptic on $\Omega$.
Let $\nu \in (0,1)$ and suppose that ${c_{kl}}_{|\Omega} \in C^{0,\nu}(\Omega)$ 
for all $k,l \in \{ 1,\ldots,d \} $.
Let $\chi \in C^\infty_{\rm b}(\Ri^d)$  with $\chi \neq 0$ and suppose 
$d(\supp \chi, \Omega^{\rm c}) > 0$.
Then there exists a continuous function $(z,x,y) \mapsto K_z(x,y)$
from $\Sigma_{\theta_\gota}^\circ \times \Ri^d \times \Ri^d$ into $\Ci$ such that 
the following is valid.
\begin{itemize}
\item
The function $K_z$ is the kernel of the operator $M_\chi \, S_z \, M_\chi$
for all $z \in \Sigma_{\theta_\gota}^\circ$, where $S$ is the semigroup generated by $-A$.
\item
The function $K_z$ is once differentiable in each variable and the 
derivative with respect to one variable is differentiable in the 
other variable.
\item
For every multi-index $\alpha,\beta$ with 
$|\alpha|,|\beta| \leq 1$, $\kappa > 0$ and $\tau \in [0,1)$
there exist $a,b > 0$ such that 
\[
|(\partial_x^\alpha \, \partial_y^\beta \, K_z)(x,y)|
\leq a \, |z|^{-d/2} \, |z|^{-(|\alpha| + |\beta|)/2} \, 
      (1 + |z|)^{\frac{d+|\alpha| +|\beta|}{2}} \,
      e^{-b \, \frac{|x-y|^2}{|z|}} 
\]
and 
\begin{eqnarray*}
\lefteqn{
|(\partial_x^\alpha \, \partial_y^\beta \, K_z)(x+h,y+k) 
     - (\partial_x^\alpha \, \partial_y^\beta \, K_z)(x,y)|
} \hspace{20mm} \\*
& \leq & a \, |z|^{-d/2} \, |z|^{-(|\alpha| + |\beta|)/2} \, 
    \left( \frac{|h| + |k|}{|x-y| + \sqrt{|z|}} \right)^\nu
   \, (1 + |z|)^{\frac{d+|\alpha| +|\beta| + \nu}{2}} \, 
      e^{-b \, \frac{|x-y|^2}{|z|}} 
\end{eqnarray*}
for all $z \in \Sigma_{\theta_\gota}^\circ$ and $x,y,h,k \in \Ri^d$ with 
$|h| + |k| \leq \tau \, |x-y| + \kappa \, \sqrt{|z|}$.
\end{itemize}
\end{thm}

The second result is for merely one multiplication operator
on the left of the semigroup, but it requires that the 
coefficients of the operator are real on $\Ri^d$ and 
uniformly Lipschitz on~$\Omega$.

\begin{thm} \label{tlocseh102}
Let $A$ be a sectorial degenerate differential operator with real
coefficients $(c_{kl})$,
where $c_{kl} \colon \Ri^d \to \Ri$ 
is a bounded measurable function for all $k,l \in \{ 1,\ldots,d \} $.
Let $S$ be the semigroup generated by $-A$.
Let $\Omega \subsetneqq \Ri^d$ be open and suppose that $(c_{kl})$ is strongly 
elliptic on $\Omega$.
Suppose that ${c_{kl}}_{|\Omega} \in W^{1,\infty}(\Omega)$ 
for all $k,l \in \{ 1,\ldots,d \} $.
Let $\chi \in C^\infty_{\rm b}(\Ri^d)$ with $\chi \neq 0$ and suppose
$d(\supp \chi, \Omega^{\rm c}) > 0$.
Then there exists a measurable function 
$(t,x,y) \mapsto K_t(x,y)$ from $(0,\infty) \times \Ri^d \times \Ri^d$ into $\Ri$ such that 
the following is valid.
\begin{itemize}
\item
The function $K_t$ is a kernel of $M_\chi \, S_t$ for all $t > 0$.
\item
The function $x \mapsto K_t(x,y)$ is continuously differentiable on $\Ri^d$
for all $t > 0$ and $y \in \Ri^d$.
\item
The function $t \mapsto K_t(x,y)$ is continuous for all $x,y \in \Ri^d$.
\item
For every multi-index $\alpha$ with 
$|\alpha| \leq 1$, $\nu \in (0,1)$, $\varepsilon > 0$, $\kappa > 0$ and $\tau \in [0,1)$
there exist $a,b > 0$ such that 
\[
|(\partial_x^\alpha K_t)(x,y)| 
   \leq a \, t^{-d/2} \, t^{-|\alpha| / 2} \, e^{\varepsilon t} \,
      e^{-b \, \frac{|x-y|^2}{t}} 
\]
and 
\[
|(\partial_x^\alpha \, K_t)(x+h,y) 
     - (\partial_x^\alpha \, K_t)(x,y)|
\leq a \, t^{-d/2} \, t^{-|\alpha|/2} \, 
    \left( \frac{|h|}{|x-y| + \sqrt{t}} \right)^\nu
   \, e^{\varepsilon t} \, 
      e^{-b \, \frac{|x-y|^2}{t}} 
\]
for all $t > 0$ and $x,y,h \in \Ri^d$ with 
$|h| \leq \tau \, |x-y| + \kappa \, \sqrt{t}$.
\end{itemize}
\end{thm}

In Theorem~\ref{tlocseh101} the function $(t,x,y) \mapsto K_t(x,y)$
is continuous, whilst it is not clear whether the 
function $(t,x,y) \mapsto K_t(x,y)$ is continuous 
in the setting of Theorem~\ref{tlocseh102}.
Likely, there even does not exists a
continuous function which is equal to this function almost 
everywhere on $(0,\infty) \times \Ri^d \times \Ri^d$.
On the other hand, we prove measurability jointly in the three variables
and do not work with an equivalent class of functions, for which 
the representative changes all the time.
Since there are uncountable many $y \in \Ri^d$ this complicates the proof.

We also investigate boundedness on $L_p$ of Riesz transform type operators. 
We obtain the following result.

\begin{thm} \label{tlocseh103}
Let $A$ be a 
sectorial degenerate differential operator with complex coefficients $(c_{kl})$.
Let $\Omega \subsetneqq \Ri^d$ be open and suppose that $(c_{kl})$ is strongly 
elliptic on $\Omega$.
Let $\chi \in C^\infty_{\rm b}(\Ri^d)$ with $\chi \neq 0$ and
$d(\supp \chi, \Omega^{\rm c}) > 0$.
Then one has the following.
\begin{tabel}
\item \label{tlocseh103-1}
Let $\nu \in (0,1)$ and suppose that ${c_{kl}}_{|\Omega} \in C^{0,\nu}(\Omega)$ 
for all $k,l \in \{ 1,\ldots,d \} $.
Then the Riesz transforms $\nabla \, M_\chi \, (I + A)^{-1/2} \, M_\chi$
are bounded on $L_p(\Ri^d)$ for all $p \in (1,\infty)$.
\item \label{tlocseh103-2}
Suppose that ${c_{kl}}_{|\Omega} \in W^{1,\infty}(\Omega)$ 
and $c_{kl}$ is real valued for all $k,l \in \{ 1,\ldots,d \} $.
Then the Riesz transforms $\nabla \, M_\chi \, (I + A)^{-1/2}$
are bounded on $L_p(\Ri^d)$ for all $p \in (1,\infty)$.
\item \label{tlocseh103-3}
Let $\nu \in (0,1)$ and suppose that ${c_{kl}}_{|\Omega} \in C^{0,\nu}(\Omega)$ 
for all $k,l \in \{ 1,\ldots,d \} $.
Then the Riesz transforms $\nabla \, M_\chi \, (I + A)^{-1/2}$
are bounded on $L_2(\Ri^d)$.
\end{tabel}
\end{thm}

Using Morrey and Campanato spaces we prove Theorem~\ref{tlocseh101} as
in \cite{EO1}, if the operator is strongly elliptic on $\Ri^d$ 
and the coefficients are H\"older continuous on $\Ri^d$.
We carefully control all the constants and show that they depend
only on the ellipticity constant on $\Omega$ and on the
H\"older continuity of the coefficients on $\Omega$.
Then the Arz\'ela--Ascoli theorem together with two approximations give the 
theorem.
We prove a quantitive version of Theorem~\ref{tlocseh101} in Section~\ref{Slocseh2}.

Since
\[
M_\chi^2 \, S_t
= M_\chi \, S_t \, M_\chi + M_\chi \, [M_\chi,S_t]
\]
and one can use Theorem~\ref{tlocseh101} to handle the first term,
it suffices to obtain good estimates on the commutator to derive the bounds
of Theorem~\ref{tlocseh102}.
This is done in Section~\ref{Slocseh3}.

Finally, in Section~\ref{Slocseh4} we prove the boundedness of the 
Riesz transforms of Theorem~\ref{tlocseh103} and the boundedness of several 
other Riesz transforms.
For strongly elliptic operators in divergence form with complex bounded
measurable coefficients the boundedness of the Riesz transforms on $L_2(\Ri^d)$
was the longstanding open Kato problem until it was solved by 
Auscher--Hofmann--Lacey--McIntosh--Tchamitchian \cite{AHLMT}.
For H\"older continuous coefficients the Kato problem was solved 
earlier by McIntosh \cite{McI3} and a simplified proof was given in \cite{ER16}.
In the proof of Theorem~\ref{tlocseh103}\ref{tlocseh103-1} we adapt
this simplified proof.
In \cite{EO3} Theorem~1.2 we proved Theorem~\ref{tlocseh103}\ref{tlocseh103-1}
for merely measurable coefficients, but with the restriction that the 
coefficients are real symmetric and $p \in (1,2]$.

In the proofs we need various times to transfer semigroup estimates
into Gaussian bounds, with control of large time behaviour,
using the Davies perturbation method. 
Note that we deduce polynomial growth for large time for the kernel 
bounds in Theorem~\ref{tlocseh101}.
The techniques are more or less folklore, however scattered over
the literature. 
In the appendix we collect them together for the convenience 
of the reader.
Finally, by decomposing $\chi$ into its real and imaginary part, 
for simplicity we may and do assume throughout the rest of this paper that 
the various cut-off functions $\chi$, $\tilde \chi$,\ldots\  are 
always real valued.

\section{Complex H\"older continuous coefficients} \label{Slocseh2}

We start this section with the definition of a number of classes of 
coefficients and operators. 
The main aim is to obtain results for elements of these classes and 
that the constants involved are uniformly for a given class.

Let $\theta \in [0,\frac{\pi}{2})$ and $M > 0$.
Define $\cs(\theta,M)$
to be the set of all measurable $C \colon \Ri^d \to \Ci^{d \times d}$ such that 
\[
\begin{array}{ll}
\langle C(x) \, \xi, \xi \rangle \in \Sigma_\theta
   & \mbox{for all } x \in \Ri^d \mbox{ and } \xi \in \Ci^d  \mbox{, and,} \\[5pt]
\|C(x)\| \leq M
   & \mbox{for all } x \in \Ri^d , \\[5pt]
\end{array}
\]
where $\|C(x)\|$ is 
the $\ell_2$-norm of $C(x)$ in $\Ci^d$ and $\langle \cdot , \cdot \rangle$ is the 
inner product on $\Ci^d$.
For all $C \in \cs(\theta,M)$ define the sectorial form
form $\gota_C \colon W^{1,2}(\Ri^d) \times W^{1,2}(\Ri^d) \to \Ci$ by
\[
\gota_C(u,v)
= \int_{\Ri^d} \sum_{k,l=1}^d c_{kl} \, (\partial_k u) \, \overline{(\partial_l v)}
\]
and let $A_C$ and $S^C$ be the associated operator and semigroup.
Here and in the sequel $c_{kl}(x)$ is the appropriate
matrix coefficient of $C(x)$.
If no confusion is possible then we drop the $C$ and write $\gota = \gota_C$, 
$A = A_C$ and $S = S^C$.
For all $C \in \cs(\theta,M)$ define $\Re C \colon \Ri^d \to \Ci$
by 
\[
(\Re C)(x) = \tfrac{1}{2} \Big( C(x) + C(x)^* \Big)
.  \]
Then $\Re C \in \cs(0,M)$ and $A_{\Re C}$ is self-adjoint.
Moreover, $\gota_{\Re C}(u) = \RRe \gota(u)$ for all $u \in W^{1,2}(\Ri^d)$.
Next, let $\ce(\theta,M)$ be the set of all 
$C \in \cs(\theta,M)$
such that there exists a $\mu_0 > 0$ such that 
$\RRe \langle C(x) \, \xi, \xi \rangle \geq \mu_0 \, |\xi|^2$ 
for all $x \in \Ri^d$ and $\xi \in \Ci^d$.
We emphasise that the constant $\mu_0$ depends on $C$.

Let $Y \subset \Ri^d$ be a set, $\theta \in [0,\frac{\pi}{2})$ and
$\mu,M > 0$.
Let $\cs(Y,\theta,\mu,M)$
be the set of all $C \in \cs(\theta,M)$ such that 
\[
\RRe \langle C(x) \, \xi, \xi \rangle \geq \mu \, |\xi|^2
   \mbox{ for all } x \in Y \mbox{ and } \xi \in \Ci^d 
\]
and define 
\[
\ce(Y,\theta,\mu,M)
= \cs(Y,\theta,\mu,M) \cap \ce(\theta,M)
.  \]
Next, let $\nu \in (0,1]$ and suppose that $Y$ contains at least two 
elements $x,y$ with $0 < |x-y| \leq 1$.
The space $C^{0,\nu}(Y)$ is the space of all H\"older continuous functions on $Y$ 
with seminorm
\[
|||u|||_{C^{0,\nu}(Y)}
= \sup \{ \frac{|u(x) - u(y)|}{|x-y|^\nu} : x,y \in Y, \; 0 < |x-y| \leq 1 \}
.  \]
Let $\cs^\nu(Y,\theta,\mu,M)$
be the set of all $C \in \cs(Y,\theta,\mu,M)$ such that 
\[
|||{c_{kl}}_{|Y}|||_{C^{0,\nu}(Y)} \leq M
   \mbox{ for all } k,l \in \{ 1,\ldots,d \} 
,  \]
and define
\[
\ce^\nu(Y,\theta,\mu,M)
= \cs^\nu(Y,\theta,\mu,M) \cap \ce(\theta,M)
.  \]
Finally, let $\ce\ch^\nu(Y,\theta,\mu,M)$ be the set of all 
$C \in \ce^\nu(Y,\theta,\mu,M)$
such that $|||c_{kl}|||_{C^{0,\nu}(\Ri^d)} < \infty$ 
for all $k,l \in \{ 1,\ldots,d \} $.
If $C \in \cs^\nu(Y,\theta,\mu,M)$ then $A_C$ is sectorial on $L_2(\Ri^d)$,
whilst $A_C$ is strongly elliptic on $\Ri^d$ if $C \in \ce^\nu(Y,\theta,\mu,M)$.
Finally, $A_C$ is strongly elliptic with H\"older continuous coefficients
on $\Ri^d$ if $C \in \ce\ch^\nu(Y,\theta,\mu,M)$.
In any case, $A_C$ is strongly elliptic on the set $Y$ with ellipticity
constant at least $\mu$.

In the proof of the theorems we frequently need the Davies perturbation.
For all $\rho \in \Ri$ and $\psi \in W^{1,\infty}(\Ri^d)$
define the multiplication
operator $U_\rho$ by $U_\rho u = e^{-\rho \psi} u$.
For all $n \in \Ni$ let
\[
\cd_n
= \{ \psi \in W^{n,\infty}(\Ri^d,\Ri) : 
          \|\sum_{1 \leq |\alpha| \leq n} |\partial^\alpha \psi|^2 \: \|_\infty \leq 1 \}
.  \]
Thus
\[
\cd_1
= \{ \psi \in W^{1,\infty}(\Ri^d,\Ri) : \|\nabla \psi\|_\infty \leq 1 \}
.  \]
Let $C \in \cs(\theta,M)$.
Define $S^{(C,\rho)}_t = S^{(\rho)}_t = U_\rho S_t U_{-\rho}$ to be the 
Davies perturbation of $S_t$ for all $t > 0$.
Let $- A_{C,\rho} = - A_\rho$
be the generator of  $S^{(C,\rho)}$.
Moreover, define the form $\gota_{C,\rho}$ by 
\[
\gota_{C,\rho}(u,v)
= \gota_C(U_{-\rho} u, U_\rho v)
\]
with form domain $D(\gota_{C,\rho}) = W^{1,2}(\Ri^d)$.
Then $A_{C,\rho}$ is the operator associated with $\gota_{C,\rho}$.

We frequently need the following lemma for estimates on~$L_2$.

\begin{lemma} \label{llocseh203}
Let $\theta \in [0,\frac{\pi}{2})$ and $\mu,M > 0$.
Then one has the following.
\begin{tabel}
\item \label{llocseh203-1}
If $C \in \cs(\theta,M)$ then 
$\|S^{(\rho)}_t\|_{2 \to 2} \leq e^{\omega \rho^2 t}$ 
for all $t > 0$, $\psi \in \cd_1$ and $\rho \in \Ri$, 
where $\omega = 3 (1 + \tan \theta)^2 (1 + d^2 M)$.
\item \label{llocseh203-2}
Let $Y \subset \Ri^d$ and $\alpha \in (-\theta_\gota, \theta_\gota)$.
Then 
\[
e^{i \alpha} \, C
\in \cs(Y, \theta + |\alpha|, 
    \mu \, \frac{ \cos(\theta + |\alpha|) }{ \cos \theta } , M)
\]
for all $C \in \cs(Y, \theta, \mu, M)$.
\item \label{llocseh203-3}
If $\chi \in W^{1,\infty}(\Ri^d)$ and 
$C \in \cs(\supp \chi, \theta,\mu,M)$, then 
$M_\chi D(A_{\Re C}^{1/2}) \subset W^{1,2}(\Ri^d)$
and 
\begin{equation}
\sum_{m=1}^d \|(\partial_m \, M_\chi - M_{\partial_m \chi})u\|_2^2
\leq \frac{\|\chi\|_\infty^2}{\mu} \, \|A_{\Re C}^{1/2} u\|_2^2
\label{ellocseh203;123}
\end{equation}
for all $u \in D(A_{\Re C}^{1/2})$.
\item \label{llocseh203-4}
If $\chi \in W^{1,\infty}(\Ri^d)$ and 
$C \in \cs(\supp \chi, \theta,\mu,M)$, then 
$M_\chi S_t \, L_2(\Ri^d) \subset W^{1,2}(\Ri^d)$
and 
\[
\|(\partial_m \, M_\chi - M_{\partial_m \chi}) S_t u\|_2
\leq \frac{\|\chi\|_\infty}{\sqrt{\mu \, \sin \theta_\gota}} \, t^{-1/2} \, \|u\|_2
\]
for all $t > 0$, $u \in L_2(\Ri^d)$ and $m \in \{ 1,\ldots,d \} $.
\end{tabel}
\end{lemma}

\begin{remarkn}
Note that  $\partial_m \, M_\chi - M_{\partial_m \chi} \supset M_\chi \, \partial_m$
in Statements~\ref{llocseh203-3} and \ref{llocseh203-4}.

The rotational invariance of Statement~\ref{llocseh203-2} allows to consider various bounds 
on kernels $K_z$ merely for $z \in (0,\infty)$.
Then the uniform bounds for $z$ in a sector $\Sigma_{\theta'}^\circ$
with $\theta' \in (0,\theta_\gota)$ follow since all bounds depend
only on $\chi$, $\Omega$, $\mu$, $M$, $\nu$ and $\theta$.
\end{remarkn}

\vertspace

\noindent
\proof\
`\ref{llocseh203-1}'.
This follows from (14) in \cite{AE2}.

`\ref{llocseh203-2}'.
Let $x \in Y$ and $\xi \in \Ci^d$.
Then 
\begin{eqnarray*}
\RRe \langle e^{i \alpha} \, C(x) \, \xi, \xi \rangle
& = & \cos \alpha \, \RRe \langle C(x) \, \xi, \xi \rangle 
   - \sin \alpha \, \IIm \langle C(x) \, \xi, \xi \rangle  \\
& \geq & \cos \alpha \, \RRe \langle C(x) \, \xi, \xi \rangle 
   - \sin |\alpha| \, \tan \theta \RRe \langle C(x) \, \xi, \xi \rangle  \\
& = & \frac{ \cos(\theta + |\alpha|) }{ \cos \theta } \, 
     \RRe \langle C(x) \, \xi, \xi \rangle  \\
& \geq & \mu \, |\xi|^2 \, 
           \frac{ \cos(\theta + |\alpha|) }{ \cos \theta } 
.
\end{eqnarray*}
This proves Statement~\ref{llocseh203-2}.

`\ref{llocseh203-3}'.
Let $\varepsilon > 0$ and $n \in \Ni$.
Set $C_n = C + \frac{1}{n} \, I$.
Then $C_n \in \ce(\supp \chi, \theta, \mu, M+1)$.
Let $u \in W^{1,2}(\Ri^d)$.
Then 
\begin{eqnarray*}
\sum_{m=1}^d \|(\partial_m \, M_\chi - M_{\partial_m \chi})u\|_2^2
& = & \sum_{m=1}^d \int |\chi \, \partial_m u|^2
\leq \tfrac{1}{\mu} \int \RRe \sum_{k,l=1}^d 
     (c_{kl} + \tfrac{1}{n} \, \delta_{kl}) \, 
           \chi^2 \, (\partial_k u) \, \overline{\partial_l u}   \\
& \leq & \frac{\|\chi\|_\infty^2}{\mu} \, \|A_{\Re C_n}^{1/2} u\|_2^2  \\
& \leq & \frac{\|\chi\|_\infty^2}{\mu} \, \|(\varepsilon I + A_{\Re C_n})^{1/2} u\|_2^2
.
\end{eqnarray*}
Now let $u \in L_2(\Ri^d)$, $v \in W^{1,2}(\Ri^d)$ and 
$m \in \{ 1,\ldots,d \} $.
Then 
\begin{eqnarray*}
\hspace*{-3pt}
|((\varepsilon I + A_{\Re C_n})^{-1/2} u, \partial_m (\chi \, v))|
& = & |( (\partial_m \, M_\chi - M_{\partial_m \chi}) (\varepsilon I + A_{\Re C_n})^{-1/2} u, v)|
\leq \frac{\|\chi\|_\infty^2}{\mu} \, \|u\|_2 \, \|v\|_2
.  
\end{eqnarray*}
Taking the limit $n \to \infty$ and using \cite{AE2} Theorem~3.7 
(cf.\ \cite{bSim5} Theorem~3.2), one deduces that 
\[
|((\varepsilon I + A_{\Re C})^{-1/2} u, \partial_m (\chi \, v))|
\leq \frac{\|\chi\|_\infty^2}{\mu} \, \|u\|_2 \, \|v\|_2
.  \]
So 
\[
|(u, \partial_m (\chi \, v))|
\leq \frac{\|\chi\|_\infty^2}{\mu} \, \|(\varepsilon I + A_{\Re C})^{1/2} u\|_2 \, \|v\|_2
\]
for all $u \in D(A_{\Re C}^{1/2})$, $v \in W^{1,2}(\Ri^d)$ and $m \in \{ 1,\ldots,d \} $.
Taking the limit $\varepsilon \downarrow 0$ one obtains
\begin{equation}
|(u,(M_\chi \, \partial_m + M_{\partial_m \chi}) v)|
\leq \frac{\|\chi\|_\infty^2}{\mu} \, \|A_{\Re C}^{1/2} u\|_2 \, \|v\|_2
.
\label{ellocseh203;12}
\end{equation}
Therefore 
\[
|(M_\chi u, \partial_m v)|
\leq \Big( \frac{\|\chi\|_\infty^2}{\mu} \, \|A_{\Re C}^{1/2} u\|_2 
            + \|\partial_m \chi\|_\infty \, \|u\|_2 \Big) 
    \|v\|_2
\]
for all $v \in W^{1,2}(\Ri^d)$ and $m \in \{ 1,\ldots,d \} $.
So $M_\chi u \in W^{1,2}(\Ri^d)$ and then the estimate (\ref{ellocseh203;123})
follows from (\ref{ellocseh203;12}).

`\ref{llocseh203-4}'.
Let $\alpha \in (0,\theta_\gota)$.
Then it follows from Statements~\ref{llocseh203-1} and \ref{llocseh203-2}
that 
$\|S^{(\rho)}_z\|_{2 \to 2} \leq e^{\omega \, \rho^2 \, |z|}$
for all $z \in \Sigma_\alpha^\circ$, $\rho \in \Ri$ and $\psi \in \cd_1$,
where $\omega = 3 (1 + \tan(\theta + \alpha))^2 (1 + d^2 M)$.
Hence the Cauchy representation formula gives
\begin{equation}
\|A_\rho \, S^{(\rho)}_t\|_{2 \to 2}
\leq \frac{1}{t \, \sin \alpha} \, e^{2 \omega \rho^2 t}
\label{ellocseh203;2}
\end{equation}
for all $t > 0$, $\rho \in \Ri$ and $\psi \in \cd_1$.
If $C \in \ce(\supp \chi, \theta, \mu, M)$ then 
\begin{eqnarray*}
\|(\partial_m \, M_\chi - M_{\partial_m \chi}) S^C_t u\|_2^2
& \leq & \frac{\|\chi\|_\infty^2}{\mu} \, \|A_{\Re C}^{1/2} \, S^C_t u\|_2^2  \\
& = & \frac{\|\chi\|_\infty^2}{\mu} \, \RRe (A_C \, S^C_t u, S^C_t u)  \\
& \leq & \frac{\|\chi\|_\infty^2}{\mu} \, \|A_C \, S^C_t u\|_2 \, \|S^C_t u\|_2  \\
& \leq & \frac{\|\chi\|_\infty^2}{\mu} \, \frac{1}{t \, \sin \alpha} \, \|u\|_2^2
.
\end{eqnarray*}
Now take the limit $\alpha \uparrow \theta_\gota$.
Finally, let $C \in \cs(\supp \chi, \theta, \mu, M)$.
For all $n \in \Ni$ define $C_n = C + \frac{1}{n} \, I$.
Since $\lim_{n \to \infty} S^{C_n}_t = S^C_t$ strongly in $\cl(L_2(\Ri^d))$
by \cite{AE2} Corollary~3.9, now
Statement~\ref{llocseh203-4} follows as in the proof of Statement~\ref{llocseh203-3}.\hfill$\Box$

\vertspace

The next theorem is a uniform version of Theorem~\ref{tlocseh101}.

\begin{thm} \label{tlocseh201}
Let $\Omega \subsetneqq \Ri^d$ be open, $\theta \in [0,\frac{\pi}{2})$,
$\mu,M > 0$, $\nu \in (0,1)$, $\kappa > 0$, $\tau \in [0,1)$, 
$\theta' \in (0,\theta_\gota)$ and 
$\chi \in C^\infty_{\rm b}(\Ri^d)$ with $\chi \neq 0$ and
$d(\supp \chi, \Omega^{\rm c}) > 0$.
Then there exist $a,b > 0$ such that 
for every $C \in \cs^\nu(\Omega,\theta,\mu,M)$ there exists 
a function $(z,x,y) \mapsto K_z(x,y)$ from 
$\Sigma_{\theta_\gota}^\circ \times \Ri^d \times \Ri^d$ into $\Ci$
such that the following is valid.
\begin{tabel}
\item \label{tlocseh201-0.5}
The function $(z,x,y) \mapsto K_z(x,y)$ is continuous 
from $\Sigma_{\theta_\gota}^\circ \times \Ri^d \times \Ri^d$ into $\Ci$.
\item \label{tlocseh201-1}
For all $z \in \Sigma_{\theta'}^\circ$ the function $K_z$ is the kernel of 
the operator $M_\chi S_z M_\chi$.
\item \label{tlocseh201-3}
For all $x,y \in \Ri^d$ the function $z \mapsto K_z(x,y)$ 
is holomorphic from $\Sigma_{\theta'}^\circ$ into $\Ci$.
\item \label{tlocseh201-2}
For all $z \in \Sigma_{\theta'}^\circ$ the
function $K_z$ is once differentiable in each variable and the 
derivative with respect to one variable is differentiable in the 
other variable.
Moreover, for every multi-index $\alpha,\beta$ with 
$0 \leq |\alpha|,|\beta| \leq 1$ one has 
\[
|(\partial_x^\alpha \, \partial_y^\beta \, K_z)(x,y)|
\leq a \, |z|^{-d/2} \, |z|^{-(|\alpha| + |\beta|)/2} \, 
      (1 + |z|)^{\frac{d+|\alpha| +|\beta|}{2}} \,
      e^{-b \, \frac{|x-y|^2}{|z|}} 
\]
and 
\begin{eqnarray*}
\lefteqn{
|(\partial_x^\alpha \, \partial_y^\beta \, K_z)(x+h,y+k) 
     - (\partial_x^\alpha \, \partial_y^\beta \, K_z)(x,y)|
} \hspace{20mm}  \\*
& \leq & a \, |z|^{-d/2} \, |z|^{-(|\alpha| + |\beta|)/2} \, 
    \left( \frac{|h| + |k|}{|x-y| + \sqrt{|z|}} \right)^\nu \, 
     (1 + |z|)^{\frac{d+|\alpha| +|\beta| + \nu}{2}} \, 
      e^{-b \, \frac{|x-y|^2}{|z|}} 
\end{eqnarray*}
for all $x,y,h,k \in \Ri^d$ with 
$|h| + |k| \leq \tau \, |x-y| + \kappa \, \sqrt{|z|}$.
\end{tabel}
\end{thm}
\proof\
We first prove the theorem with $\cs^\nu(\Omega,\theta,\mu,M)$
replaced by $\ce\ch^\nu(\Omega,\theta,\mu,M)$.
For strongly elliptic operators on $\Ri^d$ in divergence form and 
H\"older continuous coefficients all the kernels with stated
holomorphy and continuity properties are well known.
The main point is to derive the uniform bounds.
We emphasise that the constants in the proof do not depend on the ellipticity
constant $\mu_0$ for elements in $\ce\ch^\nu(\Omega,\theta,\mu,M)$,
nor on the H\"older continuity of the coefficients on $\Omega^{\rm c}$.
Then we will approximate elements of $\ce^\nu(\Omega,\theta,\mu,M)$
by elements of $\ce\ch^\nu(\Omega,\theta,\mu,M)$ and 
finally approximate elements of $\cs^\nu(\Omega,\theta,\mu,M)$
by elements of $\ce^\nu(\Omega,\theta,\mu,M)$.

Without loss of generality we may assume that $d \geq 3$.
Let $F \subset \Omega$ be a closed set with $F \neq \emptyset$ and 
$d(F,\Omega^{\rm c}) > 0$.
Let $r_0 = \frac{1}{3} d(F,\Omega^{\rm c})$.
There exist $\chi_1,\chi_2 \in C_{\rm b}^\infty(\Ri^d)$ such that 
$0 \leq \chi_1,\chi_2 \leq \one$ and 
\begin{eqnarray}
\chi_1(x) & = & \left\{ \begin{array}{ll} 
                  0 & \mbox{if } d(x,F) < r_0 ,  \\[5pt]
                  1 & \mbox{if } d(x,\Omega^{\rm c}) < r_0 ,  
                        \end{array} \right. \label{etlocseh201;5} \\
\chi_2(x) & = & \left\{ \begin{array}{ll} 
                  1 & \mbox{if } d(x,F) < 2 r_0 ,  \\[5pt]
                  0 & \mbox{if } x \notin \Omega , 
                        \end{array} \right.  \label{etlocseh201;6}
\end{eqnarray}
for all $x \in \Ri^d$.
Let $M' = 2 \|\chi_2\|_{W^{1,\infty}(\Ri^d)} \, M + \|\chi_1\|_{W^{1,\infty}(\Ri^d)} + 1$.
Let $C \in \ce\ch^\nu(\Omega,\theta,\mu,M)$.
Define $C' = \chi_2 \, C + \chi_1 \, I$.
Then $C' \in \ce\ch^\nu(\Ri^d,\theta,\mu \wedge 1,M')$.
Since the operator $A_{C'}$ is a strongly elliptic operator on $\Ri^d$
with $C^{0,\nu}$-coefficients,
it satisfies various kinds of De Giorgi estimates.
On bounded open sets these are proved by Giaquinta \cite{Gia1}, \cite{GiM},
or Xu--Zuily \cite{XZ}.
The global estimates follow  from \cite{ER15} Proposition~3.5 and 
\cite{ER19} Proposition~2.6.
Precisely, there exist $c_{DG}',c_{DG}'' > 0$ and 
for all $\nu' \in (0,1)$ a $c_{DG} > 0$,
depending only on $\mu$, $M'$, $\nu$ and $\nu'$, such that 
\begin{equation}
\int_{B(x,r)} |\nabla u|^2 
\leq c_{DG} \, \Big( \frac{r}{R} \Big)^{d-2+2\nu'} \int_{B(x,R)} |\nabla u|^2 
\label{eSlocseh2;1}
\end{equation}
and 
\begin{eqnarray}
\lefteqn{
\sum_{k=1}^d \int_{B(x,r)} |\partial_k u - \langle \partial_k u\rangle_{x,r}|^2
} \hspace{10mm} \nonumber  \\*
& \leq & c_{DG}' \, \Big( \frac{r}{R} \Big)^{d+2 \tilde \nu} 
          \sum_{k=1}^d \int_{B(x,R)} |\partial_k u - \langle \partial_k u\rangle_{x,R}|^2 
   + c_{DG}'' \, R^{2 \nu}
          \sum_{k=1}^d \int_{B(x,R)} |\nabla u|^2 
\label{eSlocseh2;2}
\end{eqnarray}
for all $R \in (0,1]$, $r \in (0,R]$ and $u \in W^{1,2}(B(x,R))$
satisfying $A_{C'} u = 0$ weakly on $B(x,R)$, where $\tilde \nu = \frac{1}{2} (1 + \nu)$.

Write $\gota = \gota_C$, $A = A_C$ and $S = S^C$.
It is well known that the semigroup generated by $S$ has a 
kernel $K^S$ satisfying continuity, holomorphy and 
Gaussian properties similar to 
\ref{tlocseh201-0.5}--\ref{tlocseh201-2} in the theorem.
Define $K_z(x,y) = \chi(x) \, K^S_z(x,y) \, \chi(y)$.
Then $K_z$ is the kernel of $M_\chi \, S_z \, M_\chi$.
Moreover, $K$ satisfies Properties~\ref{tlocseh201-0.5}--\ref{tlocseh201-3}.

Note that $\gota(u,v) = \gota_{C'}(u,v)$ for all $u,v \in W^{1,2}(\Ri^d)$ 
with $\supp u \subset F$.
For all $\gamma \in [0,d]$ let $M_{2,\gamma}(\Ri^d)$ be the Morrey space and 
for all $\gamma \in [0,d+2)$ let $\cm_{2,\gamma}(\Ri^d)$ be the 
Campanato space as defined in \cite{EO1} Section~2.

For all $\gamma \in [0,d)$ let $P(\gamma)$ be the hypothesis
\begin{quote}
For all $\chi \in C_{\rm b}^\infty(\Ri^d)$ with $\supp \chi \subset F$
there exist $a_1,\omega_1 > 0$,
depending only on $\chi$, $\Omega$, $\theta$, $\mu$, $M$ and $\nu$, such that 
\[
\|M_\chi \, S^{(\rho)}_t u\|_{M_{2,\gamma}(\Ri^d)}
\leq a_1 \, t^{-\gamma / 4} \, e^{\omega_1 (1 + \rho^2) \, t} \, \|u\|_2
\]
and 
\[
\|\nabla \, M_\chi \, S^{(\rho)}_t u\|_{M_{2,\gamma}(\Ri^d)}
\leq a_1 \, t^{-\gamma / 4} \, t^{-1/2} \, e^{\omega_1 (1 + \rho^2) \, t} \, \|u\|_2
\]
uniformly for all $t > 0$, $u \in L_2(\Ri^d)$, 
$\rho \in \Ri$ and $\psi \in \cd_1$.
\end{quote}

Arguing as in the proof of Lemma~3.3 in \cite{EO1} it follows from (\ref{ellocseh203;2}) 
and the 
De Giorgi estimates~(\ref{eSlocseh2;1}) that $P(\gamma)$ is valid for all $\gamma \in [0,d)$.

Next, for all $\gamma \in [0,d+2 \nu]$ let $P'(\gamma)$ be the hypothesis
\begin{quote}
For all $\chi \in C_{\rm b}^\infty(\Ri^d)$ with $\supp \chi \subset F$
there exist $a_1,\omega_1 > 0$,
depending only on $\chi$, $\Omega$, $\theta$, $\mu$, $M$ and $\nu$, such that 
\[
\|M_\chi \, S^{(\rho)}_t u\|_{\cm_{2,\gamma}(\Ri^d)}
\leq a_1 \, t^{-\gamma / 4} \, e^{\omega_1 (1 + \rho^2) \, t} \, \|u\|_2
\]
and 
\[
\|\nabla \, M_\chi \, S^{(\rho)}_t u\|_{\cm_{2,\gamma}(\Ri^d)}
\leq a_1 \, t^{-\gamma / 4} \, t^{-1/2} \, e^{\omega_1 (1 + \rho^2) \, t} \, \|u\|_2
\]
uniformly for all $t > 0$, $u \in L_2(\Ri^d)$, 
$\rho \in \Ri$ and $\psi \in \cd_2$.
\end{quote}
Since $M_{2,\gamma} \cap L_2 = \cm_{2,\gamma} \cap L_2$ for 
all $\gamma \in [0,d)$, with equivalent norms, one deduces from $P(\gamma)$ that 
also $P'(\gamma)$ is valid for all $\gamma \in [0,d)$.
But arguing as in \cite{ER19}, proof of Proposition~3.2 and 
the proof of Lemma~3.3 in \cite{EO1} it follows from the 
De Giorgi estimates (\ref{eSlocseh2;2}) that $P'(\gamma)$ is valid for all $\gamma \in [0,d+2\nu]$.
Hence there are $a,\omega > 0$,
depending only on $\chi$, $\Omega$, $\theta$, $\mu$, $M$ and $\nu$, such that 
\begin{equation}
\|\partial^\alpha M_\chi S^{(\rho)}_t u\|_\infty
\leq a \, t^{-d/4} \, t^{-|\alpha|/2} \, e^{\omega (1 + \rho^2) \, t} \, \|u\|_2
\label{ellocseh203;8}
\end{equation}
and 
\[
\|(I - L(h)) \partial^\alpha M_\chi S^{(\rho)}_t u\|_\infty
\leq a \, t^{-d/4} \, t^{-|\alpha|/2} \Big( \frac{|h|}{\sqrt{t}} \Big)^\nu 
     \, e^{\omega (1 + \rho^2) \, t} \, \|u\|_2
\]
for all multi-indices $\alpha$ with $|\alpha| \leq 1$, $t > 0$, $u \in L_2(\Ri^d)$, 
$\rho \in \Ri$ and $\psi \in \cd_2$.
Here $L(h)$ denotes left translation, defined by
$(L(h) u)(x) = u(x-h)$.
Next, $\|S^{(\rho)}_t\|_{2 \to 2} \leq e^{\omega_0 \rho^2 t}$
for all $t > 0$, where 
$\omega_0 = 3 (1 + \tan \theta)^2 (1 + d^2 \, M)$ by Lemma~\ref{llocseh203}\ref{llocseh203-1}.

Then the bounds of Property~\ref{tlocseh201-2} follow from Lemma~\ref{llocsehA1}
in Appendix~\ref{AlocsehA} uniformly for all 
$z \in \Sigma_{\theta'}^\circ$ and $C \in \ce\ch^\nu(\Omega,\theta,\mu,M)$.
This proves the theorem with $\cs^\nu(\Omega,\theta,\mu,M)$ replaced
by $\ce\ch^\nu(\Omega,\theta,\mu,M)$.

Let $\chi \in C_{\rm b}^\infty(\Ri^d)$ with $\chi \neq 0$ and 
$d(\supp \chi, \Omega^{\rm c}) > 0$.
Let $r_0 = \frac{1}{2} d(\supp \chi, \Omega^{\rm c})$
and let $\Omega' = \{ x \in \Omega : d(x,\supp \chi) < r_0 \} $.
Then $\Omega'$ is open, 
$d(\supp \chi, (\Omega')^{\rm c}) > 0$ and 
$d(\Omega', \Omega^{\rm c}) \geq r_0> 0$.

Let $\alpha,\beta$ be multi-indices with $|\alpha|,|\beta| \leq 1$.
Using the Cauchy representation formula on the sector $\Sigma_{\theta_\gota}^\circ$
one deduces that $\partial_t \partial_x^\alpha \partial_y^\beta K^C_z$
satisfies H\"older type Gaussian bounds uniformly for all 
$C \in \ce\ch^\nu(\Omega',\theta,\mu,M)$, where $K^C$ is the kernel 
associated with $z \mapsto M_\chi \, S^C_z \, M_\chi$.
Hence the set of functions
\[
\{ (z,x,y) \mapsto (\partial_x^\alpha \partial_y^\beta K^C_z)(x,y) 
        : C \in \ce\ch^\nu(\Omega',\theta,\mu,M) \}
\]
is equicontinuous on compact subsets of $\Sigma_{\theta_\gota}^\circ \times \Ri^d \times \Ri^d$.

Fix $\tau \in C_c^\infty(B(0,r_0))$ with $\tau \geq 0$ and $\int \tau = 1$.
For all $n \in \Ni$ define $\tau_n \in C_c^\infty(\Ri^d)$ by
$\tau_n(x) = n^d \, \tau(n \, x)$.
Now let $C \in \ce^\nu(\Omega,\theta,\mu,M)$.
For all $n \in \Ni$ and $k,l \in \{ 1,\ldots,d \} $ define $c_{kl}^{(n)} = c_{kl} * \tau_n$
and define $C^{(n)} = (c_{kl}^{(n)})$.
Then $C^{(n)} \in \ce\ch^\nu(\Omega',\theta,\mu,M)$ for all $n \in \Ni$.
Write $S^{(n)} = S^{C^{(n)}}$, etc.
Since $ \{ (z,x,y) \mapsto (\partial_x^\alpha \partial_y^\beta K^{(n)}_z)(x,y) : n \in \Ni \} $
is equicontinuous on compact subsets of $\Sigma_{\theta_\gota}^\circ \times \Ri^d \times \Ri^d$
for all $|\alpha|,|\beta| \leq 1$ it follows with a diagonal argument
from the Arz\'ela--Ascoli theorem that there exists a subsequence 
$(K^{(n_k)})_{k \in \Ni}$ of $(K^{(n)})_{n \in \Ni}$ such that 
$K^{(\alpha,\beta)} = \lim_{k \to \infty} \partial_x^\alpha \partial_y^\beta K^{(n_k)}$
exists uniformly on compact subsets of $\Sigma_{\theta_\gota}^\circ \times \Ri^d \times \Ri^d$ and every 
multi-index $\alpha,\beta$ with $|\alpha|,|\beta| \leq 1$.
Then $(z,x,y) \mapsto K^{(\alpha,\beta)}_z(x,y)$ is continuous on 
$\Sigma_{\theta_\gota}^\circ \times \Ri^d \times \Ri^d$.
Set $K = K^{(\alpha,\beta)}$ where $|\alpha| = |\beta| = 0$.
Obviously for every $z \in \Sigma_{\theta_\gota}^\circ$ the function 
$K_z$ is once differentiable in each variable and the 
derivative with respect to one variable is differentiable in the 
other variable.
Also $K_z$ satisfies all the Gaussian bounds from the theorem
and the constants in the Gaussians depends only on 
$\chi$, $\Omega$, $\theta'$, $\mu$, $M$ and $\nu$.
Let $z \in \Sigma_{\theta_\gota}^\circ$.
Let $u,v \in C_c^\infty(\Ri^d)$.
We shall prove below in Lemma~\ref{llocseh202} that 
$\lim_{n \to \infty} S^{(n)}_z u = S^C_z u$ in $L_2(\Ri^d)$.
Hence 
\begin{eqnarray*}
(S^C_z u,v)
& = & \lim_{k \to \infty} (S^{(n_k)}_z u,v)  \\
& = & \lim_{k \to \infty} \int \!\int K^{(n_k)}_z(x,y) \, u(y) \, \overline{v(x)} \, dx \, dy  \\
& = & \int \! \int K_z(x,y) \, u(y) \, \overline{v(x)} \, dx \, dy 
.
\end{eqnarray*}
Since $K$ satisfies Gaussian bounds it follows that $K_z$ is the kernel of $S^C_z$.
This proves the theorem with $\cs^\nu(\Omega,\theta,\mu,M)$ replaced
by $\ce^\nu(\Omega,\theta,\mu,M)$.

Finally, let $C \in \cs^\nu(\Omega,\theta,\mu,M)$.
For all $n \in \Ni$ define $C_n = C + \frac{1}{n} \, I$.
Then $C_n \in \ce^\nu(\Omega,\theta,\mu,M+1)$ for all $n \in \Ni$.
Since $\lim_{n \to \infty} S^{C_n}_z u = S^C_z u$ in $L_2(\Ri^d)$ 
for all $z \in \Sigma_{\theta_\gota}^\circ$ and $u \in L_2(\Ri^d)$
by \cite{AE2} Corollary~3.9, a similar approximation argument as 
in the previous step completes the proof of Theorem~\ref{tlocseh201}.\hfill$\Box$

\vertspace

It remains to show the next lemma.

\begin{lemma} \label{llocseh202}
Let  $\theta \in [0,\frac{\pi}{2})$ and $M > 0$.
Let $C \in \ce(\theta,M)$ and $\tau \in C_c^\infty(\Ri^d)$
with $\tau \geq 0$ and $\int \tau = 1$.
For all $n \in \Ni$ define $\tau_n \in C_c^\infty(\Ri^d)$ by 
$\tau_n(x) = n^d \, \tau(n \, x)$ and set $C^{(n)} = (c_{kl}^{(n)})$,
where $c_{kl}^{(n)} = \tau_n * c_{kl}$ for all $k,l \in \{ 1,\ldots,d \} $.
Then $C^{(n)} \in \ce(\theta,M)$ for all $n \in \Ni$.
Write $S^{(n)} = S^{C^{(n)}}$, $S = S^C$, etc.
Then $\lim_{n \to \infty} S^{(n)}_z = S_z$ strongly in $\cl(L_2(\Ri^d))$
for all $z \in \Sigma_{\theta_\gota}^\circ$.
\end{lemma}
\proof\
Let $\eta \in L_\infty(\Ri^d)$ and $u \in L_2(\Ri^d)$.
Then 
\[
\Big( (I - L(x)) \eta \Big) u
= ( I - L(x)) (\eta \, u) - (L(x) \eta) \, (I - L(x)) u
\]
for all $x \in \Ri^d$.
So 
\[
\|(\eta - \tau_n * \eta) u\|_2
\leq \|\eta \, u - \tau_n * (\eta \, u)\|_2
   + \|\eta\|_\infty \int_{\Ri^d} \tau_n(x) \, \|( I - L(x) ) u\|_2 \, dx
\]
for all $n \in \Ni$ and $\lim \|(\eta - \tau_n * \eta) u\|_2 = 0$.

There exists a $\mu_0 > 0$ such that 
$\RRe \langle C(x) \, \xi, \xi \rangle \geq \mu_0 \, |\xi|^2$ 
for all $x \in \Ri^d$ and $\xi \in \Ci^d$.
Then $\RRe \langle C^{(n)}(x) \, \xi, \xi \rangle \geq \mu_0 \, |\xi|^2$ 
and $\|C^{(n)}(x)\| \leq M$
for all $x \in \Ri^d$, $\xi \in \Ci^d$ and $n \in \Ni$.
Hence there exists a $c > 0$ such that $\|\nabla S^{(n)}_t u\|_2 \leq c \, t^{-1/2} \|u\|_2$
for all $t > 0$, $u \in L_2(\Ri^d)$ and $n \in \Ni$
(cf.\ the proof of (13) in \cite{EO1}).
By increasing $c$ if necessary, it follows similarly that 
$\|\nabla S_t u\|_2 \leq c \, t^{-1/2} \|u\|_2$,
and $\|\nabla S^{(n)*}_t u\|_2 \leq c \, t^{-1/2} \|u\|_2$
for all $t > 0$ and $u \in L_2(\Ri^d)$.

Without loss of generality we may assume that $z \in (0,\infty)$.
Next, let $u,v \in L_2(\Ri^d)$.
Then 
\begin{eqnarray*}
((S^{(n)}_z - S_z)u,v)
& = & \int_0^z \frac{d}{ds} \, (S_{z-s} u, S^{(n)*}_s v) \, ds  \\
& = & \int_0^z (A \, S_{z-s} u, S^{(n)*}_s v) - (S_{z-s} u, A_n^* \, S^{(n)*}_s v)\, ds  \\
& = & \int_0^z \gota((S_{z-s} u, S^{(n)*}_s v) - \gota_n((S_{z-s} u, S^{(n)*}_s v) \, ds  \\
& = & \sum_{k,l=1}^d \int_0^z ( (c_{kl} - c^{(n)}_{kl}) \partial_k S_{z-s} u,
                                  \partial_l S^{(n)*}_s v) \, ds
\end{eqnarray*}
for all $n \in \Ni$.
So 
\[
\|(S^{(n)}_z - S_z)u\|_2
\leq c \, \sum_{k,l=1}^d \int_0^z
    \|(c_{kl} - c^{(n)}_{kl}) \partial_k S_{z-s} u\|_2 \, s^{-1/2} \, ds
\]
and the lemma follows from the Lebesgue dominated convergence theorem.\hfill$\Box$

\vertspace

The proof of Theorem~\ref{tlocseh201} together with Lemma~\ref{llocseh202} gives
estimates which we need in the proof of Proposition~\ref{plocseh310}.

\begin{lemma} \label{locseh205}
Let $\Omega \subsetneqq \Ri^d$ be open, $\theta \in [0,\frac{\pi}{2})$,
$\mu,M > 0$, $\nu \in (0,1)$ and
$\chi \in C^\infty_{\rm b}(\Ri^d)$ with $\chi \neq 0$ and
$d(\supp \chi, \Omega^{\rm c}) > 0$.
Then there exist $a,\omega > 0$ such that 
for every $C \in \cs^\nu(\Omega,\theta,\mu,M)$ one has 
$M_\chi \, S_t u \in W^{1,\infty}(\Ri^d)$ and 
\[
\|\partial^\alpha \, M_\chi \, S_t u\|_\infty
\leq a \, t^{-d/4} \, t^{-|\alpha| / 2} \, e^{\omega t} \, \|u\|_2
\]
for all multi-indices $\alpha$ with $|\alpha| \leq 1$, $t > 0$ and $u \in L_2(\Ri^d)$.
\end{lemma}
\proof\
Let $\chi \in C^\infty_{\rm b}(\Ri^d)$ with $\chi \neq 0$ and
$d(\supp \chi, \Omega^{\rm c}) > 0$.
Let $r_0 = \frac{1}{2} \, d(\supp \chi, \Omega^{\rm c})$ and 
$\Omega' = \{ x \in \Omega : d(x,\supp \chi) < r_0 \} $.
By (\ref{ellocseh203;8}) there exist $a,\omega > 0$ such that 
\begin{equation}
\|\partial^\alpha \, M_\chi \, S^C_t u\|_\infty
\leq a \, t^{-d/4} \, t^{-|\alpha| / 2} \, e^{\omega t} \, \|u\|_2
\label{elocseh205;1}
\end{equation}
for all $C \in \ce\ch^\nu(\Omega,\theta,\mu,M)$, 
multi-indices $\alpha$ with $|\alpha| \leq 1$, $t > 0$ and $u \in L_2(\Ri^d)$.
Let $C \in \ce^\nu(\Omega,\theta,\mu,M)$.
Let $\tau \in C_c^\infty(B(0,r_0))$ with $\tau \geq 0$ and $\int \tau = 1$.
For all $n \in \Ni$ let $C^{(n)}$ be as in Lemma~\ref{llocseh202}.
Then $C^{(n)} \in \ce\ch^\nu(\Omega,\theta,\mu,M)$.
Let $t > 0$ and $u \in L_2(\Ri^d)$.
Then 
\[
|(S^{C^{(n)}}_t u, M_\chi \, \partial^\alpha v)|
= |(\partial^\alpha \, M_\chi \, S^{C^{(n)}}_t u, v)|
\leq a \, t^{-d/4} \, t^{-|\alpha| / 2} \, e^{\omega t} \, \|u\|_2 \, \|v\|_1
\]
for all multi-indices $\alpha$ with $|\alpha| \leq 1$, $v \in C_c^\infty(\Ri^d)$
and $n \in \Ni$.
Now take the limit $n \to \infty$ and use Lemma~\ref{llocseh202}.
It follows that 
\begin{equation}
|(M_\chi \, S^C_t u, \partial^\alpha v)|
\leq a \, t^{-d/4} \, t^{-|\alpha| / 2} \, e^{\omega t} \, \|u\|_2 \, \|v\|_1
\label{elocseh205;2}
\end{equation}
for all $v \in C_c^\infty(\Ri^d)$ and $|\alpha| \leq 1$.
Choosing $|\alpha| = 0$, it follows that 
$M_\chi \, S_t u \in L_\infty(\Ri^d)$.
Next $|\alpha| = 1$ and the density of $C_c^\infty(\Ri^d)$ in 
$W^{1,1}(\Ri^d)$ give that (\ref{elocseh205;2}) is valid for all $v \in W^{1,1}(\Ri^d)$
and $|\alpha| = 1$.
Hence $M_\chi \, S_t u \in W^{1,\infty}(\Ri^d)$ and (\ref{elocseh205;1}) is valid.

Finally, if $C \in \cs^\nu(\Omega,\theta,\mu,M)$ use the approximation
$C_n = C + \frac{1}{n} \, I$ as at the end of the proof of Theorem~\ref{tlocseh201}
and argue similarly.\hfill$\Box$

\section{Real $W^{1,\infty}$-coefficients} \label{Slocseh3}

Let $\Omega \subset \Ri^d$ be open, $\theta \in [0,\frac{\pi}{2})$ and
$\mu,M > 0$.
Define  $\cs^1(\Omega,\theta,\mu,M,{\rm real})$ to
be the set of all 
$C \in \cs^1(\Omega,\theta,\mu,M)$ such that $c_{kl}$ is real valued for all 
$k,l \in \{ 1,\ldots,d \} $.
If $C \in \cs^1(\Omega,\theta,\mu,M,{\rm real})$ then $\gota_C$ is closable.
Moreover, if $\chi \in C_{\rm b}^\infty(\Ri^d)$ with $\supp \chi \subset \Omega$
then there exists an $a > 0$ such that 
$\|M_\chi u\|_{W^{1,2}(\Ri^d)} \leq a \, \|u\|_{D(\gota_C)}$ for all $u \in W^{1,2}(\Ri^d)$.
Hence $M_\chi (D(\overline{\gota_C})) \subset W^{1,2}(\Ri^d)$ and 
$\|M_\chi u\|_{W^{1,2}(\Ri^d)} \leq a \, \|u\|_{D(\overline{\gota_C})}$
for all $u \in D(\overline{\gota_C})$.
In particular, $M_\chi S_t u \in W^{1,2}(\Ri^d)$ for all $t > 0$ and $u \in L_2(\Ri^d)$.

Throughout the remaining of this paper we set 
\[
\tilde c_{kl} = c_{kl} + c_{lk}
\]
for all $k,l \in \{ 1,\ldots,d \} $, whenever $C \in \cs(\theta,M)$.

\begin{lemma} \label{llocseh302}
Let $\Omega \subset \Ri^d$ be open, $\theta \in [0,\frac{\pi}{2})$ and
$\mu,M > 0$.
Let $C \in \ce^1(\Omega,\theta,\mu,M,{\rm real})$.
Let $\omega = 3 d^2 M$.
Then 
\[
\|S^{(\rho)}_t\|_{\infty \to \infty}
\leq e^{\omega (|\rho| + \rho^2) t}
\]
for all $t > 0$, $\rho \in \Ri$ and $\psi \in \cd_2$.
\end{lemma}
\proof\
Let $\rho \in \Ri$ and $\psi \in \cd_2$.
Obviously the form $a_\rho$ is real.
Integration by parts gives 
\begin{eqnarray*}
\gota_\rho(u,v)
& = & \gota(u,v) 
   - \rho \int \sum \tilde c_{kl} \, (\partial_k u) \, (\partial_l \psi) \, \overline v
   - \rho \int \sum (\partial_l c_{kl} ) \, (\partial_k \psi) \, u \, \overline v  \\*
& & \hspace{20mm} {}
   - \rho \int \sum c_{kl} \, (\partial_k \, \partial_l \psi) \, u \, \overline v
   - \rho^2 \int \sum c_{kl} \, (\partial_k \psi) \, (\partial_l \psi) \, u \, \overline v
\end{eqnarray*}
for all $u,v \in W^{1,2}(\Ri^d, \Ri)$.
Then the lemma follows from \cite{Ouh5} Corollary~4.10.\hfill$\Box$

\vertspace

It follows from Lemma~\ref{llocseh302} that the conditions of 
the next lemma are valid for $p=1$.

\begin{lemma} \label{llocseh301}
Let $\Omega \subset \Ri^d$ be open, $\theta \in [0,\frac{\pi}{2})$,
$\mu,M > 0$ and $p \in [1,\infty)$.
Suppose that for all $\varepsilon > 0$ and  $\chi \in C_{\rm b}^\infty(\Ri^d)$
with $\chi \neq 0$ and $d(\supp \chi, \Omega^{\rm c}) > 0$ there exist
$a,\omega > 0$ such that 
\[
\|M_\chi \, S^{(\rho)}_t u\|_p
\leq a \, t^{-\frac{d}{2} (1 - \frac{1}{p})} \, e^{\omega \rho^2 t} \, e^{\varepsilon t} \, \|u\|_1
\]
for all $C \in \ce^1(\Omega,\theta,\mu,M,{\rm real})$, 
$t > 0$, $u \in L_1 \cap L_2$, $\rho \in \Ri$ and $\psi \in \cd_2$.
Then for all $q \in [p,\infty]$, $\nu \in (0,1)$, $\varepsilon > 0$ and $\chi \in C_{\rm b}^\infty(\Ri^d)$
with $\frac{1}{p} - \frac{1}{q} < \frac{1-\nu}{d}$, 
$\chi \neq 0$ and $d(\supp \chi, \Omega^{\rm c}) > 0$ there exist
$a,\omega > 0$ such that 
\begin{eqnarray*}
\|M_\chi \, S^{(\rho)}_t u\|_q
& \leq & a \, t^{-\frac{d}{2} (1 - \frac{1}{q})} \, e^{\omega \rho^2 t} \, e^{\varepsilon t} \, \|u\|_1  \quad \mbox{and}  \\
\|(I - L(h)) M_\chi \, S^{(\rho)}_t u\|_q
& \leq & a \, |h|^\nu \, t^{-\frac{d}{2} (1 - \frac{1}{q})} \, t^{-\frac{\nu}{2}} \, 
    e^{\omega \rho^2 t} \, e^{\varepsilon t} \, \|u\|_1
\end{eqnarray*}
for all $C \in \ce^1(\Omega,\theta,\mu,M,{\rm real})$, 
$t > 0$, $u \in L_1 \cap L_2$, $\rho \in \Ri$, $\psi \in \cd_2$ and $h \in \Ri^d$.
\end{lemma}
\proof\
Let $C \in \ce^1(\Omega,\theta,\mu,M)$, $\varepsilon > 0$,
$\rho \in \Ri$, $\psi \in \cd_2$,
$\chi \in C_{\rm b}^\infty(\Ri^d)$
and suppose $\chi \neq 0$ and $d(\supp \chi, \Omega^{\rm c}) > 0$.
Then 
\begin{eqnarray}
\gota_\rho(M_\chi \, u, v) - \gota_\rho(u, M_\chi \, v)
& = & - \sum (M_{\partial_l \chi} \, \partial_k \, M_{\widetilde c_{kl}} u, v)
  + \sum (M_{ (\partial_l \chi) (\partial_k c_{kl}) } u, v) 
\label{ellocseh301;1} \\*
& & \hspace{10mm} {}
  - \sum (M_{\partial_l \partial_k \chi} \, M_{c_{kl}} u, v)
  - \rho \sum (M_{\widetilde c_{kl}} \, M_{\partial_k \psi} \, M_{\partial_l \chi} u, v)
\nonumber
\end{eqnarray}
for all $u,v \in W^{1,2}(\Ri^d)$.
Note that $(\partial_l \chi) (\partial_k c_{kl})$ is a bounded function on 
$\Ri^d$ since $\supp \chi \subset \Omega$ and $c_{kl} \in W^{1,\infty}(\Omega)$.
Since $C$ is elliptic for all $u,v \in L_2(\Ri^d)$ one deduces
\begin{eqnarray*}
(M_\chi \, S^{(\rho)}_t u - S^{(\rho)}_t \, M_\chi u, v)
& = & - \int_0^t \frac{d}{ds} (M_\chi \, S^{(\rho)}_{t-s} u, S^{(\rho)*}_s v) \, ds  \\
& = & \int_0^t \Big( - (M_\chi \, A_\rho \, S^{(\rho)}_{t-s} u, S^{(\rho)*}_s v) 
   + (M_\chi \, S^{(\rho)}_{t-s} u, A_\rho^* \, S^{(\rho)*}_s v) \Big) \, ds  \\
& = & \int_0^t \gota_\rho(M_\chi \, S^{(\rho)}_{t-s} u, S^{(\rho)*}_s v)
    - \gota_\rho(S^{(\rho)}_{t-s} u, M_\chi \, S^{(\rho)*}_s v) \, ds  \\
& = & - \sum_{k,l=1}^d \int_0^t (S^{(\rho)}_s M_{\partial_l \chi} \, \partial_k \, 
                                  M_{\widetilde c_{kl}} S^{(\rho)}_{t-s} u, v) \, ds + R,
\end{eqnarray*}
where $R$ is the contribution of the last three terms in (\ref{ellocseh301;1}),
which do not have a derivative.
Therefore 
\begin{eqnarray}
M_{\chi^2} \, S^{(\rho)}_t - M_\chi \, S^{(\rho)}_t \, M_\chi
& = & M_\chi \, [M_\chi, S^{(\rho)}_t]  \nonumber  \\
& = & - \sum_{k,l=1}^d \int_0^t 
    M_\chi \, S^{(\rho)}_s \, M_{\partial_k \chi} \partial_l \, M_{\tilde c_{kl}} \, S^{(\rho)}_{t-s} \, ds
   + R',
\label{ellocseh301;2}
\end{eqnarray}
with $R'$ the contribution of the last three terms in (\ref{ellocseh301;1}).

Fix $k,l \in \{ 1,\ldots,d \} $.
Then 
\begin{eqnarray}
\lefteqn{
\int_0^t 
    M_\chi \, S^{(\rho)}_s \, M_{\partial_k \chi} \partial_l \, M_{\tilde c_{kl}} \, S^{(\rho)}_{t-s} \, ds
} \hspace{20mm}  \nonumber \\*
& = & \int_0^{t/2} 
    M_\chi \, S^{(\rho)}_s \, M_{\partial_k \chi} \partial_l \, M_{\tilde c_{kl}} \, S^{(\rho)}_{t-s} \, ds
   + \int_{t/2}^t
    M_\chi \, S^{(\rho)}_s \, M_{\partial_k \chi} \partial_l \, M_{\tilde c_{kl}} \, S^{(\rho)}_{t-s} \, ds
\label{ellocseh301;3}
\end{eqnarray}
for all $t > 0$.
By Theorem~\ref{tlocseh201} there are $a,\omega > 0$ such that 
\begin{eqnarray*}
\|M_\chi \, S^{(\rho)}_s \, M_{\partial_k \chi} \partial_l\|_{1 \to q}
& \leq & a \, s^{-\frac{d}{2} (1 - \frac{1}{q})} \, s^{-1/2} \, e^{\omega \rho^2 s} \, e^{\varepsilon s}  
      \quad \mbox{and}  \\
\|(I - L(h)) M_\chi \, S^{(\rho)}_s \, M_{\partial_k \chi} \partial_l\|_{1 \to q}
& \leq & a \, |h|^\nu \, s^{-\frac{d}{2} (1 - \frac{1}{q})} \, s^{-\frac{1+\nu}{2}} \, 
     e^{\omega \rho^2 s} \, e^{\varepsilon s} 
\end{eqnarray*}
for all $s > 0$, $\rho \in \Ri$, $\psi \in \cd_2$ and $h \in \Ri^d$.
Suppose from now on that $C$ is real valued.
Since the matrix of coefficients is real  it follows from Lemma~\ref{llocseh302} 
that there exists an $\omega' > 0$ such that 
$\|S^{(\rho)}_t\|_{1 \to 1} \leq e^{\omega' \rho^2 t} \, e^{\varepsilon t}$ 
for all $t > 0$,
$\rho \in \Ri$ and $\psi \in \cd_2$.
Then 
\begin{eqnarray*}
\lefteqn{
\|\int_{t/2}^t
    M_\chi \, S^{(\rho)}_s \, M_{\partial_k \chi} \partial_l \, M_{\tilde c_{kl}} \, S^{(\rho)}_{t-s} \, u \, ds\|_q
} \hspace{40mm} \\*
& \leq & 2 a \, M \int_{\frac{t}{2}}^t 
    s^{-\frac{d}{2} (1 - \frac{1}{q})} \, s^{-1/2} \, 
       e^{\omega \rho^2 s} \, e^{\varepsilon s} \, e^{\omega' \rho^2 (t-s)} \, e^{\varepsilon (t-s)}
     \, \|u\|_1 \, ds \\
& \leq & 2^{d+1} \, a \, M \, t^{-\frac{d}{2} (1 - \frac{1}{q})} 
   \, t^{1/2} e^{(\omega + \omega') \rho^2 t} \, e^{\varepsilon t} \, \|u\|_1  \\
& \leq & 2^{d+1} \, a \, \varepsilon^{-1/2} \, M \, t^{-\frac{d}{2} (1 - \frac{1}{q})} 
    e^{(\omega + \omega') \rho^2 t} \, e^{2\varepsilon t} \, \|u\|_1 
\end{eqnarray*}
for all $t > 0$, $u \in L_1 \cap L_2$, $\rho \in \Ri$ and $\psi \in \cd_2$.
Similarly,
\begin{eqnarray*}
\lefteqn{
\|\int_{t/2}^t
    (I - L(h)) M_\chi \, S^{(\rho)}_s \, M_{\partial_k \chi} \partial_l \, M_{\tilde c_{kl}} \, S^{(\rho)}_{t-s} \, u \, ds\|_q
} \hspace{40mm} \\*
& \leq & 2^{d+1} \, a \, \varepsilon^{-1/2} M \, |h|^\nu \, 
    t^{-\frac{d}{2} (1 - \frac{1}{q})} \, t^{-\frac{\nu}{2}}
    e^{(\omega + \omega') \rho^2 t} \, e^{2 \varepsilon t} \, \|u\|_1 
\end{eqnarray*}
for all $t > 0$, $u \in L_1 \cap L_2$, $\rho \in \Ri$, $\psi \in \cd_2$ and $h \in \Ri^d$.

Next we estimate the $L_1 \to L_q$ norm of the first term in (\ref{ellocseh301;3}).
There exists a $\tilde \chi \in C_{\rm b}^\infty(\Ri^d)$ such that 
$\tilde \chi(x) = 1$ for all $x \in \supp \chi$ and 
$d(\supp \tilde \chi, \Omega^{\rm c}) > 0$.
Then by assumption there exist $a,\omega > 0$ such that 
\[
\|\tilde \chi \, S^{(\rho)}_t u\|_p
\leq a \, t^{-\frac{d}{2} (1 - \frac{1}{p})} \, e^{\omega \rho^2 t} \, e^{\varepsilon t} \, \|u\|_1
\]
for all $t > 0$, $u \in L_1 \cap L_2$, $\rho \in \Ri$ and $\psi \in \cd_2$.
By Theorem~\ref{tlocseh201} there are $a',\omega' > 0$ such that 
\begin{eqnarray*}
\|M_\chi \, S^{(\rho)}_s \, M_{\partial_k \chi} \partial_l\|_{p \to q}
& \leq & a' \, s^{-\frac{d}{2} (\frac{1}{p} - \frac{1}{q})} \, s^{-1/2} \, 
     e^{\omega' \rho^2 s} \, e^{\varepsilon s} 
     \quad \mbox{and}  \\
\|(I - L(h)) M_\chi \, S^{(\rho)}_s \, M_{\partial_k \chi} \partial_l\|_{p \to q}
& \leq & a' \, |h|^\nu \, s^{-\frac{d}{2} (\frac{1}{p} - \frac{1}{q})} \, 
     s^{-\frac{1+\nu}{2}} \, e^{\omega' \rho^2 s} \, e^{\varepsilon s} 
\end{eqnarray*}
for all $s > 0$, $\rho \in \Ri$, $\psi \in \cd_2$ and $h \in \Ri^d$.
Then 
\begin{eqnarray*}
\lefteqn{
\|\int_0^{t/2}
    M_\chi \, S^{(\rho)}_s \, M_{\partial_k \chi} \partial_l \, M_{\tilde c_{kl}} \, S^{(\rho)}_{t-s} \, u \, ds\|_q
} \hspace{20mm}  \\*
& \leq & \int_0^{t/2}
    \|M_\chi \, S^{(\rho)}_s \, M_{\partial_k \chi} \partial_l\|_{p \to q} \, \|M_{\tilde c_{kl}} \, M_{\widetilde \chi} \, S^{(\rho)}_{t-s} \, u\|_p \, ds  \\
& \leq & 2 a \, a' \, M \int_0^{\frac{t}{2}} 
    s^{-\frac{d}{2} (\frac{1}{p} - \frac{1}{q})} \, s^{-1/2} \, e^{\omega' \rho^2 s} \, e^{\varepsilon s} 
    (t-s)^{-\frac{d}{2} (1 - \frac{1}{q})} \,  
     e^{\omega \rho^2 (t-s)} \, e^{\varepsilon (t-s)} \, \|u\|_1 \, ds \\
& \leq & \frac{2^{\frac{d}{2} (\frac{1}{p} - \frac{1}{q}) + 1} a \, a' \, \varepsilon^{-1/2} \, M}
              {\frac{1}{2} - \frac{d}{2} (\frac{1}{p} - \frac{1}{q})}
   \, t^{-\frac{d}{2} (1 - \frac{1}{p})} e^{(\omega + \omega') \rho^2 t} \, e^{2\varepsilon t} \, \|u\|_1
\end{eqnarray*}
for all $t > 0$, $u \in L_1 \cap L_2$, $\rho \in \Ri$ and $\psi \in \cd_2$.
Similarly
\begin{eqnarray*}
\lefteqn{
\|\int_0^{t/2}
   (I - L(h)) M_\chi \, S^{(\rho)}_s \, M_{\partial_k \chi} \partial_l \, M_{\tilde c_{kl}} \, S^{(\rho)}_{t-s} \, u \, ds\|_q
} \hspace{40mm}  \\*
& \leq & \frac{2^{\frac{d}{2} (\frac{1}{p} - \frac{1}{q}) + 1} a \, a' \, \varepsilon^{- \frac{1}{2}} M \, |h|^\nu}
              {\frac{1-\nu}{2} - \frac{d}{2} (\frac{1}{p} - \frac{1}{q})}
   \, t^{-\frac{d}{2} (1 - \frac{1}{p})} t^{-\frac{\nu}{2}} \, 
     e^{(\omega + \omega') \rho^2 t}  \, e^{2\varepsilon t} \, \|u\|_1
\end{eqnarray*}
for all $t > 0$, $u \in L_1 \cap L_2$, $\rho \in \Ri$, $\psi \in \cd_2$ and $h \in \Ri^d$.

The term $R'$ in (\ref{ellocseh301;2}) can be estimated similarly.
Using Theorem~\ref{tlocseh201} to estimate the $L_1 \to L_q$ norms of 
$M_\chi \, S^{(\rho)}_t \, M_\chi$ 
and $(I - L(h)) M_\chi \, S^{(\rho)}_t \, M_\chi$
one deduces that there are $a,\omega > 0$ such that 
\begin{eqnarray*}
\|M_{\chi^2} \, S^{(\rho)}_t u\|_q
& \leq & a \, t^{-\frac{d}{2} (1 - \frac{1}{q})} \, e^{\omega \rho^2 t} \, e^{\varepsilon t} \, \|u\|_1
   \quad \mbox{and}  \\
\|(I - L(h)) M_{\chi^2} \, S^{(\rho)}_t u\|_q
& \leq & a \, |h|^\nu \, t^{-\frac{d}{2} (1 - \frac{1}{q})} \, t^{-\frac{\nu}{2}} \, 
    e^{\omega \rho^2 t} \, e^{\varepsilon t} \, \|u\|_1
\end{eqnarray*}
for all $t > 0$, $u \in L_1 \cap L_2$, $\rho \in \Ri$, $\psi \in \cd_2$ and $h \in \Ri^d$.
Then the lemma follows.\hfill$\Box$

\begin{lemma} \label{llocseh311}
Let $\Omega \subsetneqq \Ri^d$ be open, $\theta \in [0,\frac{\pi}{2})$,
$\mu,M > 0$, $\nu \in (0,1)$, $\varepsilon > 0$  and $\chi \in C_{\rm b}^\infty(\Ri^d)$
with $\chi \neq 0$ and $d(\supp \chi, \Omega^{\rm c}) > 0$.
There there exist $a > 0$ and $\omega \geq 0$ such that
\begin{eqnarray*}
\|M_\chi \, S^{(\rho)}_t  u\|_\infty
& \leq & a \, t^{-d/2} \, e^{\omega \rho^2 t} \, e^{\varepsilon t} \, \|u\|_1
   \quad \mbox{and}  \\
\|(I - L(h)) M_\chi \, S^{(\rho)}_t u\|_\infty
& \leq & a \, |h|^\nu \, t^{-d/2} \, t^{-\nu/2} \, 
    e^{\omega \rho^2 t} \, e^{\varepsilon t} \, \|u\|_1
\end{eqnarray*}
for all $C \in \cs^1(\Omega,\theta,\mu,M,{\rm real})$,
$t > 0$, $u \in L_1 \cap L_2$, $\rho \in \Ri$, $\psi \in \cd_2$ and $h \in \Ri^d$.
In particular, $M_\chi \, S_t u \in C^{0,\nu}(\Ri^d)$ for all $t > 0$ and $u \in L_1(\Ri^d)$.
\end{lemma}
\proof\
It follows by induction from Lemmas~\ref{llocseh302} and \ref{llocseh301}
that the current lemma is valid if $\cs^1(\Omega,\theta,\mu,M,{\rm real})$
is replaced by $\ce^1(\Omega,\theta,\mu,M,{\rm real})$.
Then by approximating $C$ by $C + \frac{1}{n} \, I$ the lemma follows.\hfill$\Box$

\vertspace

Next we turn to derivatives of the semigroup.

\begin{lemma} \label{llocseh305}
Let $u \in W^{1,2}(\Ri^d)$, $v \in L_2(\Ri^d)$ and $k \in \{ 1,\ldots,d \} $.
Then 
\[
(\partial_k u,v)
= \frac{1}{c_0} \, \int_0^\infty \frac{((I - L(r \, e_k))^2 u, v)}{r} \, \frac{dr}{r}
= \frac{1}{c_0} \, \int_0^\infty \frac{(I - L(r \, e_k)) u, (I - L(-r \, e_k)v)}{r} \, \frac{dr}{r}
,  \]
where $c_0 = \int_0^\infty \frac{(1-e^{-r})^2}{r} \, \frac{dr}{r}$.
\end{lemma}
\proof\
Write $D = \frac{d}{dx_k}$, the skew-adjoint operator in $L_2(\Ri^d)$.
Then $((I - L(r \, e_k))^2 u, v) = ((I - e^{-rD})^2 u, v)$.
Now the lemma follows from Fourier theory (or spectral theory).\hfill$\Box$

\begin{lemma} \label{llocseh306}
Let $\Omega \subsetneqq \Ri^d$ be open, $\theta \in [0,\frac{\pi}{2})$,
$\mu,M > 0$, $\nu \in (0,1)$, $\varepsilon \in (0,1]$ and $\chi \in C_{\rm b}^\infty(\Ri^d)$
with $\chi \neq 0$ and $d(\supp \chi, \Omega^{\rm c}) > 0$.
Then there exists an $a > 0$ such that
\begin{eqnarray*}
\| \partial_m \, M_\chi S_t\|_{1 \to \infty}
   & \leq & a \, t^{-d/2} \, t^{-1/2} \, e^{\varepsilon t}  \quad \mbox{and}  \\
\|(I - L(h)) \partial_m \, M_\chi S_t\|_{1 \to \infty}
& \leq & a \, t^{-d/2} \, t^{-1/2} \, 
    \left( \frac{|h|}{\sqrt{t}} \right)^\nu \, e^{\varepsilon t} 
\end{eqnarray*}
for all $C \in \ce^1(\Omega,\theta,\mu,M,{\rm real})$,
$t > 0$, $m \in \{ 1,\ldots,d \} $ and $h \in \Ri^d$.
\end{lemma}
\proof\
We only prove the second estimate, the proof of the first one is similar.
We argue as in the proof of Lemma~\ref{llocseh301} and use the commutator~(\ref{ellocseh301;2}).
There exists a $\tilde \chi \in C_{\rm b}^\infty(\Ri^d)$ such that 
$\tilde \chi(x) = 1$ for all $x \in \supp \chi$ and 
$d(\supp \tilde \chi, \Omega^{\rm c}) > 0$.
Now we have 
\begin{eqnarray}
\lefteqn{
(I - L(h)) \partial_m \Big( M_{\chi^2} \, S_t - M_\chi \, S_t \, M_\chi \Big)
} \hspace*{10mm}  \nonumber  \\*
& = & - \sum_{k,l=1}^d \int_0^t 
    (I - L(h)) \partial_m M_\chi \, S_s \, M_{\partial_k \chi} \partial_l 
          \, M_{\tilde c_{kl}} \, M_{\widetilde \chi} \, S_{t-s} \, ds
   + R,
\label{ellocseh306;1}
\end{eqnarray}
where $R$ is the contribution of the second and third term in (\ref{ellocseh301;1}),
which do not have a derivative. 
Note that the last term in (\ref{ellocseh301;1}) vanishes since $\rho = 0$.
Again we split the integral in two parts.
Fix $k,l \in \{ 1,\ldots,d \} $.
Then 
\begin{eqnarray*}
\lefteqn{
\Big\| \int_{t/2}^t (I - L(h)) \partial_m \, M_\chi \, S_s \, M_{\partial_k \chi} \partial_l 
          \, M_{\tilde c_{kl}} \, M_{\widetilde \chi} \, S_{t-s} \, ds \Big\|_{1 \to \infty}  
} \hspace*{10mm}  \\*
& \leq & 
\int_{t/2}^t \|(I - L(h)) \partial_m \, M_\chi \, S_s \, M_{\partial_k \chi} \partial_l \|_{1 \to \infty}
          \, \|M_{\tilde c_{kl}} \, M_{\widetilde \chi} \, S_{t-s}\|_{1 \to 1} \, ds  \\
& \leq & a \, \, t^{-d/2} \, 
    \left( \frac{|h|}{\sqrt{t}} \right)^\nu \, e^{\varepsilon t} 
\end{eqnarray*}
for a suitable $a > 0$, by the estimates of Theorem~\ref{tlocseh201} and 
Lemma~\ref{llocseh302}.

For the integral over $(0,\frac{t}{2})$ we use Lemma~\ref{llocseh305} and write
\begin{eqnarray*}
\lefteqn{
\Big\| \int_0^{t/2} (I - L(h)) \partial_m \, M_\chi \, S_s \, M_{\partial_k \chi} \partial_l 
          \, M_{\tilde c_{kl}} \, M_{\widetilde \chi} \, S_{t-s} \, ds \Big\|_{1 \to \infty}  
} \hspace*{1mm}  \\*
& = & c_0^{-1} \sup_{\scriptstyle u,v \in W^{1,2}(\Ri^d) \cap L_1(\Ri^d) \atop
               \scriptstyle \|u\|_1, \|v\|_1 \leq 1}
\Big| \int_0^{t/2} \!\int_0^\infty ( (I - L(r e_l)) M_{\tilde c_{kl}} \, M_{\widetilde \chi} \, S_{t-s} u,   \\*[-10pt]
& & \hspace{60mm}
               (I - L(-r e_l)) M_{\partial_k \chi} S^*_s M_\chi \partial_m (I - L(-h)) v )
      \, \frac{dr}{r^2} \, ds \Big|  \\
& \leq & c_0^{-1} \sup_{\scriptstyle u,v \in W^{1,2}(\Ri^d) \cap L_1(\Ri^d) \atop
               \scriptstyle \|u\|_1, \|v\|_1 \leq 1}
   \int_0^{t/2}\! \int_0^\infty \|(I - L(r e_l)) M_{\tilde c_{kl}} \, M_{\widetilde \chi} \, S_{t-s} u\|_\infty \cdot \\*[-10pt]
& & \hspace{60mm}
                 \|(I - L(-r e_l)) M_{\partial_k \chi} S^*_s M_\chi \partial_m (I - L(-h)) v \|_1
      \, \frac{dr}{r^2} \, ds 
,
\end{eqnarray*}
where $c_0$ is as in Lemma~\ref{llocseh305}.
Next split the integral over $(0,\infty)$ in two parts: $(0,1]$ and $[1,\infty)$.
There exist $\nu_1,\nu_2 \in (0,1)$ such that $\nu_2 + \nu < 1$ and $\nu_1 + \nu_2 > 1$.
By Theorem~\ref{tlocseh201} and Lemma~\ref{llocseh311} there exists an $a > 0$ such that 
\begin{eqnarray*}
\|(I - L(h_1)) M_{\tilde c_{kl}} \, M_{\widetilde \chi} \, S_s\|_{1 \to \infty}
& \leq & a \, s^{-d/2} \, 
    \left( \frac{|h_1|}{\sqrt{s}} \right)^{\nu_1} \, e^{\varepsilon s} ,  \\
\|M_{\tilde c_{kl}} \, M_{\widetilde \chi} \, S_s\|_{1 \to \infty}
& \leq & a \, s^{-d/2} \, e^{\varepsilon s} ,  \\
\|(I - L(h_1)) M_{\partial_k \chi} S^*_s M_\chi \partial_m (I - L(h_2))\|_{1 \to 1}
& \leq & a \, s^{-1/2} \left( \frac{|h_1|}{\sqrt{s}} \right)^{\nu_2}
      \left( \frac{|h_2|}{\sqrt{s}} \right)^\nu \, e^{\varepsilon s} \quad \mbox{and}  \\
\|M_{\partial_k \chi} S^*_s M_\chi \partial_m (I - L(h_2))\|_{1 \to 1}
& \leq & a \, s^{-1/2} 
      \left( \frac{|h_2|}{\sqrt{s}} \right)^\nu \, e^{\varepsilon s}   
\end{eqnarray*}
for all $s > 0$ and $h_1,h_2 \in \Ri^d$.
Let $u,v \in W^{1,2}(\Ri^d) \cap L_1(\Ri^d)$.
Then 
\begin{eqnarray*}
\lefteqn{
\int_0^{t/2}\! \int_0^1 \|(I - L(r e_l)) M_{\tilde c_{kl}} \, M_{\widetilde \chi} \, S_{t-s} u\|_\infty 
                 \|(I - L(-r e_l)) M_{\partial_k \chi} S^*_s M_\chi \partial_m (I - L(-h)) v \|_1
      \, \frac{dr}{r^2} \, ds 
} \hspace*{1mm}  \\*
& \leq & a^2 \int_0^{t/2} \!\int_0^1 \, (t-s)^{-d/2} \, \left( \frac{r}{\sqrt{t-s}} \right)^{\nu_1} \, 
     e^{\varepsilon (t-s)}
     s^{-1/2} \left( \frac{|h|}{\sqrt{s}} \right)^\nu
      \left( \frac{r}{\sqrt{s}} \right)^{\nu_2} \, e^{\varepsilon s} \, \|u\|_1 \, \|v\|_1 
      \, \frac{dr}{r^2} \, ds  \\
& = & a' \, t^{-d/2} \, t^{-1/2} \, \left( \frac{|h|}{\sqrt{t}} \right)^\nu \, e^{\varepsilon t} \, 
    t^{\frac{2 - (\nu_1 + \nu_2)}{2}}  \, \|u\|_1 \, \|v\|_1  \\
& \leq & a' \, \varepsilon^{- \frac{2 - (\nu_1 + \nu_2)}{2} } \, t^{-d/2} \, t^{-1/2} \, \left( \frac{|h|}{\sqrt{t}} \right)^\nu \, e^{2 \varepsilon t} 
     \, \|u\|_1 \, \|v\|_1  
\end{eqnarray*}
for a suitable $a' > 0$.
Similarly, since left translations are isometries, one deduces
\begin{eqnarray*}
\lefteqn{
\int_0^{t/2}\! \int_1^\infty \|(I - L(r e_l)) M_{\tilde c_{kl}} \, M_{\widetilde \chi} \, S_{t-s} u\|_\infty 
                 \|(I - L(-r e_l)) M_{\partial_k \chi} S^*_s M_\chi \partial_m (I - L(-h)) v \|_1
      \, \frac{dr}{r^2} \, ds 
} \hspace*{40mm}  \\*
& \leq & 4 a^2 \int_0^{t/2}\! \int_1^\infty (t-s)^{-d/2} \, e^{\varepsilon (t-s)} 
    s^{-1/2} 
      \left( \frac{|h|}{\sqrt{s}} \right)^\nu \, e^{\varepsilon s} \, \|u\|_1 \, \|v\|_1
      \, \frac{dr}{r^2} \, ds  \\
& = & a'' \, t^{-d/2} \, t^{1/2} \, \left( \frac{|h|}{\sqrt{t}} \right)^\nu \, e^{\varepsilon t} 
   \|u\|_1 \, \|v\|_1  \\
& \leq & a'' \, \varepsilon^{-1} \, t^{-d/2} \, t^{-1/2} \, \left( \frac{|h|}{\sqrt{t}} \right)^\nu \, e^{2 \varepsilon t} 
     \, \|u\|_1 \, \|v\|_1  
\end{eqnarray*}
for a suitable $a'' > 0$.
As before, the contribution of the term $R$ in (\ref{ellocseh306;1}) can be estimated
similarly and we leave the rest of the proof to the reader.\hfill$\Box$

\vertspace

We next replace $\ce^1(\Omega,\theta,\mu,M,{\rm real})$ by $\cs^1(\Omega,\theta,\mu,M,{\rm real})$
in Lemma~\ref{llocseh306}.

\begin{lemma} \label{llocseh313}
Let $\Omega \subsetneqq \Ri^d$ be open, $\theta \in [0,\frac{\pi}{2})$,
$\mu,M > 0$, $\nu \in (0,1)$, $\varepsilon \in (0,1]$ and $\chi \in C_{\rm b}^\infty(\Ri^d)$
with $\chi \neq 0$ and $d(\supp \chi, \Omega^{\rm c}) > 0$.
There there exists an $a > 0$ such that
\begin{eqnarray*}
\| \partial_m \, M_\chi S_t\|_{1 \to \infty}
   & \leq & a \, t^{-d/2} \, t^{-1/2} \, e^{\varepsilon t}  \quad \mbox{and}  \\
\|(I - L(h)) \partial_m \, M_\chi S_t\|_{1 \to \infty}
& \leq & a \, t^{-d/2} \, t^{-1/2} \, 
    \left( \frac{|h|}{\sqrt{t}} \right)^\nu \, e^{\varepsilon t} 
\end{eqnarray*}
for all $C \in \cs^1(\Omega,\theta,\mu,M,{\rm real})$,
$t > 0$, $m \in \{ 1,\ldots,d \} $ and $h \in \Ri^d$.
In particular, $M_\chi S_t u \in W^{1+\nu,\infty}(\Ri^d)$ 
for all $t > 0$ and $u \in L_1(\Ri^d)$.
\end{lemma}
\proof\
By Lemma~\ref{llocseh306} there exists an $a > 0$ such that
\begin{equation}
\|\partial_m \, M_\chi S_t\|_{1 \to \infty}
\leq a \, t^{-d/2} \, t^{-1/2} \, e^{\varepsilon t}
\label{ellocseh313;1}
\end{equation}
for all $C \in \ce^1(\Omega,\theta,\mu,M+1,{\rm real})$,
$t > 0$ and $m \in \{ 1,\ldots,d \} $.
Fix $C \in \cs^1(\Omega,\theta,\mu,M,{\rm real})$.
For all $n \in \Ni$ define $C^{(n)} = C + \frac{1}{n} \, I$.
Let $t > 0$, $m \in \{ 1,\ldots,d \} $ and $u,v \in C_c^\infty(\Ri^d)$.
Then it follows from (\ref{ellocseh313;1}) that 
\[
|(M_\chi \, S^{C^{(n)}}_t u, \partial_m v)|
\leq a \, t^{-d/2} \, t^{-1/2} \, e^{\varepsilon t} \, \|u\|_1 \, \|v\|_1
\]
for all $n \in \Ni$.
Since $\lim S^{C^{(n)}}_t = S^C_t$ strongly in $L_2(\Ri^d)$ by 
\cite{AE2} Corollary~3.9 it follows that 
\begin{equation}
|(M_\chi \, S^C_t u, \partial_m v)|
\leq a \, t^{-d/2} \, t^{-1/2} \, e^{\varepsilon t} \, \|u\|_1 \, \|v\|_1
\label{ellocseh313;2}
\end{equation}
for all $u,v \in C_c^\infty(\Ri^d)$.
By continuity it then follows that (\ref{ellocseh313;2}) is 
valid for all $u \in L_1(\Ri^d)$ and $v \in C_c^\infty(\Ri^d)$.
This implies that $M_\chi \, S^C_t u \in W^{1,\infty}(\Ri^d)$ and
the first estimate of the lemma is valid.
The second one follows similarly.\hfill$\Box$

\vertspace

It follows from Lemma~\ref{llocseh311} that for each $t > 0$ there exists 
a measurable function $L_t \colon \Ri^d \times \Ri^d \to \Ri$ such that 
$L_t$ is a kernel of $M_\chi \, S_t$ and 
$L_t$ satisfies Gaussian bounds. 
But then it is unclear whether $x \mapsto L_t(x,y)$ is H\"older 
continuous or differentiable for some $y \in \Ri^d$.
Even worse, there is no reason that the combined map
$(t,x,y) \mapsto L_t(x,y)$ from $(0,\infty) \times \Ri^d \times \Ri^d$
into $\Ri$ is measurable.
In order to circumvent this measurability problem with the 
uncountable many null sets, 
we first obtain a measurable  map on 
$(0,\infty) \times \Ri^d \times \Ri^d$  for the kernels
of $M_\chi \, S_t$ and its derivatives $\partial_m \, M_\chi \, S_t$
and then consider continuity and differential 
properties in Theorem~\ref{tlocseh309}.

\begin{prop} \label{plocseh310}
Let $\Omega \subsetneqq \Ri^d$ be open, $\theta \in [0,\frac{\pi}{2})$,
$\mu,M > 0$ and $\chi \in C_{\rm b}^\infty(\Ri^d)$
with $\chi \neq 0$, $\varepsilon > 0$ and $d(\supp \chi, \Omega^{\rm c}) > 0$.
There there exist $a,b > 0$ such that
for all $C \in \cs^1(\Omega,\theta,\mu,M,{\rm real})$ 
and $m \in \{ 1,\ldots,d \} $
there exist measurable functions $(t,x,y) \mapsto K_t(x,y)$
and $(t,x,y) \mapsto K^{(m)}_t(x,y)$
from $(0,\infty) \times \Ri^d \times \Ri^d$ into $\Ri$
such that $K_t$ is a kernel of the operator $M_\chi \, S_t$
and $K^{(m)}_t$ is a kernel of the operator $\partial_m \, M_\chi \, S_t$
for all $t > 0$.
Moreover, 
\begin{eqnarray*}
|K_t(x,y)|
& \leq & a \, t^{-d/2} \, e^{\varepsilon t} \,
      e^{-b \, \frac{|x-y|^2}{t}} \quad \mbox{and}  \\
|K^{(m)}_t(x,y)|
& \leq & a \, t^{-d/2} \, t^{-1/2} \, e^{\varepsilon t} \,
      e^{-b \, \frac{|x-y|^2}{t}} 
\end{eqnarray*}
for all $t > 0$ and $(x,y) \in \Ri^d \times \Ri^d$.
\end{prop}
\proof\
Let $a,\omega$ be as in Lemma~\ref{llocseh311} (with $\nu = \frac{1}{2}$).
Let $C \in \cs^1(\Omega,\theta,\mu,M,{\rm real})$.
Then $M_\chi \, S_t$ is a continuous operator from $L_2(\Ri^d)$ into 
$L_\infty(\Ri^d)$ for all $t > 0$ by Lemma~\ref{locseh205}.
Moreover, $z \mapsto S_z$ is a holomorphic contraction semigroup 
on $\Sigma_{\theta_\gota}^\circ$.
Hence $M_\chi \, S_z$ maps $L_2(\Ri^d)$ continuously into 
$L_\infty(\Ri^d)$ for all $z \in \Sigma_{\theta_\gota}^\circ$ and 
$z \mapsto M_\chi S_z u$ is holomorphic for all $u \in L_2(\Ri^d)$.
It follows from the discussion after Definition~1.8 in \cite{ABu} and \cite{ABu} Theorem~3.1 
that there exists a 
measurable function $(z,x,y) \mapsto K_z(x,y)$ from 
$\Sigma_{\theta_\gota}^\circ \times \Ri^d \times \Ri^d$ into $\Ci$
such that $z \mapsto K_z(x,y)$ is holomorphic for all $x,y \in \Ri^d$ 
and $K_z$ is the kernel of $M_\chi \, S_z$ for all $z \in \Sigma_{\theta_\gota}^\circ$.
In particular, $K_t$ is a kernel of $M_\chi \, S_t$ for all $t > 0$.
By the usual minimising argument the estimates of Lemma~\ref{llocseh311}
give Gaussian bounds for the kernel $K_t$ for each $t > 0$.
Precisely, for each $t > 0$ one has $|K_t| \leq G_t$ a.e.\ on 
$\Ri^d \times \Ri^d$, where 
$G_t(x,y) = a \, t^{-d/2} \, e^{\varepsilon t} \, e^{-b \, \frac{|x-y|^2}{t}} $
for all $t > 0$ and $x,y \in \Ri^d$, where $b$ depends only on $\omega$ and~$d$.
Obviously $(t,x,y) \mapsto G_t(x,y)$ is a continuous function 
from $(0,\infty) \times \Ri^d \times \Ri^d$ into $\Ri$, therefore it is 
measurable.
Then $(-G_t) \vee K_t \wedge G_t$ is also a kernel of $M_\chi \, S_t$ for 
each $t > 0$.
Now the proposition follows for 
$M_\chi \, S_t$ with $K_t$ replaced by $(-G_t) \vee K_t \wedge G_t$.

The argument for $\partial_m \, M_\chi \, S_t$ is similar.\hfill$\Box$

\vertspace

The next lemma is also valid for complex coefficients.
The complex version will be used in the proof of Proposition~\ref{plocseh407}.

\begin{lemma} \label{llocseh308}
Let $\Omega \subsetneqq \Ri^d$ be open, $\theta \in [0,\frac{\pi}{2})$,
$\mu,M > 0$, $C \in \cs^1(\Omega,\theta,\mu,M)$ and 
$\chi \in C_{\rm b}^\infty(\Ri^d)$
with $\chi \neq 0$ and $d(\supp \chi, \Omega^{\rm c}) > 0$.
Let $\tilde \chi \in C_{\rm b}^\infty(\Ri^d)$
be such that $d(\supp \tilde \chi, \Omega^{\rm c}) > 0$ and 
$\tilde \chi(x) = 1$ for all $x \in \supp \chi$.
Then 
\begin{eqnarray}
M_\chi^2 S_t 
& = & M_\chi S_t M_\chi
   - \sum_{k,l=1}^d \int_0^t 
          M_\chi \, S_s \, M_{\partial_k \chi} \, \partial_l \, 
                M_{\tilde c_{kl}} \, M_{\widetilde \chi} \, S_{t-s} \, ds  
\label{ellocseh308;1} \\*
& & \hspace{30mm} {}
   + \sum_{k,l=1}^d \int_0^t 
           M_\chi \, S_s \, 
          \Big( M_{\partial_l \tilde c_{kl}} \, M_{\partial_k \chi} \, S_{t-s}
                - M_{c_{kl}} \, M_{\partial_k \partial_l \chi} \, S_{t-s} \Big)
          \, ds
\nonumber
\end{eqnarray}
for all $t > 0$, where the operators act in $L_2(\Ri^d)$.
\end{lemma}
\proof\
It follows from Lemma~\ref{llocseh203}\ref{llocseh203-4} that the integrals on the 
right hand side of (\ref{ellocseh308;1}) are convergent.

If $C \in \ce^1(\Omega,\theta,\mu,M)$ then
(\ref{ellocseh308;1}) follows from (\ref{ellocseh301;2}).
For all $n \in \Ni$ let $C_n = C + \frac{1}{n} I$.
Then $C_n \in \ce^1(\Omega,\theta,\mu,M+1)$.
We use (\ref{ellocseh308;1}) with respect to $C^{(n)}$.
Write $S = S^C$ and $S^{(n)} = S^{C_n}$.
It follows from Lemma~\ref{llocseh203}\ref{llocseh203-4}
that there exist $a,\omega > 0$ such that 
\[
\|M_\chi \, S^{(n)}_s \, M_{\partial_k \chi} \, \partial_l\|_{2 \to 2}
   \leq a \, s^{-1/2} \, e^{\omega s}
\]
for all $s > 0$ and $n \in \Ni$.
Obviously $\|S^{(n)}_s\|_{2 \to 2} \leq 1$ for all $s > 0$ and $n \in \Ni$.

Let $u,v \in L_2(\Ri^d)$.
Then $\lim_{n \to \infty} S^{(n)}_s u = S_s u$ in $L_2(\Ri^d)$
for all $s \in (0,t]$ by \cite{AE2} Corollary~3.9.
Note that  $M_{\widetilde \chi} \,S_s u \in W^{1,2}(\Ri^d)$ for all $s > 0$.
Let $k,l \in \{ 1,\ldots,d \} $.
Then 
\begin{eqnarray*}
\lefteqn{
|(M_\chi \, S^{(n)}_s \, M_{\partial_k \chi} \, \partial_l \, M_{\widetilde{c^{(n)}}_{kl}} \, M_{\widetilde \chi} \, S^{(n)}_{t-s} u, v)
- (M_\chi \, S_s \, M_{\partial_k \chi} \, \partial_l \, M_{\tilde c_{kl}} \, M_{\widetilde \chi} \, S_{t-s} u, v)|
} \hspace{10mm} \\*
& \leq & \tfrac{1}{n} |(M_\chi \, S^{(n)}_s \, M_{\partial_k \chi} \, \partial_l \, S^{(n)}_{t-s} u, v)|  \\*
& & \hspace{5mm} {} 
   + |(M_{\tilde c_{kl}} \, (S^{(n)}_{t-s} u - S_{t-s} u), \partial_l \, M_{\partial_k \chi} \, {S^{(n)}_s}^* \, M_\chi v)|  \\*
& & \hspace{5mm} {} 
   + |(\partial_l \, M_{\tilde c_{kl}} \, M_{\widetilde \chi} \,S_{t-s} u, 
              M_{\partial_k \chi} \, ({S^{(n)}_s}^* \, M_\chi v - S_s^* M_\chi v))|
\end{eqnarray*}
for all $s \in (0,t)$ and $n \in \Ni$.
So 
\[
\lim_{n \to \infty}
(M_\chi \, S^{(n)}_s \, M_{\partial_k \chi} \, \partial_l \, 
    M_{\widetilde{c^{(n)}}_{kl}} \, M_{\widetilde \chi} \, S^{(n)}_{t-s} u, v)
= (M_\chi \, S_s \, M_{\partial_k \chi} \, \partial_l \, M_{\tilde c_{kl}} \, M_{\widetilde \chi} \, S_{t-s} u, v)
\]
for all $s \in (0,t)$.
Moreover, 
\[
|(M_\chi \, S^{(n)}_s \, M_{\partial_k \chi} \, \partial_l \, 
     M_{\widetilde{c^{(n)}}_{kl}} \, M_{\widetilde \chi} \, S^{(n)}_{t-s} u, v)|
\leq a \, s^{-1/2} \, e^{\omega s} \, \|u\|_2 \, \|v\|_2
\]
for all $s \in (0,t)$ and $n \in \Ni$.
So 
\[
\lim_{n \to \infty}
\int_0^t (M_\chi \, S^{(n)}_s \, M_{\partial_k \chi} \, \partial_l \, 
      M_{\widetilde{c^{(n)}}_{kl}} \, M_{\widetilde \chi} \, S^{(n)}_{t-s} u, v) \, ds
= \int_0^t (M_\chi \, S_s \, M_{\partial_k \chi} \, \partial_l \, M_{\tilde c_{kl}} \, M_{\widetilde \chi} \, S_{t-s} u, v) \, ds
\]
by the Lebesgue dominated convergence theorem.

One can treat similarly the last term on the 
right hand side of (\ref{ellocseh308;1}) and the lemma follows.\hfill$\Box$

\vertspace

The next theorem is a uniform version of Theorem~\ref{tlocseh102}.

\begin{thm} \label{tlocseh309}
Let $\Omega \subsetneqq \Ri^d$ be open, $\theta \in [0,\frac{\pi}{2})$,
$\mu,M > 0$, $\nu \in (0,1)$, $\kappa > 0$, $\tau \in [0,1)$, $\varepsilon \in (0,1]$
and $\chi \in C_{\rm b}^\infty(\Ri^d)$
with $\chi \neq 0$ and $d(\supp \chi, \Omega^{\rm c}) > 0$.
Then there exist $a,b > 0$ such that for all 
$C \in \cs^1(\Omega,\theta,\mu,M,{\rm real})$
there exists a measurable function 
$(t,x,y) \mapsto K_t(x,y)$ from $(0,\infty) \times \Ri^d \times \Ri^d$ into $\Ri$ such that 
the following is valid.
\begin{tabel}
\item \label{tlocseh309-1}
For all $t > 0$ the function $K_t$ is a kernel of $M_\chi \, S_t$.
\item \label{tlocseh309-2}
For all $t > 0$ and $y \in \Ri^d$
the function $x \mapsto K_t(x,y)$ is continuously differentiable on $\Ri^d$.
\item \label{tlocseh309-3}
The function $t \mapsto (\partial_x^\alpha K_t)(x,y)$ is continuous for all $x,y \in \Ri^d$
and multi-index $\alpha$ with $|\alpha| \leq 1$.
\item \label{tlocseh309-4}
For every multi-index $\alpha$ with $|\alpha| \leq 1$ one has
\begin{equation}
|(\partial_x^\alpha K_t)(x,y)| 
   \leq a \, t^{-d/2} \, t^{-|\alpha| / 2} \, e^{\varepsilon t} \,
      e^{-b \, \frac{|x-y|^2}{t}} 
\label{etlocseh309;3}
\end{equation}
and 
\[
|(\partial_x^\alpha \, K_t)(x+h,y) 
     - (\partial_x^\alpha \, K_t)(x,y)|
\leq a \, t^{-d/2} \, t^{-|\alpha|/2} \, 
    \left( \frac{|h|}{|x-y| + \sqrt{t}} \right)^\nu
   \, e^{\varepsilon t} \, 
      e^{-b \, \frac{|x-y|^2}{t}} 
\]
for all $t > 0$ and $x,y,h \in \Ri^d$ with 
$|h| \leq \tau \, |x-y| + \kappa \, \sqrt{t}$.
\item \label{tlocseh309-5}
$ \displaystyle
(\partial^\alpha \, M_\chi \, S_t u)(x)
= \int_{\Ri^d} (\partial_x^\alpha K_t)(x,y) \, u(y) \, dy
$
for all $t > 0$, $u \in L_1(\Ri^d)$, $x \in \Ri^d$ and multi-index $\alpha$
with $|\alpha| \leq 1$.
\end{tabel}
\end{thm}
\proof\
Without loss of generality we may prove the proposition with 
$M_\chi \, S_t$ replaced by $M_\chi^2 \, S_t$.
We consider each of the operators in the terms on the right hand side 
of (\ref{ellocseh308;1}).
It is possible to differentiate all terms in (\ref{ellocseh308;1}) at least once
in $L_2$-sense.
Let $\alpha$ be a multi-index with $|\alpha| \leq 1$.
Clearly $\partial^\alpha \, M_\chi S_t M_\chi$ has a kernel satisfying the stated 
requirements by Theorem~\ref{tlocseh201}.
Let $k,l \in \{ 1,\ldots,d \} $.
For all $t > 0$ consider the operator
\begin{equation}
\int_0^t \partial^\alpha \, M_\chi \, S_s \, M_{\partial_k \chi} \, \partial_l \, 
                M_{\tilde c_{kl}} \, M_{\widetilde \chi} \, S_{t-s} \, ds 
.  
\label{etlocseh309;1}
\end{equation}
By Theorem~\ref{tlocseh201} 
there exist a continuous function $(t,x,y) \mapsto L^{(1,\alpha)}_t(x,y)$
from $(0,\infty) \times \Ri^d \times \Ri^d$ into $\Ri$ and suitable constants
$a_1,b_1 > 0$ 
such that for all $t > 0$ the function $L^{(1,\alpha)}_t$ is a 
kernel of the operator $\partial^\alpha \, M_\chi \, S_t \, M_{\partial_k \chi}$,
\begin{eqnarray*}
|L^{(1,\alpha)}_t(x,y)| 
   & \leq & a_1 \, t^{-d/2} \, t^{-|\alpha|/2} \, e^{\varepsilon t} \,
      e^{-b_1 \, \frac{|x-y|^2}{t}} \mbox{ and}  \\
|L^{(1,\alpha)}_t(x+h,y) - L^{(1,\alpha)}_t(x,y)| 
   & \leq & a_1 \, t^{-d/2} \, t^{-|\alpha|/2} \, 
    \left( \frac{|h|}{\sqrt{t}} \right)^\nu
   \, e^{\varepsilon t} \, 
      e^{-b_1 \, \frac{|x-y|^2}{t}} 
\end{eqnarray*}
for all $t > 0$ and $x,y,h \in \Ri^d$ with 
$|h| \leq \frac{1}{2} |x-y| + \sqrt{t}$.
Next, by Proposition~\ref{plocseh310} 
there exist a measurable function $(t,x,y) \mapsto L^{(2)}_t(x,y)$
from $(0,\infty) \times \Ri^d \times \Ri^d$ into $\Ri$ and suitable constants
$a_2,b_2 > 0$ 
such that for all $t > 0$ the function $L^{(2)}_t$ is a 
kernel of the operator $\partial_l \, M_{\tilde c_{kl}} \, M_{\widetilde \chi} \, S_t$
and
\[
|L^{(2)}_t(x,y)| 
   \leq a_2 \, t^{-d/2} \, t^{-1/2} \, e^{\varepsilon t} \,
      e^{-b_2 \, \frac{|x-y|^2}{t}} 
\]
for all $t > 0$ and $x,y \in \Ri^d$.
Define the function $(t,x,y) \mapsto L^{(3,\alpha)}_t(x,y)$
from $(0,\infty) \times \Ri^d \times \Ri^d$ into $\Ri$ by
\[
L^{(3,\alpha)}_t(x,y)
= \int_0^t \!\int_{\Ri^d} L^{(1,\alpha)}_s(x,z) \, L^{(2)}_{t-s}(z,y) \, dz \, ds
.  \]
Since the convolution of two Gaussians is a Gaussian, it follows that 
the integral is convergent and $L^{(3,\alpha)}$ has appropriate Gaussian bounds.
Moreover, $(t,x,y) \mapsto L^{(3,\alpha)}_t(x,y)$ is measurable and 
$L^{(3,\alpha)}_t$ is a kernel of the operator (\ref{etlocseh309;1})
for all $t > 0$.
Define $L^{(3)}_t = L^{(3,\beta)}$ if $|\beta| = 0$.
It is an elementary exercise in integration theory to prove that 
$x \mapsto L^{(3,\alpha)}_t(x,y)$ is continuous for all $t > 0$ and $y \in \Ri^d$,
that $t \mapsto L^{(3,\alpha)}_t(x,y)$ is continuous for all $x,y \in \Ri^d$
and that for all $t > 0$ and $y \in \Ri^d$
the function $x \mapsto L^{(3)}_t(x,y)$ is differentiable 
and $(\partial_x^\alpha L^{(3)}_t)(x,y) = L^{(3,\alpha)}_t(x,y)$ for all $x \in \Ri^d$.

The last term in (\ref{ellocseh308;1}) can be treated in a 
similar way and Statements~\ref{tlocseh309-1}--\ref{tlocseh309-3} 
and the first part of Statement~\ref{tlocseh309-4} follow.
Let $K$ be the so obtained kernel.

Next, let $t > 0$ and $u \in L_1(\Ri^d)$.
Since $(x,y) \mapsto (\partial_x^\alpha K_t)(x,y)$ is a kernel of 
$\partial^\alpha \, M_\chi \, S_t$ it follows that 
$(\partial^\alpha \, M_\chi \, S_t u)(x) 
   = \int_{\Ri^d} (\partial_x^\alpha K_t)(x,y) \, u(y) \, dy$ 
for a.e.\ $x \in \Ri^d$.
Hence this is valid for all $x \in\Ri^d$ since 
$\partial^\alpha \, M_\chi \, S_t u$ is continuous by Lemmas~\ref{llocseh311}
and \ref{llocseh313}.
This proves Statement~\ref{tlocseh309-5}.

Finally we the H\"older Gaussian bounds of Statement~\ref{tlocseh309-4}.
Set $\tilde \nu = \frac{1}{2}(1 + \nu)$.
By Lemma~\ref{llocseh313} there exists an $a_1 > 0$, depending only 
on $\Omega$, $\theta$, $\mu$, $M$, $\nu$, $\varepsilon$ and $\chi$, such that 
\[
\|(I - L(h)) \partial_m \, M_\chi \, S_t\|_{1 \to \infty}
\leq a_1 \, t^{-d/2} \, t^{-1/2} \, 
    \left( \frac{|h|}{\sqrt{t}} \right)^{\tilde \nu} \, e^{\varepsilon t} 
\]
for all $t > 0$, $m \in \{ 1,\ldots,d \} $ and $h \in \Ri^d$.
Let $x,h \in \Ri^d$, $t > 0$ and  $m \in \{ 1,\ldots,d \} $.
Then it follows from Statement~\ref{tlocseh309-5} that 
\begin{eqnarray*}
\Big| \int \Big( (\partial_{x,m} K_t)(x,y) - (\partial_{x,m} K_t)(x-h,y) \Big) u(y) \, dy \Big|
& = & \Big| \Big( (I - L(h)) \partial_m \, M_\chi\, S_t u \Big) (x) \Big|  \\
& \leq & a_1 \, t^{-d/2} \, t^{-1/2} \, 
    \left( \frac{|h|}{\sqrt{t}} \right)^{\tilde \nu} \, e^{\varepsilon t} \, \|u\|_1
\end{eqnarray*}
for all $u \in L_1(\Ri^d)$.
Hence there exists a null set $N_{x,h,t,m} \subset \Ri^d$ such that 
\begin{equation}
|(\partial_{x,m} \, K_t)(x,y) 
     - (\partial_{x,m} \, K_t)(x-h,y)|
\leq a_1 \, t^{-d/2} \, t^{-1/2} \, 
    \left( \frac{|h|}{\sqrt{t}} \right)^{\tilde \nu} \, e^{\varepsilon t}
\label{etlocseh309;2}
\end{equation}
for all $y \in \Ri^d \setminus N_{x,h,t,m}$.
Then 
\[
N = \bigcup_{x,h \in \Qi^d} \:
    \bigcup_{t \in (0,\infty) \cap \Qi} \:
    \bigcup_{m=1}^d
      N_{x,h,t,m}
\]
is a null set in $\Ri^d$ and (\ref{etlocseh309;2})
is valid for all $y \in \Ri^d \setminus N$, $x,h \in \Qi^d$, $t \in (0,\infty) \cap \Qi$
and $m \in \{ 1,\ldots,d \} $.
Then by density and continuity (Statements~\ref{tlocseh309-2} and \ref{tlocseh309-3})
one deduces that (\ref{etlocseh309;2}) is valid for all $y \in \Ri^d \setminus N$, 
$x,h \in \Ri^d$, $t > 0$ and  $m \in \{ 1,\ldots,d \} $.

Let $a,b$ be as in (\ref{etlocseh309;3}).
Then it follows as in the proof of Step~\ref{locsehAStep7} of Lemma~\ref{llocsehA1} that 
 there are $a_2,b_2 > 0$, depending only on 
$a$, $b$, $a_1$, $\kappa$, $\tau$, $\varepsilon$, $\nu$ and $\tilde \nu$,
such that 
\[
|(\partial_{x,m} \, K_t)(x,y) 
     - (\partial_{x,m} \, K_t)(x-h,y)|
\leq a_2 \, t^{-d/2} \, t^{-|\alpha|/2} \, 
    \left( \frac{|h|}{|x-y| + \sqrt{t}} \right)^\nu
   \, e^{\varepsilon t} \, 
      e^{-b_2 \, \frac{|x-y|^2}{t}} 
\]
for all $y \in \Ri^d \setminus N$, 
$x,h \in \Ri^d$, $t > 0$ and  $m \in \{ 1,\ldots,d \} $
with $|h| \leq \tau \, |x-y| + \kappa \, \sqrt{t}$.
Now the theorem follows by replacing 
$(t,x,y) \mapsto K_t(x,y)$ by the function
$(t,x,y) \mapsto K_t(x,y) \, \one_{\Ri^d \setminus N}(y)$.\hfill$\Box$

\begin{cor} \label{clocseh314}
Let $\Omega \subsetneqq \Ri^d$ be open, $\theta \in [0,\frac{\pi}{2})$,
$\mu,M > 0$ and 
$\chi \in C_{\rm b}^\infty(\Ri^d)$
with $\chi \neq 0$ and $d(\supp \chi, \Omega^{\rm c}) > 0$.
Then there exists a $c > 0$ such that 
$M_\chi \, (I + A)^{-1} \, L_p(\Ri^d) \subset W^{1,p}(\Ri^d)$ and 
\[
\|\partial_m \, M_\chi u\|_p
\leq c \, \|(I + A) u\|_p
\]
for all $C \in \cs^1(\Omega,\theta,\mu,M,{\rm real})$, $p \in [1,\infty]$,
$u \in D(A)$ and $m \in \{ 1,\ldots,d \} $.
\end{cor}
\proof\
The Gaussian bounds of Theorem~\ref{tlocseh309} imply bounds
$\|\partial_m \, M_\chi \, S_t\|_{p \to p} \leq c \, t^{-1/2} \, e^{t/2}$.
Then the corollary follows by a Laplace transform.\hfill$\Box$

\section{Riesz transforms} \label{Slocseh4}

In this section we shall prove that various Riesz transforms like 
$\partial_k \, M_\chi (I + A)^{-1/2}$ and 
$\partial_k \, M_\chi (I + A)^{-1/2} \, M_\chi $ are bounded on $L_2$ or $L_p$.
The first results on $L_2$ for the Riesz transform merely use
that $\chi \in W^{1,\infty}(\Ri^d)$.

If $\eta \in W^{1,\infty}(\Ri^d)$, then 
\[
\RRe \gota_C(\eta \, u)
\leq 2 \|\eta\|_\infty \RRe \gota_C(u) + 2 M \, \|\nabla \eta\|_\infty^2 \, \|u\|_2^2
\]
for all $u \in W^{1,2}(\Ri^d)$, $\theta \in [0,\frac{\pi}{2})$, $M > 0$
and $C \in \cs(\theta,M)$.
Therefore for self-adjoint operators the boundedness of the 
Riesz transforms on $L_2(\Ri^d)$ is trivial.
Throughout this section let $L = - \sum_{k=1}^d \partial_k^2$ be the Laplacian
and let $H = A_{\Re C}$.

The first lemma is a variation of Lemma~\ref{llocseh203}\ref{llocseh203-3}.

\begin{lemma} \label{llocseh402}
Let $\theta \in [0,\frac{\pi}{2})$, $\mu,M > 0$, 
$\chi \in W^{1,\infty}(\Ri^d)$ and 
$C \in \cs(\supp \chi,\theta,\mu,M)$.
Then $M_\chi D(H^{1/2}) \subset W^{1/2}(\Ri^d)$ and 
\[
\|\partial_k \, M_\chi u\|_2
\leq \|(I + L)^{1/2} \, M_\chi u\|_2
\leq c_1 \, \|(I + H)^{1/2} u\|_2
\]
for all $u \in D(H^{1/2})$ and $k \in \{ 1,\ldots,d \} $,
where 
\[
c_1 = \Big( \|\chi\|_\infty^2 + 2 \frac{\|\chi\|_\infty^2}{\mu} 
            + 2 \|\nabla \chi\|_\infty^2
      \Big)^{1/2}
.  \]
\end{lemma}
\proof\
We only have to prove the last estimate.
It follows from Lemma~\ref{llocseh203}\ref{llocseh203-3} that 
\begin{eqnarray*}
\|(I + L)^{1/2} \, M_\chi u\|_2^2
& = & \|M_\chi u\|_2^2 + \sum_{k=1}^d \|\partial_k \, M_\chi u\|_2^2  \\
& \leq & \|\chi\|_\infty^2 \, \|u\|_2^2
   + 2 \sum_{k=1}^d \|(\partial_k \, M_\chi - M_{\partial_k \chi}) u\|_2^2
   + 2 \sum_{k=1}^d \|M_{\partial_k \chi} u\|_2^2  \\
& \leq & \|\chi\|_\infty^2 \, \|u\|_2^2
   + 2 \frac{\|\chi\|_\infty^2}{\mu} \, \|H^{1/2} u\|_2^2
   + 2 \|(\nabla \chi) u\|_2^2  \\
& \leq & c_1^2 \, \|(I + H)^{1/2} u\|_2^2
\end{eqnarray*}
as required.\hfill$\Box$

\begin{lemma} \label{llocseh402.5}
Let $\theta \in [0,\frac{\pi}{2})$, $\mu, M > 0$,
$\gamma \in (0,1)$
and $\chi \in W^{1,\infty}(\Ri^d)$.
Then there exists a $c_2 > 0$ such that 
$M_\chi D(A^{\gamma/2}) \subset D(L^{\gamma/2})$ and 
\[
\|(I + L)^{\gamma/2} \, M_\chi u\|_2
\leq c_2 \, \|(I + A)^{\gamma/2} u\|_2
\]
for all $C \in \cs(\supp \chi,\theta,\mu,M)$ 
and $u \in D((I + A)^{\gamma/2})$. 
\end{lemma}
\proof\
Let $c_1$ be as in Lemma~\ref{llocseh402}.
Then $M_\chi$ is continuous from $D((I + H)^{1/2})$ into $D((I + L)^{1/2})$
with norm bounded by $c_1$.
Then by interpolation, Proposition~G together with Theorem~G in \cite{ADM},
it follows that $M_\chi D((I + H)^{\gamma/2}) \subset D((I + L)^{\gamma/2})$
and 
\[
\|(I + L)^{\gamma/2} \, M_\chi u\|_2
\leq c_1' \, \|(I + H)^{\gamma/2} u\|_2
\]
for all $u \in D((I + L)^{\gamma/2})$, where
$c_1' = c_1^\gamma \, \|\chi\|_\infty^{1-\gamma}$.
But $D((I + H)^{\gamma/2}) = D((I + A)^{\gamma/2})$ with equivalent norms 
by \cite{Kat2} Theorem~3.1.
Explicitly, 
\[
\|(I + H)^{\gamma/2} u\|_2 
\leq \frac{1}{1 - \tan \frac{\pi \, \gamma}{4}} \, \|(I + A)^{\gamma/2} u\|_2 
\]
for all $u \in D((I + A)^{\gamma/2})$.
Then the lemma follows.\hfill$\Box$

\begin{lemma} \label{llocseh403}
Let $\theta \in [0,\frac{\pi}{2})$, $\mu, M > 0$,
$\chi \in W^{1,\infty}(\Ri^d)$ and $C \in \ce(\supp \chi,\theta,\mu,M)$.
Then 
\[
\|(I + H)^{1/2} M_\chi u\|_2
\leq 2 \|(I + H)^{-1/2} M_\chi \, (I + A) u\|_2
   + 4 M^{1/2} \, \|\nabla \chi\|_\infty \, \|u\|_2
\]
for all $u \in D(A)$.
\end{lemma}
\proof\
The proof is a variation of the proof of Lemma~1 in \cite{Kat3}.
If $u \in D(A)$ then 
\begin{eqnarray*}
\|(I + H)^{1/2} M_\chi u\|_2^2
& = & \RRe \gota(M_\chi u) + \|M_\chi u\|_2^2  \\
& = & \RRe \gota(u, M_\chi^2 u) 
   - \RRe \sum (c_{kl} \, \partial_k(\chi \, u) , (\partial_l \chi) \, u)  \\*
& & {} \hspace*{0mm}
   + \RRe \sum (c_{kl} \, (\partial_k \chi) \, u , (\partial_l \chi) \, u)
   + \RRe \sum (c_{kl} \, (\partial_k \chi) \, u , \partial_l (\chi \, u))
   + \|M_\chi u\|_2^2  
.  
\end{eqnarray*}
But 
\begin{eqnarray*}
\RRe \gota(u, M_\chi^2 u) + \|M_\chi u\|_2^2  
& = & \RRe (A u, M_\chi^2 u) + \|M_\chi u\|_2^2  \\
& = & \RRe ((I + H)^{-1/2} M_\chi \, (I + A) u, (I + H)^{1/2} M_\chi u)  \\
& \leq & \|(I + H)^{-1/2} M_\chi \, (I + A) u\|_2 \, \|(I + H)^{1/2} M_\chi u\|_2  \\
& \leq & \tfrac{1}{4} \,  \|(I + H)^{1/2} M_\chi u\|_2^2 + \|(I + H)^{-1/2} M_\chi \, (I + A) u\|_2^2
\end{eqnarray*}
and 
\begin{eqnarray*}
\lefteqn{
- \RRe \sum (c_{kl} \, \partial_k(\chi \, u) , (\partial_l \chi) \, u) 
   + \RRe \sum (c_{kl} \, (\partial_k \chi) \, u , (\partial_l \chi) \, u)
   + \RRe \sum (c_{kl} \, (\partial_k \chi) \, u , \partial_l (\chi \, u))
} \hspace{10mm} \\*
& \leq & 2 (\RRe \gota(M_\chi u) )^{1/2} 
    \Big( \RRe \sum \int c_{kl} \, (\partial_k \chi) \, (\partial_l \chi) \, |u|^2 \Big)^{1/2} 
   + \RRe \sum \int c_{kl} \, (\partial_k \chi) \, (\partial_l \chi) \, |u|^2  \\
& \leq & 2 (\RRe \gota(M_\chi u) )^{1/2} \, M^{1/2} \, \|\nabla \chi\|_\infty \, \|u\|_2
   + M \, \|\nabla \chi\|_\infty^2 \, \|u\|_2^2  \\
& \leq & \tfrac{1}{4} \|(I + H)^{1/2} M_\chi u\|_2^2
   + 5 M \, \|\nabla \chi\|_\infty^2 \, \|u\|_2^2
.
\end{eqnarray*}
So 
\[
\|(I + H)^{1/2} M_\chi u\|_2^2
\leq 2 \|(I + H)^{-1/2} M_\chi \, (I + A) u\|_2^2
   + 10 M \, \|\nabla \chi\|_\infty^2 \, \|u\|_2^2
\]
and the lemma follows.\hfill$\Box$

\vertspace

The next lemma is well known, but we need uniform constants.

\begin{lemma} \label{llocseh404}
Let $0 < \gamma < \nu < 1$ and $\eta \in C^{0,\nu}(\Ri^d) \cap L_\infty(\Ri^d)$.
Let $u \in D(L^{\gamma/2})$.
Then $M_\eta \, u \in D(L^{\gamma/2})$ and 
\[
\|(I + L)^{\gamma/2} \, M_\eta u\|_2
\leq c_4 \, (\|\eta\|_\infty \vee |||\eta|||_{C^{0,\nu}}) \, \|(I + L)^{\gamma/2} u\|_2
,  \]
where $c_4 = 1 + c_\gamma^{-1} \, 2^{d+3} \, \Gamma(\frac{\nu - \gamma}{2})$
and $c_\gamma = \int_0^\infty t^{-1 - \frac{\gamma}{2}} (1 - e^{-t}) \, dt$.
\end{lemma}
\proof\
Let $T$ be the semigroup generated by $-L$ and for all $t > 0$ 
let $G_t$ be the kernel of $T_t$.
If $u \in L_2(\Ri^d)$ and $t > 0$ then 
\[
([T_t, M_\eta] u)(x)
= \int G_t(y) \, \Big( \eta(x-y) - \eta(x) \Big) u(x-y) \, dy
\]
for all $x \in \Ri^d$.
But
\[
\Big|G_t(y) \, \Big( \eta(x-y) - \eta(x) \Big) \Big|
\leq 2 M \, G_t(y) \, |y|^\nu
\leq 2^{d+3} M \, t^{\nu/2} \, G_{2t}(y)
\]
for all $x,y \in \Ri^d$, 
where $M = \|\eta\|_\infty \vee |||\eta|||_{C^{0,\nu}}$.
Therefore $\|[T_t, M_\eta] u\|_2 \leq 2^{d+3} M \, t^{\nu/2} \, \|u\|_2$.
Next, 
\[
(I + L)^{\gamma/2} 
= \frac{1}{c_\gamma} \, \int_0^\infty t^{-1 - \frac{\gamma}{2}} ( I - e^{-t} \, T_t ) \, dt
.  \]
So for all $u \in C_c^\infty(\Ri^d)$ one obtains
\[
[(I + L)^{\gamma/2}, M_\eta] u
= - \frac{1}{c_\gamma} \, \int_0^\infty t^{-1 - \frac{\gamma}{2}} \, e^{-t} [T_t, M_\eta] u \, dt
.  \]
Therefore 
\[
\|[(I + L)^{\gamma/2}, M_\eta] u\|_2
\leq \frac{2^{d+3} M}{c_\gamma} \int_0^\infty t^{-1 + \frac{\nu - \gamma}{2}} \, e^{-t} \, dt \, \|u\|_2
= \frac{2^{d+3} M \, \Gamma(\frac{\nu - \gamma}{2})}{c_\gamma} \, \|u\|_2
\]
and 
\[
\|(I + L)^{\gamma/2} \, M_\eta u\|_2
\leq c_4 \, \|(I + L)^{\gamma/2} u\|_2
.  \]
Then by density the lemma follows.\hfill$\Box$

\begin{lemma} \label{llocseh405}
Let $\theta \in [0,\frac{\pi}{2})$, $\mu,M > 0$,
$\nu \in (0,1)$
and $\chi \in W^{1,\infty}(\Ri^d)$.
Set $\gamma = \frac{\nu}{2}$.
Then there exists a $c_5 > 0$ such that 
$M_\chi^2 D(L^{(1+\gamma)/2}) \subset D(A^{(1+\gamma)/2})$ and 
\[
\|(I + A)^{(1+\gamma)/2} \, M_\chi^2 u\|_2
\leq c_5 \, \|(I + L)^{(1+\gamma)/2} u\|_2
\]
for all $C \in \ce^\nu(\supp \chi,\theta,\mu,M)$ and $u \in D(L^{(1+\gamma)/2})$.
\end{lemma}
\proof\
The proof is a variation of the proof in \cite{ER16}.
Let $u \in D(L^{(1+\gamma)/2}) \subset W^{1,2}(\Ri^d)$ and $v \in D(A^*) \subset W^{1,2}(\Ri^d)$.
Then 
\begin{eqnarray}
\lefteqn{
(M_\chi^2 u, (I + A^*)^{(1+\gamma)/2} v)
} \hspace{10mm} \nonumber  \\*
& = & (M_\chi^2 u, (I + A^*) (I + A^*)^{-(1-\gamma)/2} v)  \nonumber  \\
& = & \sum (c_{kl} \, \partial_k \, M_\chi^2 u, \partial_l \, (I + A^*)^{-(1-\gamma)/2} v)
   + (M_\chi^2 u, (I + A^*)^{-(1-\gamma)/2} v) \nonumber  \\
& = & \sum (c_{kl} \, M_\chi \, \partial_k u, \partial_l \, M_\chi \, (I + A^*)^{-(1-\gamma)/2} v)
\label{ellocseh405-1}
   \\*
& & \hspace{10mm} {}
   - \sum (c_{kl} \, M_\chi \, \partial_k u, (\partial_l \, \chi) \, (I + A^*)^{-(1-\gamma)/2} v)
   \nonumber  \\*
& & \hspace{10mm} {}
   + 2 \sum (c_{kl} \, (\partial_k \chi) \, u, \partial_l \, M_\chi \, (I + A^*)^{-(1-\gamma)/2} v)
   \nonumber  \\*
& & \hspace{10mm} {}
   - 2 \sum (c_{kl} \, (\partial_k \chi) \, u, (\partial_l \, \chi) \, (I + A^*)^{-(1-\gamma)/2} v)
   \nonumber  \\*
& & \hspace{10mm} {}
   + (M_\chi u, M_\chi \, (I + A^*)^{-(1-\gamma)/2} v) \nonumber 
.
\end{eqnarray}
Fix $k,l \in \{ 1,\ldots,d \} $.
Then $\chi \, c_{kl} \in C^{0,\nu}(\Ri^d) \cap L_\infty(\Ri^d)$
and $|||\chi \, c_{kl}|||_{C^{0,\nu}(\Ri^d)} \leq 2 M \, \|\chi\|_{W^{1,\infty}}$.
So 
\begin{eqnarray*}
\lefteqn{
|(c_{kl} \, M_\chi \, \partial_k u, \partial_l \, M_\chi \, (I + A^*)^{-(1-\gamma)/2} v)|
} \hspace{10mm} \\*
& = & |((I + L)^{\gamma/2} \, M_{\chi \, c_{kl}} \, \partial_k u, 
   (I + L)^{-\gamma/2} \, \partial_l \, M_\chi \, (I + A^*)^{-(1-\gamma)/2} v)|  \\
& \leq & \|(I + L)^{\gamma/2} \, M_{\chi \, c_{kl}} \, \partial_k u\|_2 \,
   \|(I + L)^{-\gamma/2} \, \partial_l \, M_\chi \, (I + A^*)^{-(1-\gamma)/2} v\|_2
.
\end{eqnarray*}
If $c_4$ is as in Lemma~\ref{llocseh404}, then 
\begin{eqnarray*}
\|(I + L)^{\gamma/2} \, M_{\chi \, c_{kl}} \, \partial_k u\|_2
& \leq & 2 c_4 \, M \, \, \|\chi\|_{W^{1,\infty}} \, \|(I + L)^{\gamma/2} \, \partial_k u\|_2  \\
& \leq & 2 c_4 \, M \, \, \|\chi\|_{W^{1,\infty}} \, \|(I + L)^{(1+\gamma)/2} u\|_2  
.
\end{eqnarray*}
Alternatively, 
\[
\|(I + L)^{-\gamma/2} \, \partial_l \, M_\chi \, (I + A^*)^{-(1-\gamma)/2} v\|_2
\leq \|(I + L)^{(1-\gamma)/2} \, M_\chi \, (I + A^*)^{-(1-\gamma)/2} v\|_2
.  \]
By Lemma~\ref{llocseh402.5} there exists a $c_2 > 0$, 
depending only on $\theta$, $\mu$, $M$, $\nu$ and $\chi$,
such that 
\[
\|(I + L)^{(1-\gamma)/2} \, M_\chi w\|_2
\leq c_2 \|(I + A^*)^{(1-\gamma)/2} w\|_2
\]
for all $w \in D((I + L)^{(1-\gamma)/2})$.
Hence
\[
\|(I + L)^{-\gamma/2} \, \partial_l \, M_\chi \, (I + A^*)^{-(1-\gamma)/2} v\|_2
\leq c_2 \, \|(I + A^*)^{(1-\gamma)/2} \, (I + A^*)^{-(1-\gamma)/2} v\|_2
= c_2 \, \|v\|_2
.   \]
The other four terms in (\ref{ellocseh405-1}) can be estimated similarly.

Combining the contributions, it follows that there exists a $c > 0$,
depending only on $\theta$, $\mu$, $M$, $\nu$ and $\chi$, such that 
\[
|(M_\chi^2 u, (I + A^*)^{(1+\gamma)/2} v)|
\leq c \, \|(I + L)^{(1+\gamma)/2} u\|_2 \, \, \|v\|_2
\]
for all $v \in D(A^*)$.
Since $D(A^*)$ is a core for 
$(I + A^*)^{(1+\gamma)/2}$, one deduces that $M_\chi^2 u \in D(A^{(1+\gamma)/2})$
and 
\[
\|(I + A)^{(1+\gamma)/2} \, M_\chi^2 u\|_2
\leq c \, \|(I + L)^{(1+\gamma)/2} u\|_2
\]
as required.\hfill$\Box$

\vertspace

Now we are able to prove a uniform version of Theorem~\ref{tlocseh103}\ref{tlocseh103-3}.

\begin{thm} \label{tlocseh406}
Let $\theta \in [0,\frac{\pi}{2})$, $\mu,M > 0$,
$\nu \in (0,1)$
and $\chi \in W^{1,\infty}(\Ri^d)$.
Then there exists a $c > 0$ such that 
$M_\chi^3 u \in W^{1,2}(\Ri^d)$ and 
\[
\|\partial_k \, M_\chi^3 u\|_2
\leq c \, \|(I + A)^{1/2} u\|_2
\]
for all $C \in \cs^\nu(\supp \chi,\theta,\mu,M)$,  $u \in D(A^{1/2})$
and $k \in \{ 1,\ldots,d \} $.
\end{thm}
\proof\
First suppose that $C \in \ce^\nu(\Omega,\theta,\mu,M)$.
Let $c_5 > 0$ be as in Lemma~\ref{llocseh405}.
Then by interpolation, Proposition~G together with Theorem~G in \cite{ADM},
one establishes that $M_\chi^2 D((I + L)^{1/2}) \subset D((I + A)^{1/2})$ and 
\[
\|(I + A)^{1/2} \, M_\chi^2 \, u\|_2
\leq c_5' \,  \|(I + L)^{1/2} u\|_2
\]
for all $u \in W^{1,2}(\Ri^d) = D(L^{1/2})$, 
where $c_5' = c_5^{\frac{1}{1+\gamma}} \, \|\chi\|_\infty^{\frac{\gamma}{1+\gamma}}$.
Hence if $c_1 > 0$ is as in Lemma~\ref{llocseh402} then 
\[
\|(I + A)^{1/2} \, M_\chi^3 \, u\|_2
\leq c_5' \|(I + L)^{1/2} \, M_\chi \, u\|_2
\leq c_1 \, c_5' \, \|(I + H)^{1/2} u\|_2
.  \]
So 
\[
\|(I + A)^{1/2} \, M_\chi^3 \, (I + H)^{-1/2}\|_{2 \to 2} 
\leq c_1 \, c_5'
\]
and then by duality
\begin{equation}
\|(I + H)^{-1/2} \, M_\chi^3 \, (I + A^*)^{1/2}\|_{2 \to 2} 
\leq c_1 \, c_5'
. 
\label{etlocseh406;1}
\end{equation}
Since $\ce^\nu(\Omega,\theta,\mu,M)$ is invariant under taking adjoints,
one may replace $A^*$ by $A$ in (\ref{etlocseh406;1}).
Then Lemma~\ref{llocseh403} gives
\begin{eqnarray*}
\|(I + H)^{1/2} \, M_\chi^3 u\|_2
& \leq & 2 \|(I + H)^{-1/2} \, M_\chi^3 \, (I + A) u\|_2
   + c_6 \, \|u\|_2  \\
& \leq & 2 c_1 \, c_5' \, \|(I + A)^{1/2} u\|_2 + c_6 \, \|u\|_2  \\
& \leq & (2 c_1 \, c_5' + c_6) \, \|(I + A)^{1/2} u\|_2
\end{eqnarray*}
for all $u \in D(A)$, where $c_6 = 4 M \, \|\nabla(\chi^3)\|_\infty$.
Therefore
\begin{eqnarray*}
\|\partial_k \, M_\chi^3 u\|_2
& \leq & \mu^{-1/2} \, \|(I + H)^{1/2} \, M_\chi^3 u\|_2
\leq \mu^{-1/2} \, (2 c_1 \, c_5' + c_6) \, \|(I + A)^{1/2} u\|_2
.
\end{eqnarray*}
This extends to all $u \in D(A^{1/2})$ by density.
It follows that 
\begin{equation}
|(M_\chi^3 \, (I + A)^{-1/2} u, \partial_k v)|
\leq \mu^{-1/2} \, (2 c_1 \, c_5' + c_6) \, \|u\|_2 \, \|v\|_2
\label{etlocseh406;2}
\end{equation}
for all $u \in L_2(\Ri^d)$, $v \in W^{1,2}(\Ri^d)$, $k \in \{ 1,\ldots,d \} $
and $C \in \ce^\nu(\supp \chi,\theta,\mu,M)$.
By approximating $C$ by $C + \frac{1}{n} \, I$ it follows as before that 
(\ref{etlocseh406;2}) extends to all $C \in \cs^\nu(\supp \chi,\theta,\mu,M)$
and the theorem follows.\hfill$\Box$

\vertspace

The theorem has many corollaries.

\begin{cor} \label{clocseh402}
Let $\Omega \subsetneqq \Ri^d$ be open, $\theta \in [0,\frac{\pi}{2})$,
$\mu,M > 0$, $\nu \in (0,1)$, $p \in (1,\infty)$ and $\chi \in C_{\rm b}^\infty(\Ri^d)$
with $\chi \neq 0$ and $d(\supp \chi, \Omega^{\rm c}) > 0$.
Then there exists a $c > 0$ such that 
\[
\|\nabla \, M_\chi \, (I + A)^{-1/2} \, M_\chi\|_{p \to p} \leq c
\]
for all $C \in \cs^\nu(\Omega,\theta,\mu,M)$.
\end{cor}
\proof\
For all $t > 0$ let $K_t$ be the continuous kernel of the operator $M_\chi \, S_t \, M_\chi$.
Fix $m \in \{ 1,\ldots,d \} $.
By Theorem~\ref{tlocseh201} there are $a > 0$ and $b \in (0,1)$ such that 
\[
|(\partial_{x,m} K_t)(x,y)| 
   \leq a \, t^{-d/2} \, t^{-1 / 2} \, e^{t/2} \,
      e^{-b \, \frac{|x-y|^2}{t}} 
\]
and 
\[
|(\partial_{x,m} \, K_t)(x+h,y+k) 
     - (\partial_{x,m} \, K_t)(x,y)|
\leq a \, t^{-d/2} \, t^{-|\alpha|/2} \, 
    \left( \frac{|h| + |k|}{\sqrt{t}} \right)^\nu
   \, e^t \, 
      e^{-b \, \frac{|x-y|^2}{t}} 
\]
for all $t > 0$ and $x,y,h,k \in \Ri^d$ with 
$|h| + |k| \leq \frac{1}{2} \, |x-y|$.
For all $x,y \in \Ri^d$ with $x \neq y$ define 
\[
L(x,y) 
= \frac{1}{\sqrt{\pi}} \, \int_0^\infty t^{-1/2} \, e^{-t} \, 
    (\partial_{x,m} K_t)(x,y) \, dt
.  \]
Then 
\[
(\partial_m \, M_\chi \, (I + A)^{-1/2} \, M_\chi u)(x)
= \int_{\Ri^d} L(x,y) \, u(y) \, dy
\]
for all $u \in L_2(\Ri^d)$ and $x \in \Ri^d \setminus \supp u$.
Moreover, for all $x,y,h \in \Ri^d$ with $|h| \leq \frac{1}{2} \, |x-y|$ one has
\begin{eqnarray*}
|L(x+h,y) - L(x,y)|
& \leq & \frac{a}{\sqrt{\pi}} \, \int_0^\infty t^{-d/2} \, t^{-1} \, 
    \left( \frac{|h|}{\sqrt{t}} \right)^\nu \, e^{-b \, \frac{|x-y|^2}{t}} \, dt  
= c  \, \left( \frac{|h|}{|x-y|} \right)^\nu \, \frac{1}{|x-y|^d} ,
\end{eqnarray*}
where $c = \frac{a}{\sqrt{\pi}}
            \int_0^\infty t^{-1 - \frac{d+\nu}{2}} \, e^{- \frac{b}{t}} \, dt$.
Similarly, 
\[
|L(x,y+k) - L(x,y)|
\leq c \, \left( \frac{|h|}{|x-y|} \right)^\nu \, \frac{1}{|x-y|^d}
\]
for all $x,y,k \in \Ri^d$ with $|k| \leq \frac{1}{2} \, |x-y|$.
In addition,
\[
|L(x,y)|
\leq \frac{a}{\sqrt{\pi}} \, \frac{1}{|x-y|^d} \, 
     \int_0^\infty t^{-1 - \frac{d+\nu}{2}} \, e^{- \frac{b}{t}} \, dt
\]
for all $x,y \in \Ri^d$ with $x \neq y$.
Therefore $\partial_m \, M_\chi \, (I + A)^{-1/2} \, M_\chi$ is 
a Calder\'on--Zygmund operator.
Since it is bounded on $L_2(\Ri^d)$ by Theorem~\ref{tlocseh406},
it is also bounded in $L_p(\Ri^d)$ for all $p\in (1,\infty)$
by Theorem~L in \cite{ADM}, or \cite{Ste3}.\hfill$\Box$

\begin{cor} \label{clocseh403}
Let $\Omega \subsetneqq \Ri^d$ be open, $\theta \in [0,\frac{\pi}{2})$,
$\mu,M > 0$, $\nu \in (0,1)$, $p \in (1,\infty)$, 
$k,l \in \{ 1,\ldots,d \} $ and $\chi \in C_{\rm b}^\infty(\Ri^d)$
with $\chi \neq 0$ and $d(\supp \chi, \Omega^{\rm c}) > 0$.
Then there exists a $c > 0$ such that 
\[
\|\partial_k \, M_\chi \, (I + A)^{-1} \, M_\chi \, \partial_l\|_{p \to p} \leq c
\]
for all $C \in \cs^\nu(\Omega,\theta,\mu,M)$.
\end{cor}
\proof\
Note that in $L_2$ the operator $\partial_k \, M_\chi \, (I + A)^{-1/2}$
is bounded by Theorem~\ref{tlocseh406}.
So by duality the theorem follows for $p = 2$.
Then the rest of the proof is similar to the proof of Corollary~\ref{clocseh402}.\hfill$\Box$

\begin{cor} \label{clocseh404}
Let $\Omega \subsetneqq \Ri^d$ be open, $\theta \in [0,\frac{\pi}{2})$,
$\mu,M > 0$, $p \in (1,\infty)$, 
$m \in \{ 1,\ldots,d \} $ and $\chi \in C_{\rm b}^\infty(\Ri^d)$
with $\chi \neq 0$ and $d(\supp \chi, \Omega^{\rm c}) > 0$.
Then there exists a $c > 0$ such that 
\[
\|\partial_m \, M_\chi \, (I + A)^{-1/2}\|_{p \to p} \leq c
\]
for all $C \in \cs^1(\Omega,\theta,\mu,M,{\rm real})$.
\end{cor}
\proof\
It suffices to prove the corollary with $M_\chi$ replaced by $M_\chi^2$.
We use again a commutator.
It follows from Lemma~\ref{llocseh308} that on $L_2(\Ri^d)$ one has
\begin{eqnarray*}
\lefteqn{
\partial_m \, M_\chi^2 \, (I + A)^{-1/2} - \partial_m \, M_\chi \, (I + A)^{-1/2} \, M_\chi
} \hspace{10mm} \\*
& = & \frac{1}{\sqrt{\pi}} \, 
   \int_0^\infty t^{-1/2} \, e^{-t} \, \partial_m \, M_\chi \, [M_\chi, S_t] \, dt  \\
& = & - \frac{1}{\sqrt{\pi}} \sum_{k,l=1}^d 
   \int_0^\infty\! \int_0^t t^{-1/2} \, e^{-t} \, \partial_m \, M_\chi \, 
    S_s \, M_{\partial_k \chi} \, \partial_l \, 
                M_{\tilde c_{kl}} \, M_{\widetilde \chi} \, S_{t-s} \, ds  \, dt
   + R
,  
\end{eqnarray*}
where $R$ is the contribution of the last term in (\ref{ellocseh308;1}) and 
$\tilde \chi \in C_{\rm b}^\infty(\Ri^d)$
is such that $d(\supp \tilde \chi, \Omega^{\rm c}) > 0$ and 
$\tilde \chi(x) = 1$ for all $x \in \supp \chi$.
Using the Gaussian bounds of Theorem~\ref{tlocseh309} it follows that there 
exists a $c > 0$ such that 
\begin{eqnarray*}
\|\partial_m \, M_\chi \, S_t \, M_{\partial_k \chi}\|_{p \to p} 
& \leq & c \, t^{-1/2} \, e^{t/2} \quad \mbox{and} \\
\|\partial_l \, M_{\tilde c_{kl}} \, M_{\widetilde \chi} \, S_t\|_{p \to p} 
& \leq & c \, t^{-1/2} \, e^{t/2} 
\end{eqnarray*}
for all $k,l \in \{ 1,\ldots,d \} $ and $t > 0$.
Then 
\begin{eqnarray*}
\lefteqn{
\Big\| \int_0^\infty\! \int_0^t t^{-1/2} \, e^{-t} \, \partial_m \, M_\chi \, 
    S_s \, M_{\partial_k \chi} \, \partial_l \, 
        M_{\tilde c_{kl}} \, M_{\widetilde \chi} \, S_{t-s} \, ds  \, dt \Big\|_{p \to p}
} \hspace*{60mm} \\*
& \leq & c^2 \int_0^\infty\! \int_0^t t^{-1/2} \, e^{-t/2} \, s^{-1/2} \, (t-s)^{-1/2} \, ds  \, dt  \\
& = & c^2 \, \pi \, \sqrt{2 \pi}
.  
\end{eqnarray*}
The contribution of $R$ can be estimated similarly and the current corollary
follows from Corollary~\ref{clocseh402}.\hfill$\Box$

\vertspace

We end with two propositions on second-order Riesz transforms.

\begin{prop} \label{plocseh406}
Let $\Omega \subsetneqq \Ri^d$ be open, $\theta \in [0,\frac{\pi}{2})$,
$\mu,M > 0$, $p \in (1,\infty)$, 
$m,n \in \{ 1,\ldots,d \} $ and $\chi \in C_{\rm b}^\infty(\Ri^d)$
with $\chi \neq 0$ and $d(\supp \chi, \Omega^{\rm c}) > 0$.
Then there exists a $c > 0$ such that 
$M_\chi \, (I + A)^{-1} L_p(\Ri^d) \subset W^{2,p}(\Ri^d)$
and 
\[
\|\partial_m \, \partial_n \, M_\chi \, (I + A)^{-1}\|_{p \to p} \leq c
\]
for all $C \in \cs^1(\Omega,\theta,\mu,M,{\rm real})$.
\end{prop}
\proof\
First let $C \in \ce\ch^1(\Omega,\theta,\mu,M,{\rm real})$.
Let $F = \supp \chi$ and define $\chi_1,\chi_2$ as in (\ref{etlocseh201;5})
and (\ref{etlocseh201;6}).
Let $M' = 2 \|\chi_2\|_{W^{1,\infty}(\Ri^d)} \, M + \|\chi\|_{W^{1,\infty}(\Ri^d)} + 1$
and  $C' = \chi_2 \, C + \chi_1 \, I$.
Then $C' \in \ce\ch^1(\Ri^d,\theta,\mu \wedge 1,M')$.
By \cite{ER19} Proposition~5.1 there exists a suitable $c > 0$ such that 
$(I + A')^{-1} L_p(\Ri^d) \subset W^{2,p}(\Ri^d)$ and 
\[
\|\partial_m \, \partial_n u\|_{p \to p} \leq c \, \|(I + A') u\|_p
\]
for all $u \in W^{2,p}(\Ri^d)$, where $A' = A_{C'}$.
Then 
\[
\|\partial_m \, \partial_n \, M_\chi u\|_{p \to p} 
\leq c \, \|(I + A') \, M_\chi u\|_p
= c \, \|(I + A) \, M_\chi u\|_p
.  \]
But 
\[
(I + A) \, M_\chi u
= M_\chi (I + A) u
   - \sum_{k,l=1}^d \Big( (\partial_l \chi) \, c_{kl} \, \partial_k u
                          + \partial_l \, c_{kl} \, (\partial_k \chi) u \Big)
.  \]
Hence it follows from Corollary~\ref{clocseh314} that there exists a suitable $c' > 0$
such that 
\[
\|\partial_m \, \partial_n \, M_\chi u\|_{p \to p} 
\leq c' \, \|(I + A) u\|_p
\]
for all $u \in W^{2,p}(\Ri^d) = D(A)$.
Then the proposition follows by approximation.\hfill$\Box$

\vertspace

On $L_2$ the same argument works for operators with complex coefficients.

\begin{prop} \label{plocseh407}
Let $\Omega \subsetneqq \Ri^d$ be open, $\theta \in [0,\frac{\pi}{2})$,
$\mu,M > 0$, $p \in (1,\infty)$, 
$m,n \in \{ 1,\ldots,d \} $ and $\chi \in C_{\rm b}^\infty(\Ri^d)$
with $\chi \neq 0$ and $d(\supp \chi, \Omega^{\rm c}) > 0$.
Then there exists a $c > 0$ such that 
$M_\chi \, (I + A)^{-1} L_2(\Ri^d) \subset W^{2,2}(\Ri^d)$
and 
\[
\|\partial_m \, \partial_n \, M_\chi \, (I + A)^{-1}\|_{2 \to 2} \leq c
\]
for all $C \in \cs^1(\Omega,\theta,\mu,M)$.
\end{prop}

\appendix

\section{Gaussian bounds} \label{AlocsehA}

The main aim of this appendix is to transfer weighted semigroup 
bounds into Gaussian kernel bounds, 
with optimal large time behaviour and optimal control 
of the constant in the Gaussian, including the H\"older bounds.
Throughout this appendix we write $C^\nu = C^{0,\nu}(\Ri^d)$ for all 
$\nu \in (0,1)$.

\begin{lemma} \label{llocsehA1}
Let $N,N^* \in \Ni_0$ and $\nu,\nu^* \in (0,1)$.
Let $S$ be a $C_0$-semigroup on $L_2(\Ri^d)$ and $T_1,T_2 \in \cl(L_2(\Ri^d))$
such that $T_1 \, S_t (C_c^\infty(\Ri^d)) \subset W^{N + \nu,\infty}(\Ri^d)$
and $T_2^* \, S_t^* (C_c^\infty(\Ri^d)) \subset W^{N^* + \nu^*,\infty}(\Ri^d)$
for all $t > 0$.
Assume that $[T_1,U_\rho] = [T_2,U_\rho] = 0$ for all $\rho \in \Ri$ and 
$\psi \in \cd_1$.
Let $a_0, a_1,\omega_1,\omega > 0$ and suppose that 
\[
\|U_\rho \, S_t \, U_{-\rho}\|_{2 \to 2} \leq a_0 \, e^{\omega \rho^2 t}
\]
for all $t > 0$, $\rho \in \Ri$ and $\psi \in \cd_1$.
Moreover, suppose that 
\begin{eqnarray*}
|||\partial^\alpha \, U_\rho \, T_1 \, S_t \, U_{-\rho} u|||_{C^\nu}
& \leq & a_1 \, t^{- d/4} \, t^{-|\alpha| / 2} \, t^{-\nu / 2} \, e^{\omega_1 (1 + \rho^2) t} \, \|u\|_2 ,  \\
|||\partial^\beta \, U_\rho \, T_2^* \, S_t^* \, U_{-\rho} u|||_{C^{\nu^*}}
& \leq & a_1 \, t^{- d/4} \, t^{-|\beta| / 2} \, t^{-\nu^* / 2} \, e^{\omega_1 (1 + \rho^2) t} \, \|u\|_2 
\end{eqnarray*}
for all $t > 0$, $\rho \in \Ri$, $\psi \in \cd_{|\alpha| \vee |\beta| \vee 1}$,
$u \in C_c^\infty(\Ri^d)$ and multi-indices $\alpha,\beta$ with 
$|\alpha| \leq N$ and $|\beta| \leq N^*$.
Let $c > 0$ be such that 
\begin{equation}
\sup \{ \psi(x) - \psi(y) : \psi \in \cd_{N \vee N^* \vee 1} \}
\geq c \, |x-y|
\label{ellocsehA1;3}
\end{equation}
for all $x,y \in \Ri^d$ and set $b = \frac{c^2}{4 \omega}$.

Then for all $t > 0$ the operator $T_1 \, S_t \, T_2$ has a continuous 
kernel $K_t$ which is $N_1$-times differentiable in the first variable and 
$N_2$-times in the second one, in any order.
Moreover, there exists an $a > 0$, depending only on 
$a_0$, $a_1$, $\omega$, $\omega_1$, $N$, $N^*$, $\nu$, $\nu^*$, $\|T_1\|$ and $\|T_2\|$
such that 
\begin{equation}
|(\partial_x^\alpha \, \partial_y^\beta K_t)(x,y)|
\leq a \, t^{- \frac{d + |\alpha| + |\beta|}{2}} \, 
    \Big( 1 + t + \frac{|x-y|^2}{t} \Big)^{\frac{d + |\alpha| + |\beta|}{2}} \, 
    e^{- b \frac{|x-y|^2}{t}}
\label{ellocsehA1;3.3}
\end{equation}
for all $x,y \in \Ri^d$, $t > 0$ and 
multi-indices $\alpha,\beta$ with $|\alpha| \leq N$ and 
$|\beta| \leq N^*$.

Finally, let $\gamma,\gamma^* \in (0,1)$, $\kappa > 0$, $\tau \in [0,1)$,
\[
0 < b_1 < b \, (1 - \tau)^2 \, \frac{c}{c + 2 \tau}
\]
and 
$\alpha,\beta$  multi-indices with $|\alpha| \leq N$,
$|\beta| \leq N^*$, $|\alpha| + \gamma \leq N + \nu$ and 
$|\beta| + \gamma^* \leq N^* + \nu^*$.
Then there exists an $a > 0$, depending only on 
$a_0$, $a_1$, $b_1$, $N$, $N^*$, $\nu$, $\nu^*$, $\gamma$, $\gamma^*$, $\kappa$, 
$\tau$, $\|T_1\|$ and $\|T_2\|$
such that 
\begin{eqnarray}
\lefteqn{
|(\partial_x^\alpha \, \partial_y^\beta K_t)(x+h,y+k) 
    - (\partial_x^\alpha \, \partial_y^\beta K_t)(x,y)|
} \hspace{9mm}   \label{ellocsehA1;30} \\*
& \leq & a \, t^{- \frac{d + |\alpha| + |\beta|}{2}} \, 
    (1 + t)^{\frac{d + |\alpha| + |\beta| + \gamma + \gamma^*}{2}} \, 
   \Bigg( \Big( \frac{|h|}{\sqrt{t} + |x-y|} \Big)^\gamma
          + \Big( \frac{|k|}{\sqrt{t} + |x-y|} \Big)^{\gamma^*} \Bigg)
    e^{- b_1 \frac{|x-y|^2}{t}}
\nonumber
\end{eqnarray}
for all $x,y,h,k \in \Ri^d$ and $t > 0$ with 
$|h| + |k| \leq \kappa \, \sqrt{t} + \tau \, |x-y|$.
\end{lemma}

In the proof of Lemma~\ref{llocsehA1} we need some estimates which are 
of independent interest.

\pagebreak[2]

\begin{lemma} \label{llocsehA2}
Let $\nu \in (0,1)$.
\begin{tabel}
\item \label{llocsehA2-1}
If $u \in C^\nu \cap L_2$ then $u \in L_\infty$ and 
\[
\|u\|_\infty
\leq \tfrac{d}{d+\nu} \, \varepsilon^\nu \, |||u|||_{C^\nu}
   + |B(1)|^{-1/2} \, \varepsilon^{-d/2} \, \|u\|_2
\]
for all $\varepsilon \in (0,1]$.
\item \label{llocsehA2-2}
If $k \in \{ 1,\ldots,d \} $ and $u \in W^{1+\nu,\infty}(\Ri^d)$
then 
\[
\|\partial_k u\|_\infty
\leq \tfrac{1}{1+\nu} \, \varepsilon^\nu \, |||\partial_k u|||_{C^\nu}
   + \varepsilon^{-(1-\nu)} \, |||u|||_{C^\nu}
\]
for all $\varepsilon \in (0,1]$.
\item \label{llocsehA2-3}
If $u \in W^{1,\infty}(\Ri^d)$, then 
\[
|||u|||_{C^\gamma}
\leq 2 \, \|\nabla u\|_\infty^\gamma \, \|u\|_\infty^{1-\gamma}
\]
for all $\gamma\in (0,1)$.
\item \label{llocsehA2-4}
If $0 < \gamma < \nu$, then 
\[
|||u|||_{C^\gamma} 
\leq 2 \, |||u|||_{C^\nu}^{\frac{\gamma}{\nu}} \, \|u\|_\infty^{1 - \frac{\gamma}{\nu}}
\]
for all $u \in C^\nu \cap L_\infty$.
\end{tabel}
\end{lemma}
\proof\
Let $x \in \Ri^d$ and $h \in B(\varepsilon)$.
Then $|u(x)| \leq |u(x) - u(x+h)| + |u(x+h)| \leq |h|^\nu \, |||u|||_{C^\nu} + |u(x+h)|$.
Integration over $h$ gives
\begin{eqnarray*}
\varepsilon^d \, |B(1)| \, |u(x)|
& \leq & |||u|||_{C^\nu} \int_{B(\varepsilon)} |h|^\nu \, dh
   + \int_{B(\varepsilon)} |u(x+h)| \, dh  \\
& \leq & |||u|||_{C^\nu} \, d \, |B(1)| \int_0^\varepsilon \rho^\nu \, \rho^{d-1} \, d\rho
   + |B(\varepsilon)|^{1/2} \, \|u\|_2  \\
& = & |||u|||_{C^\nu} \, \frac{d}{d+\nu} \, |B(1)| \, \varepsilon^{d+\nu}
   + |B(1)|^{1/2} \, \varepsilon^{d/2} \, \|u\|_2
\end{eqnarray*}
from which Statement~\ref{llocsehA2-1} follows.

For Statement~\ref{llocsehA2-2} note that 
\begin{eqnarray*}
u(x+\varepsilon \, e_k) - u(x)
& = & \int_0^\varepsilon (\partial_k u)(x + t \, e_k) \, dt  \\
& = & \varepsilon \, (\partial_k u)(x)  
   + \int_0^\varepsilon \Big( (\partial_k u)(x + t \, e_k) - (\partial_k u)(x) \Big)  \, dt
.
\end{eqnarray*}
So 
\begin{eqnarray*}
|(\partial_k u)(x)| 
& \leq & \frac{1}{\varepsilon} \Big| u(x+\varepsilon \, e_k) - u(x) \Big|
   + \frac{1}{\varepsilon} \int_0^\varepsilon \Big| (\partial_k u)(x + t \, e_k) - (\partial_k u)(x) \Big|  \, dt  \\
& \leq & \varepsilon^{-(1 - \nu)} |||u|||_{C^\nu} 
   + \frac{1}{\varepsilon} \int_0^\varepsilon t^\nu \, |||\partial_k u|||_{C^\nu} \, dt  \\
& = & \varepsilon^{-(1 - \nu)} |||u|||_{C^\nu} 
   + \frac{1}{1+\nu} \, \varepsilon^\nu \, |||\partial_k u|||_{C^\nu}
\end{eqnarray*}
for all $x \in \Ri^d$.

The proof of Statements~\ref{llocsehA2-3} and \ref{llocsehA2-4} is easy.\hfill$\Box$

\vertspace

\noindent
{\bf Proof of Lemma~\ref{llocsehA1}.\ }\
We follow arguments as in \cite{ER15}, \cite{ER21} and \cite{Ouh6}.
Set $N_0 = N \vee N^* \vee 1$.

\firststep\
Since $T_1$ and $T_2$ commute with $U_\rho$ one has estimates
\[
\|U_\rho \, T_1 \, S_t \, U_{-\rho}\|_{2 \to 2} 
\leq a_0 \, \|T_1\| \, e^{\omega \rho^2 t}
\] 
and similarly for $T_2$.
It then follows from the first two statements of 
Lemma~\ref{llocsehA2} with $\varepsilon = t^{1/2} \, e^{-t}$ that
\begin{eqnarray*}
\|\partial^\alpha \, U_\rho \, T_1 \, S_t \, U_{-\rho} u\|_\infty
& \leq & \tfrac{1}{2} \, a_2 \, t^{- d/4} \, t^{-|\alpha| / 2} \, e^{\omega_2 (1 + \rho^2) t} \, \|u\|_2 ,  \\
\|\partial^\beta \, U_\rho \, T_2^* \, S_t^* \, U_{-\rho} u\|_\infty
& \leq & \tfrac{1}{2} \, a_2 \, t^{- d/4} \, t^{-|\beta| / 2} \, e^{\omega_2 (1 + \rho^2) t} \, \|u\|_2 
\end{eqnarray*}
for all $t > 0$, $\rho \in \Ri$, $\psi \in \cd_{N_0}$,
$u \in C_c^\infty(\Ri^d)$ and multi-indices $\alpha,\beta$ with 
$|\alpha| \leq N$ and $|\beta| \leq N^*$,
where $a_2 = 4 \, a_1 + 2 |B(1)|^{-1/2} \, (\|T_1\| \vee \|T_2\|) a_0$
and $\omega_2 = (\omega_1 + 1-\nu) \vee (\omega + \frac{d}{2})$.
Then using the last two statements of 
Lemma~\ref{llocsehA2} one deduces that 
\begin{eqnarray*}
|||\partial^\alpha \, U_\rho \, T_1 \, S_t \, U_{-\rho} u|||_{C^\gamma}
& \leq & a_2 \, t^{- d/4} \, t^{-|\alpha| / 2} \, t^{-\gamma / 2} \, 
      e^{\omega_2 (1 + \rho^2) t} \, \|u\|_2 ,  \\
|||\partial^\beta \, U_\rho \, T_2^* \, S_t^* \, U_{-\rho} u|||_{C^\gamma}
& \leq & a_2 \, t^{- d/4} \, t^{-|\beta| / 2} \, t^{-\gamma^* / 2} \, 
      e^{\omega_2 (1 + \rho^2) t} \, \|u\|_2
\end{eqnarray*}
for all $t > 0$, $\rho \in \Ri$, $\psi \in \cd_{N_0}$,
$u \in C_c^\infty(\Ri^d)$, multi-indices $\alpha,\beta$ and $\gamma,\gamma^* \in (0,1)$ with 
$|\alpha| + \gamma \leq N + \nu$ and $|\beta| + \gamma^* \leq N^* + \nu^*$.

\nextstep\
Let $\gamma \in (0,1)$ and $\alpha$ be a multi-index with $|\alpha| + \gamma \leq N + \nu$.
Note that $U_\rho \, \partial_k = \partial_k \, U_\rho + \rho \, M_{\partial_k \psi} \, U_\rho$
and $|\rho| \, t^{(n+1)/2} \leq n! \, e^{(1 + \rho^2) t}$ for all $n \in \Ni_0$.
Hence it follows by induction to $|\alpha''|$ that 
\[
|||\partial^{\alpha'} \, U_\rho \, \partial^{\alpha''} \, T_1 \, S_t \, U_{-\rho} u|||_{C^\gamma}
\leq (1 + N! \, 2^N)^{|\alpha''|} \, a_2 \, 
    t^{- d/4} \, t^{-(|\alpha'| + |\alpha''|) / 2} \, t^{-\gamma / 2} \,
    e^{(\omega_2 + |\alpha''|) (1 + \rho^2) t} \, \|u\|_2 
\]
for all $t > 0$, $\rho \in \Ri$, $\psi \in \cd_{N_0}$, $u \in L_2$ and 
multi-indices $\alpha',\alpha''$ with $|\alpha'| + |\alpha''| \leq |\alpha|$.
In particular, 
\[
|||U_\rho \, \partial^\alpha \, T_1 \, S_t \, U_{-\rho} u|||_{C^\gamma}
\leq a_2 \, c_1 \, 
    t^{- d/4} \, t^{-|\alpha| / 2} \, t^{-\gamma / 2} \,
    e^{(\omega_2 + N) (1 + \rho^2) t} \, \|u\|_2 
,  \]
where $c_1 = (1 + N! \, 2^N)^N$.
Similarly, 
\[
\|U_\rho \, \partial^\alpha \, T_1 \, S_t \, U_{-\rho} u\|_\infty
\leq a_2 \, c_1 \, 
    t^{- d/4} \, t^{-|\alpha| / 2}  \,
    e^{(\omega_2 + N) (1 + \rho^2) t} \, \|u\|_2 
.  \]
If $h \in \Ri^d$ with $|h| \geq 1$ then 
\begin{eqnarray*}
\|(I - L(h)) U_\rho \, \partial^\alpha \, T_1 \, S_t \, U_{-\rho} u\|_\infty
& \leq & 2 a_2 \, c_1 \, 
    t^{- d/4} \, t^{-|\alpha| / 2}  \,
    e^{(\omega_2 + N) (1 + \rho^2) t} \, \|u\|_2 \\
& \leq & 2 a_2 \, c_1 \, 
    t^{- d/4} \, t^{-|\alpha| / 2}  \, \Big( \frac{|h|}{\sqrt{t}} \Big)^\gamma
    e^{(\omega_2 + N + 1) (1 + \rho^2) t} \, \|u\|_2
.  
\end{eqnarray*}
So 
\[
\|(I - L(h)) U_\rho \, \partial^\alpha \, T_1 \, S_t \, U_{-\rho} u\|_\infty
\leq  a_3 \, t^{- d/4} \, t^{-|\alpha| / 2}  \, \Big( \frac{|h|}{\sqrt{t}} \Big)^\gamma
    e^{\omega_3 (1 + \rho^2) t} \, \|u\|_2
\]
for all $t > 0$, $h \in \Ri^d$, $\rho \in \Ri$, $\psi \in \cd_{N_0}$ and $u \in L_2$, where 
$a_3 = 2 a_2 \, c_1$ and $\omega_3 = \omega_2 + N + 1$.

\nextstep\
Let $\gamma \in (0,1)$ and $\alpha$ be a multi-index with $|\alpha| + \gamma \leq N + \nu$.
Let $h \in \Ri^d$, $\psi \in \cd_{N_0}$, $\rho \in \Ri$ and $u \in L_2$.
Set $t_0 = \frac{1}{\omega_3 (1+\rho^2)}$.
If $t \in (0,t_0]$ then 
\[
\|(I - L(h)) U_\rho \, \partial^\alpha \, T_1 \, S_t \, U_{-\rho} u\|_\infty
\leq a_3 \, e \, 
    t^{- d/4} \, t^{-|\alpha| / 2}  \, \Big( \frac{|h|}{\sqrt{t}} \Big)^\gamma \, \|u\|_2 
\]
and if $t \in (t_0,\infty)$ then 
\begin{eqnarray*}
\lefteqn{
\|(I - L(h)) U_\rho \, \partial^\alpha \, T_1 \, S_t \, U_{-\rho} u\|_\infty
} \hspace{20mm} \\*
& = & \|(I - L(h)) U_\rho \, \partial^\alpha \, T_1 \, S_{t_0} \, U_{-\rho} \, U_\rho \, S_{t-t_0} \, U_{-\rho} u\|_\infty  \\
& \leq & a_3 \, e \, 
    t_0^{- d/4} \, t_0^{-|\alpha| / 2} \, \Big( \frac{|h|}{\sqrt{t_0}} \Big)^\gamma 
       \, \|U_\rho \, S_{t-t_0} \, U_{-\rho} u\|_2  \\
& \leq & a_0 \, a_3 \, e \, 
    \Big( \omega_3 (1+\rho^2) t \Big)^{\frac{d}{4} + \frac{|\alpha| + \gamma}{2}} \, 
    t^{- d/4} \, t^{-|\alpha| / 2} \, \Big( \frac{|h|}{\sqrt{t}} \Big)^\gamma \, 
    e^{\omega \rho^2 (t-t_0)} \|u\|_2
.
\end{eqnarray*}
Hence 
\[
\|(I - L(h)) U_\rho \, \partial^\alpha \, T_1 \, S_t \, U_{-\rho} u\|_\infty
\leq a_4 \, \Big( 1 + \omega_3 (1+\rho^2) t \Big)^{\frac{d}{4} + \frac{|\alpha| + \gamma}{2}} \,
   t^{- d/4} \, t^{-|\alpha| / 2} \, \Big( \frac{|h|}{\sqrt{t}} \Big)^\gamma \, 
    e^{\omega \rho^2 t} \|u\|_2
\]
for all $t > 0$, where $a_4 = a_3 (1 \vee a_0) e$.
Similarly, 
\[
\|U_\rho \, \partial^\alpha \, T_1 \, S_t \, U_{-\rho} u\|_\infty
\leq a_4 \, \Big( 1 + \omega_3 (1+\rho^2) t \Big)^{\frac{d}{4} + \frac{|\alpha|}{2}} \,
   t^{- d/4} \, t^{-|\alpha| / 2} \, 
    e^{\omega \rho^2 t} \|u\|_2
\]
for all $\psi \in \cd_{N_0}$, $\rho \in \Ri$, $t > 0$, $|\alpha| \leq N$ and $u \in L_2$.
Also similar bounds are valid with $T_2^*$ and $S_t^*$.

\nextstep\
Let $h,k \in \Ri^d$, $\psi \in \cd_{N_0}$, $\rho \in \Ri$, $t > 0$, 
$\gamma,\gamma^* \in (0,1)$ and $\alpha,\beta$ be multi-indices with 
$|\alpha| + \gamma \leq N + \nu$ and $|\beta| + \gamma^* \leq N^* + \nu^*$.
Then 
\begin{eqnarray}
\lefteqn{
\|L(h) \, U_\rho \, \partial^\alpha \, T_1 \, S_t \, T_2 \, \partial^\beta \, U_{-\rho} \, L(-k)
   - U_\rho \, \partial^\alpha \, T_1 \, S_t \, T_2 \, \partial^\beta \, U_{-\rho}\|_{1 \to \infty}
} \hspace{10mm} \nonumber  \\*
& \leq & \|(I - L(h)) U_\rho \, \partial^\alpha \, T_1 \, S_{t/2} \, U_{-\rho}\|_{2 \to \infty} \, 
   \|(I - L(k)) U_\rho \, \partial^\beta \, T_2^* \, S_{t/2}^* \, U_{-\rho}\|_{2 \to \infty}
\nonumber  \\*
& & \hspace{5mm} {}
  + \|(I - L(h)) U_\rho \, \partial^\alpha \, T_1 \, S_{t/2} \, U_{-\rho}\|_{2 \to \infty} \, 
   \|U_\rho \, \partial^\beta \, T_2^* \, S_{t/2}^* \, U_{-\rho}\|_{2 \to \infty}
\nonumber  \\*
& & \hspace{5mm} {}
  + \| U_\rho \, \partial^\alpha \, T_1 \, S_{t/2} \, U_{-\rho}\|_{2 \to \infty} \, 
   \|(I - L(k)) U_\rho \, \partial^\beta \, T_2^* \, S_{t/2}^* \, U_{-\rho}\|_{2 \to \infty} \nonumber   \\
& \leq & a_5 \, \Big( 1 + \omega_3 (1+\rho^2) t \Big)^{\frac{d + |\alpha| + |\beta| + \gamma + \gamma^*}{2}} \,
   t^{- \frac{d + |\alpha| + |\beta|}{2}} \,  e^{\omega \rho^2 t} \, \cdot
\label{ellocsehA1;1}  \\*
& & \hspace{45mm} {} \cdot
   \Bigg( \Big( \frac{|h|}{\sqrt{t}} \Big)^\gamma \Big( \frac{|k|}{\sqrt{t}} \Big)^{\gamma^*}
         + \Big( \frac{|h|}{\sqrt{t}} \Big)^\gamma
         + \Big( \frac{|k|}{\sqrt{t}} \Big)^{\gamma^*}
   \Bigg)
,  \nonumber  
\end{eqnarray}
where $a_5 = 2^{(d + |\alpha| + |\beta| + 2)/2} \, a_4^2$.
Similarly, 
\begin{equation}
\|U_\rho \, \partial^\alpha \, T_1 \, S_t \, T_2 \, \partial^\beta \, U_{-\rho}\|_{1 \to \infty}
\leq a_5 \, \Big( 1 + \omega_3 (1+\rho^2) t \Big)^{\frac{d + |\alpha| + |\beta| }{2}} \,
   t^{- \frac{d + |\alpha| + |\beta|}{2}} \,  e^{\omega \rho^2 t}
.  
\label{ellocsehA1;2} 
\end{equation}

\nextstep
Let $t > 0$ and $\alpha,\beta$ be multi-indices with 
$|\alpha| \leq N$ and $|\beta| \leq N^*$.
Choosing $\rho = 0$ it follows from (\ref{ellocsehA1;2}) and the Dunford--Pettis theorem
that the operator
$\partial^\alpha \, T_1 \, S_t \, T_2 \, \partial^\beta$ has a kernel 
$K^{(\alpha,\beta)}_t \in L_\infty(\Ri^d \times \Ri^d)$.
If $\widetilde L$ denotes the w$^*$-continuous left regular representation 
of $\Ri^d \times \Ri^d$ in $L_\infty(\Ri^d \times \Ri^d)$, then it follows from 
(\ref{ellocsehA1;1}) that 
\begin{eqnarray*}
\lefteqn{
\|(I - \widetilde L(h,k)) K^{(\alpha,\beta)}_t\|_\infty
} \hspace{20mm} \\*
& \leq & a_5 \, t^{- \frac{d + |\alpha| + |\beta|}{2}} \, 
    \Big( 1 + \omega_3 \, t \Big)^{\frac{d + |\alpha| + |\beta| + \nu + \nu^*}{2}}
   \Bigg( \Big( \frac{|h|}{\sqrt{t}} \Big)^\nu \Big( \frac{|k|}{\sqrt{t}} \Big)^{\nu^*}
         + \Big( \frac{|h|}{\sqrt{t}} \Big)^\nu
         + \Big( \frac{|k|}{\sqrt{t}} \Big)^{\nu^*}
   \Bigg)
\end{eqnarray*}
for all $(h,k) \in \Ri^d \times \Ri^d$.
So $\lim_{(h,k) \to (0,0)} \|(I - \widetilde L(h,k)) K^{(\alpha,\beta)}_t\|_\infty = 0$
and $K^{(\alpha,\beta)}_t$ is uniformly continuous on $\Ri^d \times \Ri^d$.

Define $K_t = K^{(\alpha,\beta)}_t$ if $|\alpha| = |\beta| = 0$.
Thus $K_t$ is the kernel of $T_1 \, S_t \, T_2$.

Let $|\alpha| \leq N$, $|\beta| \leq N^*$ and $t > 0$.
Then for all $u,v \in C_c^\infty(\Ri^d)$ one has
\begin{eqnarray*}
\lefteqn{
(-1)^{|\alpha| + |\beta|} \int_{\Ri^d} \int_{\Ri^d} 
    K_t(x,y) \, (\partial^\alpha u)(x) \, (\partial^\beta v)(y) \, dx \, dy    
} \hspace{45mm} \\*
& = & (-1)^{|\alpha| + |\beta|}
   (T_1 \, S_t \, T_2 \, \partial^\beta v, \partial^\alpha \, \overline u)
= (-1)^{|\beta|} (\partial^\alpha \, T_1 \, S_t \, T_2 \, \partial^\beta v, \overline u)  \\
& = & (-1)^{|\beta|} \int_{\Ri^d} u(x) \, 
    (\partial^\alpha \, T_1 \, S_t \, T_2 \, \partial^\beta v)(x) \, dx  \\
& = & (-1)^{|\beta|} \int_{\Ri^d} \int_{\Ri^d} K^{(\alpha,\beta)}_t(x,y) \, u(x) \, v(y) \, dx \, dy
.
\end{eqnarray*}
So by density
\[
(-1)^{|\alpha| + |\beta|} \int_{\Ri^d \times \Ri^d} \hspace*{-2pt}
    K_t(x,y) \, (\partial^\alpha_x \, \partial^\beta_y \, w)(x,y) \, d(x,y)    
= (-1)^{|\beta|} \int_{\Ri^d \times \Ri^d}  \hspace*{-2pt}
    K^{(\alpha,\beta)}_t(x,y) \, w(x,y) \, d(x,y) 
\]
for all $w \in C_c^\infty(\Ri^d \times \Ri^d)$ and the
$(-1)^{|\beta|} K^{(\alpha,\beta)}_t$ are the successive distributional
derivatives of~$K_t$.
Since the $K^{(\alpha,\beta)}_t$ are continuous one deduces from the
lemma of Du Bois--Reymond that $K_t$ is $N$ times differentiable in
the first variable, the derivatives are $N^*$-times differentiable
in the second variable and all derivatives are continuous.

\nextstep\
Let $|\alpha| \leq N$, $|\beta| \leq N^*$, $t > 0$ and $x,y \in \Ri^d$.
Then it follows from (\ref{ellocsehA1;2}) that 
\[
|K^{(\alpha,\beta)}_t(x,y)|
\leq a_5 \, \Big( 1 + \omega_3 (1+\rho^2) t \Big)^{\frac{d + |\alpha| + |\beta| }{2}} \,
   t^{- \frac{d + |\alpha| + |\beta|}{2}} \,  e^{\omega \rho^2 t} \, e^{-\rho(\psi(x) - \psi(y))}
\]
for all $\rho \geq 0$ and $\psi \in \cd_{N_0}$.
Minimizing over $\psi$ and using (\ref{ellocsehA1;3}) gives
\[
|K^{(\alpha,\beta)}_t(x,y)|
\leq a_5 \, \Big( 1 + \omega_3 (1+\rho^2) t \Big)^{\frac{d + |\alpha| + |\beta| }{2}} \,
   t^{- \frac{d + |\alpha| + |\beta|}{2}} \,  e^{\omega \rho^2 t} \, e^{-\rho c |x-y|}
\]
and with the choice $\rho = \frac{c \, |x-y|}{2 \omega t}$ one deduces that 
\begin{equation}
|K^{(\alpha,\beta)}_t(x,y)|
\leq a_5 \, \Big( 1 + \omega_3 (t + \frac{c^2 \, |x-y|^2}{4 \omega^2 \, t}) \Big)^{\frac{d + |\alpha| + |\beta| }{2}} \,
   t^{- \frac{d + |\alpha| + |\beta|}{2}} \,  e^{- b \frac{|x-y|^2}{t}}
.  
\label{llocsehA1;4}
\end{equation}
This proves the bounds (\ref{ellocsehA1;3.3}).

\nextstep\ \label{locsehAStep7}
Let $\alpha,\beta$ be multi-indices, $\gamma,\gamma^* \in (0,1)$
and suppose that $|\alpha| + |\gamma| \leq N + \nu$ and 
$|\beta| + \gamma^* \leq N^* + \nu^*$.
Let $\kappa > 0$ and $\tau \in [0,1)$.
There exists a $\tau_1 \in (\tau,1)$ such that 
\[
b_1 =  b \, (1 - \tau_1)^2 \, \frac{c}{c + 2 \tau}
.  \]
Set $\lambda = \frac{b}{b_1} \geq 1$.
Further, let $\varepsilon,\eta \in (0,1)$.
Let $x,y,h,k \in \Ri^d$, $t > 0$ and suppose that 
$|h| + |k| \leq \kappa \, \sqrt{t} + \tau \, |x-y|$.
Let $\rho \in \Ri$ and $\psi \in \cd_{N_0}$.

If $RHS = a_5 \ldots$ denotes the right hand side of (\ref{ellocsehA1;1}),
then it follows from (\ref{ellocsehA1;1}) that 
\[
\Big| e^{-\rho \psi(x)} K^{(\alpha,\beta)}_t(x,y) \, e^{\rho \psi(y)}
   - e^{-\rho \psi(x-h)} K^{(\alpha,\beta)}_t(x-h,y-k) \, e^{\rho \psi(y-k)} \Big|
\leq RHS
.  \]
So 
\[
\Big| K^{(\alpha,\beta)}_t(x,y)
   - e^{\rho (\psi(x) - \psi(x-h))} K^{(\alpha,\beta)}_t(x-h,y-k) \, e^{- \rho (\psi(y) - \psi(y-k))} \Big|
\leq RHS \cdot e^{\rho (\psi(x) - \psi(y))}
\]
and 
\begin{eqnarray*}
\lefteqn{
\Big| K^{(\alpha,\beta)}_t(x,y) - K^{(\alpha,\beta)}_t(x-h,y-k) \Big|
} \hspace{10mm}    \\*
& \leq & RHS \cdot e^{\rho (\psi(x) - \psi(y))}
   + \Big| 1 - e^{\rho (\psi(x) - \psi(x-h))} \, e^{-\rho (\psi(y) - \psi(y-k))} \Big|
     \, |K^{(\alpha,\beta)}_t(x-h,y-k)|  \\
& \leq & RHS \cdot e^{\rho (\psi(x) - \psi(y))}
   + |\rho| \, (|h| + |k|) \, e^{|\rho| (|h| + |k|)} \, |K^{(\alpha,\beta)}_t(x-h,y-k)|
.
\end{eqnarray*}
Suppose $\rho \geq 0$.
Optimizing over $\psi$ gives
\begin{eqnarray}
\lefteqn{
\Big| K^{(\alpha,\beta)}_t(x,y) - K^{(\alpha,\beta)}_t(x-h,y-k) \Big|
} \hspace{10mm}   \nonumber  \\*
& \leq & RHS \cdot e^{- c \rho |x-y|}
   + \rho \, (|h| + |k|) \, e^{\rho (|h| + |k|)}
     \, |K^{(\alpha,\beta)}_t(x-h,y-k)|  \nonumber  \\
& \leq & a_5 \, \Big( 1 + \omega_3 (1+\rho^2) t \Big)^E \,
   t^{- \widetilde E} \,  e^{\omega \rho^2 t} \, e^{- c \rho |x-y|}
   \Bigg( \Big( \frac{|h|}{\sqrt{t}} \Big)^\gamma \Big( \frac{|k|}{\sqrt{t}} \Big)^{\gamma^*}
         + \Big( \frac{|h|}{\sqrt{t}} \Big)^\gamma
         + \Big( \frac{|k|}{\sqrt{t}} \Big)^{\gamma^*}
   \Bigg) \nonumber  \\*
& & \hspace{10mm} {}
+ \rho \, (|h| + |k|) \, e^{\rho (|h| + |k|)} \, |K^{(\alpha,\beta)}_t(x-h,y-k)|
,
\label{ellocsehA1;6}
\end{eqnarray}
where for briefety we set
$E = \frac{d + |\alpha| + |\beta| + \gamma + \gamma^*}{2}$ and 
$\widetilde E = \frac{d + |\alpha| + |\beta|}{2}$.
Choose $\rho = \frac{c |x-y|}{2 \lambda \omega t}$.
Before we estimate both terms in (\ref{ellocsehA1;6}) we need one more estimate
to replace the denominator $\sqrt{t}$ by $\sqrt{t} + |x-y|$.

Since 
$\frac{|x-y|}{\sqrt{t}} 
\leq \frac{1}{\sqrt{\varepsilon}} \, e^{\varepsilon \frac{|x-y|^2}{t}}$
it follows that 
\[
t^{-1/2} (\sqrt{t} + |x-y|)
= 1 + \frac{|x-y|}{\sqrt{t}}
\leq \frac{2}{\sqrt{\varepsilon}} \, e^{\varepsilon \frac{|x-y|^2}{t}}
\]
and 
\[
t^{-1/2}
\leq \frac{2}{\sqrt{\varepsilon}} \, \frac{1}{\sqrt{t} + |x-y|} \, e^{\varepsilon \frac{|x-y|^2}{t}}
.  \]
Therefore 
\[
\Big( \frac{|h|}{\sqrt{t}} \Big)^\gamma
    + \Big( \frac{|k|}{\sqrt{t}} \Big)^{\gamma^*}
\leq \frac{2}{\sqrt{\varepsilon}} 
   \Bigg(  \Big( \frac{|h|}{\sqrt{t} + |x-y|} \Big)^\gamma
         + \Big( \frac{|k|}{\sqrt{t} + |x-y|} \Big)^{\gamma^*}
       \Bigg) e^{\varepsilon \frac{|x-y|^2}{t}}
\]
and 
\begin{eqnarray*}
\Big( \frac{|h|}{\sqrt{t}} \Big)^\gamma
     \Big( \frac{|k|}{\sqrt{t}} \Big)^{\gamma^*}
& \leq & \frac{4}{\varepsilon} \Big( \frac{|h|}{\sqrt{t} + |x-y|} \Big)^\gamma
   \Big( \frac{|k|}{\sqrt{t} + |x-y|} \Big)^{\gamma^*} \, e^{2 \varepsilon \frac{|x-y|^2}{t}}  \\
& \leq & \frac{4}{\varepsilon} \, (\kappa + \tau)^{\gamma^*} \, 
    \Big( \frac{|h|}{\sqrt{t} + |x-y|} \Big)^\gamma \, e^{2 \varepsilon \frac{|x-y|^2}{t}} 
.  
\end{eqnarray*}
We estimate both terms in (\ref{ellocsehA1;6}) separately. 

For the first term note that 
\[
\omega \, \rho^2 \, t - c \, \rho \, |x-y|
= - \frac{c^2 \, |x-y|^2}{4 \omega \, t} \, \frac{2\lambda - 1}{\lambda^2}
\leq - \frac{b}{\lambda} \, \frac{|x-y|^2}{t}
= - b_1 \, \frac{|x-y|^2}{t}
.  \]
Therefore the first term in (\ref{ellocsehA1;6}) can be estimated by
\[
\frac{6}{\varepsilon} \, a_5 \, (1 + \kappa)
   \Big( 1 + \omega_3 \, t + \frac{\omega_3 \, c^2 \, |x-y|^2}{4 \lambda^2 \, \omega^2 \, t } \Big)^E \,
   t^{-\widetilde E} \, 
         \Bigg( 
         \Big( \frac{|h|}{\sqrt{t} + |x-y|} \Big)^\gamma
         + \Big( \frac{|k|}{\sqrt{t} + |x-y|} \Big)^{\gamma^*}
       \Bigg)
     \makebox[0pt][l]{$\displaystyle e^{- (b_1 - 2 \varepsilon) \frac{|x-y|^2}{t} }$}
. \hspace*{10mm}  
\]
For the second term we use (\ref{llocsehA1;4}) to estimate
\[
|K^{(\alpha,\beta)}_t(x-h,y-k)|
\leq a_5 \, \Big( 1 + \omega_3 (t + \frac{c^2 \, |x-y-h+k|^2}{4 \omega^2 \, t}) \Big)^{\widetilde E} \,
   t^{- \widetilde E} \,  e^{- b \frac{|x-y-h+k|^2}{t}}
.  
\]
Clearly
\[
|x-y-h+k|^2
\leq 2 |x-y|^2 + 2 (|h| + |k|)^2
\leq 6 |x-y|^2 + 4 \kappa^2 \, t
.  \]
For the exponential set 
$\eta = \frac{\tau_1}{1 - \tau_1} > 0$ and $\delta = \frac{\tau_1^2}{\tau^2} - 1 > 0$.
Then 
\begin{eqnarray*}
- b \, \frac{|x-y - h + k|^2}{t}
& \leq & - b \, \frac{|x-y|^2}{(1 + \eta) t} + b \, \frac{|h-k|^2}{\eta \, t}  \\
& \leq & - b \, \frac{|x-y|^2}{(1 + \eta) t} 
    + \frac{b}{\eta \, t} \, 
     \Big( (1 + \delta) \, \tau^2 \, |x-y|^2 + (1 + \delta^{-1}) \, \kappa^2 \, t \Big)  \\
& = & - b \, (1 - \tau_1)^2 \, \frac{|x-y|^2}{t} 
   + b \, \kappa^2 \, \frac{\tau_1 - \tau_1^2}{\tau_1^2 - \tau^2}
.
\end{eqnarray*}
So 
\begin{eqnarray*}
\lefteqn{
|K^{(\alpha,\beta)}_t(x-h,y-k)|
}  \hspace{10mm} \\*
& \leq & a_5 \, \Big( 1 + \frac{\omega_3 \, c^2 \, \kappa^2}{\omega^2} 
                      + \omega_3 \, t
                      + \frac{3 \omega_3 \, c^2 \, |x-y|^2}{2 \omega^2 t} \Big)^{\widetilde E} \,
   t^{- \widetilde E} \,  e^{- b (1 - \tau_1)^2 \frac{|x-y|^2}{t}}
    \, \exp( b \, \kappa^2 \, \frac{\tau_1 - \tau_1^2}{\tau_1^2 - \tau^2} )
.  
\end{eqnarray*}
Next we estimate the factor $\rho \, (|h| + |k|) \, e^{\rho (|h| + |k|)}$.
One has 
\[
\rho (|h| + |k|)
\leq \frac{c \, |x-y|}{2 \lambda \, \omega \, t} (\kappa \, \sqrt{t} + \tau \, |x-y|)
\leq \frac{c \, (\tau + \eta \, \kappa)}{2 \lambda \, \omega} \, \frac{|x-y|^2}{t}
   + \frac{c \, \kappa}{2 \eta \, \lambda \, \omega}
.  \]
and alternatively
\begin{eqnarray*}
\rho (|h| + |k|)
& \leq & \frac{c \, |x-y|}{2 \lambda \, \omega \, \sqrt{t}} 
   \frac{|h| + |k|}{\sqrt{t}}  \\
& \leq & \frac{c \, |x-y|}{\sqrt{\varepsilon} \, \lambda \, \omega \, \sqrt{t}} 
   \frac{|h| + |k|}{\sqrt{t} + |x-y|} \, e^{\varepsilon \frac{|x-y|^2}{t}}  \\
& \leq & \frac{c \, (1 + \kappa)^2}{\varepsilon \, \lambda \, \omega} \, 
      e^{2 \varepsilon \frac{|x-y|^2}{t}}
       \Bigg( 
         \Big( \frac{|h|}{\sqrt{t} + |x-y|} \Big)^\gamma
         + \Big( \frac{|k|}{\sqrt{t} + |x-y|} \Big)^{\gamma^*}
       \Bigg)
.
\end{eqnarray*}
So 
\begin{eqnarray*}
\rho \, (|h| + |k|) \, e^{\rho (|h| + |k|)}
& \leq & \frac{c \, (1 + \kappa)^2}{\varepsilon \, \lambda \, \omega} \, 
      \exp \Big( \frac{c \, \kappa}{2 \eta \, \lambda \, \omega} \Big) \, 
      \exp \Big( (\frac{c \, (\tau + \eta \, \kappa)}{2 \lambda \, \omega} + 2 \varepsilon) 
                     \frac{|x-y|^2}{t} \Big) \cdot
\\* 
& & \hspace{40mm} {} \cdot
       \Bigg( 
         \Big( \frac{|h|}{\sqrt{t} + |x-y|} \Big)^\gamma
         + \Big( \frac{|k|}{\sqrt{t} + |x-y|} \Big)^{\gamma^*}
       \Bigg)
.  
\end{eqnarray*}
Using the identity $b \, (1 - \tau_1)^2 - \frac{c \tau}{2 \lambda \omega} = b_1$
one deduces that the second term in (\ref{ellocsehA1;6}) can be estimated by 
\begin{eqnarray*}
\lefteqn{
\frac{a_5 \, c \, (1 + \kappa)^2}{\varepsilon \, \lambda \, \omega} \,
\exp \Big( \frac{c \, \kappa}{2 \eta \, \lambda \, \omega} \Big) \, 
\exp( b \, \kappa^2 \, \frac{\tau_1 - \tau_1^2}{\tau_1^2 - \tau^2} ) \, 
\Big( 1 + \frac{\omega_3 \, c^2 \, \kappa^2}{\omega^2} 
                      + \omega_3 \, t
                      + \frac{3 \omega_3 \, c^2 \, |x-y|^2}{2 \omega^2 t} \Big)^{\widetilde E} \cdot
} \hspace{45mm} \\
& & \cdot t^{- \widetilde E} \,  
       \Bigg( 
         \Big( \frac{|h|}{\sqrt{t} + |x-y|} \Big)^\gamma
         + \Big( \frac{|k|}{\sqrt{t} + |x-y|} \Big)^{\gamma^*}
       \Bigg)
e^{- (b_1 - 2 \varepsilon - \frac{c \kappa \eta}{2 \lambda \omega}) \, \frac{|x-y|^2}{t}}
\end{eqnarray*}
Then (\ref{ellocsehA1;30}) follows.\hfill$\Box$

\subsection*{Acknowledgements} 
Most of this work was carried out whilst 
the first named author visited the University of Bordeaux~I and the second named
 author was visiting the  Department of Mathematics at the University of  Auckland. 
Both authors wish to thank the University of Bordeaux~I,
the University of Auckland and  the CNRS for financial support. 
 Part of this work is supported by the Marsden Fund Council from Government funding, 
administered by the Royal Society of New Zealand.

\newcommand{\etalchar}[1]{$^{#1}$}

\end{document}